\documentclass[11pt,twoside]{article}
\usepackage{amsmath}
\usepackage{amsthm}
\usepackage{amssymb}
\usepackage{amsfonts}
\makeatletter

\newtheorem{The}{Theorem}[section]
\newtheorem{Cor}[The]{Corollary}
\newtheorem{Pro}[The]{Proposition}

\newtheorem{Rem}[The]{Remark}
\newtheorem{Lem}[The]{Lemma}
\newtheorem{Def}[The]{Definition}
\newtheorem{Con}[The]{Condition}

\def\endproof{\relax\ifmmode\expandafter\endproofmath\else
\unskip\nobreak\hfil\penalty50\hskip.75em\hbox{}\nobreak\hfil\bull
{\parfillskip=0pt \finalhyphendemerits=0 \bigbreak}\fi}
\def\endproofmath$${\eqno\bull$$\bigbreak}
\def\bull{\vbox{\hrule\hbox{\vrule\kern3pt\vbox{\kern6pt}\kern3pt\vrule}\hrule}}

\def\cancel#1#2{\ooalign{$\hfil#1\mkern1mu/\hfil$\crcr$#1#2$}}
\def\Dirac{\mathpalette\cancel D}
\def\dirac{\mathpalette\cancel\partial}

\newcommand{\ba}{\begin{eqnarray}}
\newcommand{\na}{\end{eqnarray}}
\newcommand{\ban}{\begin{eqnarray*}}
\newcommand{\nan}{\end{eqnarray*}}
\newcommand{\spinc}{\mathrm{Spin}^c}

\newcommand{\sta}{\stackrel}

\newcommand{\scr}{\mathcal}
\newcommand{\C}{\mathbb{C}}

\newcommand{\R}{\mathbb{R}}
\newcommand{\Z}{\mathbb{Z}}
\newcommand{\N}{\mathbb{N}}

\renewcommand{\AA}{\mathbb{A}}

\newcommand{\A}{\scr{A}}
\newcommand{\G}{\scr{G}}
\newcommand{\B}{\scr{B}}
\newcommand{\M}{\scr{M}}

\newcommand{\la}{\langle}
\newcommand{\ra}{\rangle}

\newcommand{\nab}{\nabla}
\newcommand{\ind}{\mathfrak{i}}
\newcommand{\s}{\mathfrak{s}}

\newcommand{\p}{\mathfrak{p}}
\renewcommand{\t}{\mathfrak{t}}
\renewcommand{\d}{\partial}

\begin{document}

\title{Seiberg-Witten-Floer  Homology and Gluing Formulae \\
}
\author{Alan L. Carey,  Bai-Ling Wang }
\date{}
\maketitle

\abstract
This paper gives a detailed  construction of Seiberg-Witten-Floer homology
for a closed oriented 3-manifold with a non-torsion $\spinc$ structure.
Gluing formulae for certain 4-dimensional manifolds
splitting along an embedded 3-manifold are obtained.

\tableofcontents

\newpage
\section{Introduction and statement of results}

In this paper we give a detailed account of the construction of
Seiberg-Witten-Floer homology for a closed oriented 3-manifold $Y$ 
with a non-torsion $\spinc$ structure $\s$ and discuss some applications.
The basic idea is to think of the Seiberg-Witten equations 
as defining critical points for
the Chern-Simons-Dirac functional. This functional is 
 on the configuration space (which consists of
pairs of $U(1)$
connections and spinors modulo the gauge transformations) but 
may have degenerate critical points. These can be removed
by  a suitable perturbation. To construct the Seiberg-Witten-Floer
homology one proceeds to study
 the time independent gradient flow lines from one critical 
point to another. These flow lines or trajectories
may be thought of as solutions to the
Seiberg-Witten equations
in temporal gauge on the four manifold $Y\times \R$.
A suitable choice of perturbation of 
the  Chern-Simons-Dirac functional gives the desired properties
(such as smoothness and transversality)
for the moduli space of trajectory flowlines. 

We now summarise the contents of the paper.
In section 2, we review the basic properties of 
Seiberg-Witten monopoles on a closed, oriented 3-manifold $(Y, g, \s)$
with $\spinc$ structure $\s$ and metric $g$.
Section 3 gives
a complete account of the infinite dimensional Morse
theory for the Chern-Simons-Dirac functional whose critical points
satisfy the Seiberg-Witten equations on $(Y, g, \s)$. The 
downward gradient flowline
of the Chern-Simons-Dirac functional, thought as a solution to
the Seiberg-Witten equations on a cylinder $\R\times Y$, is endowed with 
an energy function (\ref{E=2C}), such that when the energy is sufficiently
small, the flowline stays in a small neighbourhood of some
critical point (see Lemma \ref{smallenergy:smalldistance}). This
enables us to
analyse the local behaviour of any gradient flowline near a critical point
in Proposition  \ref{decay:expo}, under the small energy condition.  

In order to achieve transversality for the moduli space of
flowline, we adopt the sort of perturbations 
suggested by Kronheimer \cite{Kro} with an additional 
condition \ref{condition}. 
We also show that, for generic perturbations, the moduli space of
flowlines with finite energy is smooth for each connected component
(cf. Theorem \ref{transversal:M}), which can be compactified to
a smooth manifold with corners (cf. Proposition \ref{compactification}).

The analysis developed in this paper can be applied to more
general cases but we only discuss one namely 
Seiberg-Witten-Floer homology theory 
for a closed oriented 3-manifold $Y$ with $b_1(Y)>0$ and 
a $\spinc$ structure $\s$ such that $c_1(\s)$ is non-torsion.
For any rational homology 3-sphere, there exists
an equivariant Seiberg-Witten Floer homology in \cite{MW}, which
is a topological invariant of the underlying manifold and $\spinc$
structure. For the
case of $(Y, \s)$ with $b_1(Y)>0$ and $c_1(\s)=  0$, 
Seiberg-Witten-Floer homology theory will be investigated
elsewhere \cite{MW2}. The results on Seiberg-Witten-Floer homology for a closed
 3-manifold $Y$ with $b_1(Y)>0$ and a non-torsion $\spinc$ structure 
in this paper are summarized in the following theorem, whose
proof will be given in section 3.

\begin{The}\label{Theorem:1}
For any closed oriented 3-manifold $Y$ with $b_1(Y)>0$ and 
a $\spinc$ structure $\s$ such that $c_1(\s)\neq 0$ in 
$H^2(Y, \Z)/Torsion$, let $d(\s)$ be the divisibility of
$c_1(\s)$ in $H^2(Y, \Z)/Torsion$, that is
\[
d(\s) = g.c.d.\{ <c_1(\s), \sigma> : \ for \ \sigma \in 
H_2(Y,\Z)\}.
\]
Then there exists a finitely generated
 Seiberg-Witten-Floer complex whose homology
$HF^{SW}_*(Y, \s)$ satisfies the following properties:
\begin{enumerate}
\item $HF^{SW}_*(Y, \s)$ is a topological invariant of $(Y, \s)$ and 
is a $\Z_{d(\s)}$-graded Abelian group. 
\item There is an action of
\[
\AA (Y) = Sym^*(H_0(Y, \Z))\otimes \Lambda ^*\bigl(H_1(Y, \Z)/Torsion\bigr)\]
on $HF^{SW}_*(Y, \s)$ with elements in
$H_0(Y, \Z)$ and $H_1(Y, \Z)/Torsion$ decreasing degree
in $HF^{SW}_*(Y, \s)$ by $2$ and $1$ respectively.
\item For $(-Y, -\s)$, where $-Y$ is $Y$ with the reversed orientation and
 $-\s$ is the induced $\spinc$ structure,
the corresponding Seiberg-Witten-Floer complex $C_*(-Y, -\s)$
is the dual complex of $C_*(Y, \s)$. 
There is a natural pairing
\[
\la \ , \ \ra : \qquad  HF^{SW}_*(Y, \s) \times HF^{SW}_{-*}(-Y, -\s) 
\longrightarrow \Z
\]
such that $<z.\Xi_1, \Xi_2> = <\Xi_1, z.\Xi_2>$
 for any $z\in \AA (Y)\cong \AA(-Y)$ and any cycles
$\Xi_1 \in HF^{SW}_*(Y, \s)$ and $\Xi_2\in HF^{SW}_{-*}(-Y, -\s)$
respectively. 
\item For any subgroup $K\subseteq Ker (c_1(\s))\subset H^1(Y, \Z)$, 
there is a variant of Seiberg-Witten-Floer homology denoted by
$HF^{SW}_{*, [K]}(Y, \s)$. $HF^{SW}_{*, [K]}(Y, \s)$ 
is a topological invariant and a
$\Z$-graded $\AA (Y)$ module.
For any element $[u]\in H^1(Y, \Z)/K$, there is an action
of $[u]$ on $HF^{SW}_{*, [K]}(Y, \s)$ which decreases 
degree by $<[u]\wedge c_1(\s), [Y]>.$
There exists a $\AA (Y)$-equivariant homomorphism:
\ba
\label{homo:1}
\pi:\qquad
HF^{SW}_{*, [K]}(Y, \s) \longrightarrow HF^{SW}_{*}(Y, \s).
\na
  There is also a  natural pairing
\ba\label{pair:K}
\la \ , \ \ra : \qquad  HF^{SW}_{*, [K]}(Y, \s) \times
HF^{SW}_{-*, [K]}(-Y, -\s)\longrightarrow \Z
\na
 satisfying $<z.\Xi_1, \Xi_2> = <\Xi_1, z.\Xi_2>$
 for any $z\in \AA (Y)\cong \AA(-Y)$ and any cycles
$\Xi_1 \in HF^{SW}_{*, [K]}(Y, \s) $ and
$\Xi_2 \in HF^{SW}_{-*, [K]}(-Y, -\s)$ respectively. 
If $K_1\subset K_2$
are two subgroups in $Ker (c_1(\s))$, there is a $\AA (Y)$-equivariant
homomorphism:
\[
HF^{SW}_{*, [K_1]}(Y, \s)  \longrightarrow HF^{SW}_{*, [K_2]}(Y, \s).
\]
Moreover, $HF^{SW}_{*, [Ker (c_1(\s))]}(Y, \s)$ satisfies the following
periodicity property:
\ba\label{periodicity}
HF^{SW}_{m, [Ker (c_1(\s))]}(Y, \s) \cong HF^{SW}_{m (mod \ d(\s))}
(Y, \s),
\na
for any $m\in \Z$. 
\end{enumerate}
\end{The}

Section 4 gives an application to the problem of associating invariants
to a 4-manifold with boundary and gluing formulae for
4-dimensional monopole invariants. The setup here is to start with
 a 4-manifold $(X_+, \s_+)$  with cylindrical end modelled on $(Y, \t)$:
that is, over the end $[-2, \infty)\times Y$, there is a fixed isomorphism
between the restrction of $\s_+$ and the pull-back $\spinc$ structure
of $\t$.  In addition, we assume $c_1(det( \t)) $ is non-torsion. 

Finite energy solutions to the Seiberg-Witten equations on $X_+$ define a  
moduli space with finite variations of the
 perturbed Chern-Simons-Dirac functional on the end. We associate
to this moduli space a boundary asymptotic value map
\[
\partial_\infty: \qquad \M_{X_+}(\s_+) \longrightarrow
\M_{Y, X_+}(\t, \eta)
\]
where $\M_{Y, X_+}(\t, \eta)$ is the quotient of solutions to
the perturbed Seiberg-Witten equations on $(Y, \t, \eta)$ by the
action of those gauge transformations which can be extended to $X_+$. In fact,
\[
\pi_+: \qquad \M_{Y, X_+}(\t, \eta) \longrightarrow \M_Y(\t, \eta)
\]
 is a covering map with
fiber an $H^1(Y, \Z)/Im (i_+^*)$-homogeneous space. Here 
$Im (i_+^*)$ is the range of the map 
$i^*_+: H^1(X_+, \Z) \to H^1(Y, \Z)$. 

Over the end $[-2, \infty) \times Y$,
if $(X_+, \s_+)$ is modelled on $(Y, \t)$ up to an isomorphism $u\in 
C^\infty (Y, U(1))$, then the corresponding 
 asymptotic value map is given by $[u]\circ \partial_\infty$
with $[u] \in H^1(Y, \Z)/Im (i_+^*)$ determined by the connected
component of $C^\infty (Y, U(1))$ which $u$ belongs to.

Fix an orientations on $\Lambda^{top} H^1(X_+(0), Y; \R)
\otimes \Lambda^{top} H^{2, +}(X_+(0), Y; \R)$. The structure of the
fiber for $\partial_\infty$ and its compactification
are discussed in Proposition \ref{smooth:cyl} and Proposition 
\ref{compact:cyl}. These propositions tell us that for each $\Gamma_\alpha
\in \M_{Y, X_+}(\t, \eta)$, there is a Baire set of self-dual 2-forms with
compact support such that  the fiber of $\partial_\infty$ (if non-empty), 
denoted
by $ \M_{X_+}(\s_+, \Gamma_\alpha)$, is an oriented, smooth manifold of
dimension $\ind_{X_+}(\Gamma_\alpha)$ given by (\ref{index:cyl}). It can be
compactified to a smooth manifold with corners. In particular, for any 
$d\geq 0$, the components of dimension $d$ in $\M_{X_+}(\s_+) $,
denoted by $\M_{X_+}^d(\s_+)$, can be described as
\ba
\M_{X_+}^d(\s_+) = \cup_{\alpha\in \M_Y(\t, \eta)}\Bigl( 
\bigcup_{\small
\Gamma_\alpha\in \pi_+^{-1}(\alpha),
\ind_{X_+}(\Gamma_\alpha) =d} \M_{X_+}(\s_+, \Gamma_\alpha)\Bigr).
\label{M^d:cyl}
\na
Here $\M_{X_+}(\s_+, \Gamma_\alpha)= \partial_\infty^{-1}(\Gamma_\alpha).$

We then define a relative Seiberg-Witten invariant of $X_+$.
These relative Seiberg-Witten invariants take values in the 
homology groups $HF^{SW}_{*, [Im (i_+^*)]}(Y, \t)$.
 We denote this invariant of $X_+$ by $SW_{X_+}(\s_+)$, which
is a linear functional
\[
SW_{X_+}(\s_+, \cdot ): \qquad \AA(X_+) \longrightarrow
 HF^{SW}_{*, [Im (i_+^*)]}(Y, \t),
\]
where $\AA(X_+) = Sym^*(H_0(X_+, \Z)) \otimes \Lambda^*(H_1(X_1, \Z)/Torsion)$
the free graded algebra generated by the class
of elements in $H_0(X_+, \Z)$ and 
$H_1(X_+, \Z)$ with degree 2 and 1 respectively. 

$SW_{X_+}(\s_+, z)$ with $deg (z) =d$
can be represented in terms of the Seiberg-Witten invariant for
the moduli space $\M_{X_+}^d(\s_+)$ in (\ref{M^d:cyl}), that is,
for each non-empty component $\M_{X_+}(\s_+, \Gamma_\alpha)$ with
$\ind_{X_+}(\Gamma_\alpha) =d$, there is a topological invariant
$SW_{X_+}(\s_+, z, \Gamma_\alpha)$ defined in subsection \ref{Relative},
with
\[
SW_{X_+}(\s_+, z) = \sum_{\alpha \in \M_Y (\t, \eta)} \sum_{\Gamma_\alpha
\in \pi_+^{-1}(\alpha)} SW_{X_+}(\s_+, z, \Gamma_\alpha) <\Gamma_\alpha>.
\]

Section 4 also describes results 
on closed smooth 4-manifolds $(X, \s)$ with $\spinc$ structure $\s$
and a smooth embedded separating 3-manifold $Y$ such that there is an embedding
$[-2, 2]\times Y$ in $X$ and 
$\t=\s|_Y$ is a non-torsion $\spinc$ structure on $Y$. 
Consider a 1-parameter family of metrics $\{g_R\}_{R>0}$ on
$X$ such that for each $X(R) = (X, g_R)$, there are an isometrically embedded
submanifold $([-R-2, R+2]\times Y, dt^2 + g_Y)$ and two 4-manifolds
$X_\pm (R)$ obtained by setting $X(R) = X_+(R)\cup  X_-(R)$. As $R\to \infty$, 
$X(R)$ has a geometric limit, two  4-manifolds 
with cylindrical ends, denoted by $X_\pm(\infty)$. Note that there are
induced $\spinc$ structures $\s_\pm$ on  $X_\pm(\infty)$ from the
$\spinc$ structure $\s$ on $X$. 

Denote by $i_\pm$ the boundary embedding maps of $Y$ in $X_\pm(0)$, and 
by
$Im (i_\pm^*)$ the ranges of the maps $H^1(X_\pm(0), \Z) \to H^1(Y, \Z)$. Then
we have the following commutative diagram relating various subgroups
$Im(i_+^*)\cap Im(i_-^*)$, $Im(i_\pm^*)$ and $Im(i_+^*)+Im(i_-^*)$ of
$Ker (c_1(\t))$ in $H^1(Y, Z)$:
\ba
\begin{array}{ccccc}
&& \displaystyle{\frac {H^1(Y, Z)}{Im(i_+^*)}}&& \\[-2mm]
&\nearrow&&\searrow&\\[-2mm]
 \displaystyle{\frac {H^1(Y, Z)}{Im(i_+^*)\cap Im(i_-^*)}}&&&& 
\displaystyle{\frac {H^1(Y, Z)}{Im(i_+^*)+Im(i_-^*)}}\\[-2mm]
&\searrow&&\nearrow&\\[-2mm]
&&\displaystyle{\frac {H^1(Y, Z)}{Im(i_-^*)}}&& \end{array}
\label{diag:subgroup}
\na 
which induces a commutative diagram of $\AA(Y)$-equivariant homomorphisms
between various variants of the Seiberg-Witten-Floer homologies for $(Y, \t)$:
\ba
\begin{array}{ccccc}
&& HF^{SW}_{*, [Im(i_+^*)]}(Y, \t)&& \\[-2mm]
&\nearrow&&\stackrel{\quad\pi_+}{\searrow}&\\[-2mm]
HF^{SW}_{*, [Im(i_+^*]\cap Im(i_-^*))}(Y, \t)&&&& 
HF^{SW}_{*, [Im(i_+^*)+Im(i_-^*)]}(Y, \t)\\[-2mm]
&\searrow&&\stackrel{\pi_-\quad}{\nearrow} &\\[-2mm]
&&HF^{SW}_{*, [Im(i_-^*)]}(Y, \t).&& \end{array}
\label{homo:diag}
\na
With these notations understood, the relative invariants for
$(X_\pm (\infty), \s_\pm)$ are linear functionals
\[
SW_{X_\pm}(\s_\pm, \cdot ) : \qquad \AA(X_\pm) \longrightarrow
HF^{SW}_{*, [Im(i_\pm^*)]}(\pm Y, \pm\t).
\]

Fix an orientation on the line $\Lambda^{top} H^1(X, \R)
\otimes \Lambda^{top} H^{2, +}(X, \R)$ which is induced
from tensoring the orientations on
$\Lambda^{top} H^1(X_+(0), Y; \R)
\otimes \Lambda^{top} H^{2, +}(X_+(0), Y; \R)$ and
$\Lambda^{top} H^1(X_-(0), Y; \R)
\otimes \Lambda^{top} H^{2, +}(X_-(0), Y; \R)$.
When $b^+_2 =1$, we also need to fix an  orientation on the line
$\Lambda^{top} H^{2, +}(X, \R)$.

The Seiberg-Witten invariant for a closed manifold $(X, \s)$
is a linear functional 
\[
SW_X(\s, \cdot): \AA(X) \mapsto \Z
\]
 as defined in \cite{Tau4} where 
$\AA(X) = Sym^*(H_0(X)) \otimes \Lambda^*(H_1(X)/Torsion)$.
Note that  there is indeterminacy in
the $\spinc$ structures on $X(R)$ obtained by gluing $\s_\pm$ along
$(Y, \t)$. The set of these $\spinc$ structures,
say $\spinc (X, \s_\pm)$, 
can be represented by
\[
\spinc (X, \s_\pm) = \{ \s_+\#_{[u]} \s_- | [u] \in 
\displaystyle{\frac {H^1(Y, Z)}{Im(i_+^*)+Im(i_-^*)}}
\}.
\]
Here $\s_+\#_{[u]} \s_-$ is the resulting
$\spinc$ structure on $X$ obtained
by gluing $\s_\pm$ using the gauge transformation
$u$ representing $[u]$.

Now our main gluing formulae along $(Y, \t)$ with $b_1(Y)>0$ and
$c_1(\t)\neq 0$   
are contained in the following theorem.

\vspace{2mm}

\begin{The}\label{Theorem:2}
Let $X$ be a closed 4-manifold with $b^+_2\geq 1$
which can be 
split as $X= X_+ \cup_Y X_-$, where $X_+$ and
$X_-$ are 4-manifolds with boundary $\partial X_+ = -\partial X_- =Y$.
Suppose that we have $\spinc$ structures $\s_+$ and $\s_-$ on $X_+$
and $X_-$ respectively which restrict to $\t$ on the boundary.
Assume that $b_1(Y)>0$ and $c_1(\t)$ is non-torsion then the Seiberg-Witten
invariants for $(X, \s_+\#_{[u]} \s_-)$ can be expressed as
\ba
SW_X (\s_+\#_{[u]} \s_-, z_+ z_-) = 
\la [u]\bigl(\pi_+(SW_{X_+} (\s_+, z_+))\bigr), \pi_-(SW_{X_-} (\s_-, z_-) )\ra,
\na
where $[u]$ acts on
$HF^{SW}_{*,[Im(i_+^*)+Im(i_-^*)]}(Y, \t)$, 
$z_\pm \in \AA(X_\pm)$, $\pi_\pm$ are the homomorphisms given in
(\ref{homo:diag}), and the pairing on the right hand side 
is the natural pairing on 
\[
HF^{SW}_{*,[Im(i_+^*)+Im(i_-^*)]}(Y, \t)  \times 
HF^{SW}_{_*, [Im(i_+^*)+Im(i_-^*)]}(-Y, -\t),
\] with the degrees in $HF^{SW}_{*, [Im(i_+^*)+Im(i_-^*)]}(-Y, -\t)$
shifted by 
\[ d_X(\s) =  \frac 14 (c_1(\s)^2 - (2 \chi (X) + 3\sigma (X) ))
= deg (z_1) + deg (z_2).\]
When $b_2^+=1$, the
Seiberg-Witten invariants $SW_X (\s, \cdot)$ are defined with a fixed
orientation on $H^{2, +}(X, \R)$ such that 
$c_1(\s)\cdot \omega^+ >0$ for an oriented
generator $\omega^+$ of $H^{2, +}(X, \R)$.
In particular, there is a set $\cal S$ of $\spinc$ structures
$\s \in \spinc (X, \s_\pm)$ with 
$d_X(\s)= deg (z_1) + deg (z_2)$ such that
\[
\sum_{\s\in {\cal S}}
 SW_X (\s, z_+ z_-) = \la\pi( SW_{X_+} (\s_+, z_+)), 
\pi(SW_{X_-} (\s_-, z_-))\ra. \]
Here $\pi(SW_{X_\pm}(\s_\pm, z_\pm))$ are elements in
$HF^{SW}_*(\pm Y, \pm \t)$ under the homomorphism (\ref{homo:1})
and  the pairing on the right hand side is the natural pairing on
$HF^{SW}_*(Y, \t) \times HF^{SW}_*(-Y, -\t)\to \Z$. 
\end{The}

This gluing theorem generalizes the gluing result of 
Morgan-Szabo-Taubes (Theorem 9.5 in \cite{MST}),
has already found  applications to the 
study of Seiberg-Witten-Floer homology as
in \cite{Munoz}. We like to thank Vicente Munoz for his 
help in providing examples to test our gluing formulae.

We would like to  acknowledge our gratitude to Matilde Marcolli, Tom Mrowka
and Cliff Taubes for many useful correspondences.  The authors also
thank the Max-Planck-Institut fur Mathematik in Bonn
for hospitality and support.  
This research was supported in part by the Australian Research 
Council.

\section{Seiberg-Witten monopoles on 3-manifolds: review}

Let $Y$ be a closed, oriented 3-manifold equipped  with a Riemannian metric
$g$. To introduce the Seiberg-Witten equations on $(Y, g)$, we need to
choose a $\spinc$ structure on $Y$. This
is a lift of the $SO(3)$-frame bundle on $Y$ to a $\spinc(3)$-bundle $P$.
Note that $\spinc(3) = U(2)$, the determinant homomorphism 
$\spinc(3) \to U(1)$ determines a principal $U(1)$-bundle, the corresponding
complex line bundle is called the determinant line bundle 
of the $\spinc$ structure $P$. 

The isomorphism classes of $\spinc$ structures are classified by
the first Chern class of the determinant line bundle. Given a
$\spinc$ structure $P$, there is an associated $\spinc$ bundle
\[
 W= P \times _{U(2)} \C^2
\]
which is complex vector bundle whose sections form a module for the
Clifford bundle on $(Y, g)$. As a vector space the
Clifford bundle can be identified with
the exterior bundle on $T^*Y$ but with a different
algebra structure: the Clifford relation 
\[
v_1 v_2 + v_2 v_1 = -2 g(v_1, v_2)
\]
holds for two sections $v_1$ and $ v_2$ of $T^*Y$. To end this discussion,
we give a more practical definition of a $\spinc$ structure on $(Y, g)$.

\begin{Def}\label{spinc:3d} A $\spinc$ structure $\s$
on an oriented, closed 3-manifold
$(Y, g)$ is a  pair $(W, \rho)$ consisting of a  $U(2)$ bundle $W$
and a  Clifford multiplication homomorphism 
$\rho: T^*(Y) \to  End (W)$ satisfying the Clifford relation:
$\rho(v_1) \rho (v_2) + \rho (v_2) \rho(v_1) = -2 g(v_1, v_2)$.
and satisfying $\rho (e_1) \rho(e_2) \rho(e_3) =Id_{W}$.
\end{Def}

Locally, let $\{e^1, e^2, e^3\}$ be an oriented orthonormal
co-frame for $T^*Y$, then the Clifford multiplication can be written
as follows 
\[
\rho (\sum_{i=1}^3 a_ie^i ) = \left( \begin{array}{lr} 
 a_1i & - a_2 + i a_3 \\
 a_2 + i a_3& - a_1i \end{array}
\right).
\]

The Clifford multiplication homomorphism 
$\rho$ can be employed to define 
an $\R$-bilinear form $\sigma: W\otimes W \to T^*Y \otimes i\R$, given by
\ba\begin{array}{lll}
\sigma(\psi, \psi) &=& 
-\rho^{-1} (\psi\otimes \psi^* - \displaystyle{\frac {|\psi|^2}{2}} Id)\\[2mm]
&=&\displaystyle{\frac i2} Im(\la e_i.\psi, \psi \ra )e^i \in \Omega^1(Y, i\R).
\end{array}
\label{sigma:bil}
\na
Here $\psi$ is a section of $W$ and 
$<, >$ is the Hermitian product on $W$. It is easy to check that
$\sigma(\cdot, \cdot)$ satisfies the following properties.

\begin{Lem} Use $.$ to denote Clifford
multiplication, $<, >$ is the Hermitian product on $W$ and
$g$ is the Riemannian metric on $Y$, extended linearly to
$T^*Y\otimes_\R \C$, Then
\begin{enumerate}
\item $\sigma (\psi, \psi) .\psi = - \frac 12 |\psi|^2 \psi$, and
$|\sigma (\psi, \psi)|^2 = \frac 14 |\psi|^4$.
\item For any imaginary-valued 1-form $\alpha$ on $Y$, 
$\la \alpha. \psi, \psi \ra =   2 g( \alpha, \sigma (\psi, \psi))$
and $\sigma (\alpha. \psi, \phi) + \sigma (\psi, \alpha.\phi)
= - (Re \langle \psi, \phi\rangle ) \alpha$.
\item $\sigma (\psi, \phi) =0$ if and only if on $Y\backslash\psi^{-1}(0)$,
$\phi = i r \psi$
for a real-valued function $r$  on $Y\backslash\psi^{-1}(0)$.
\item If $\psi$ is a nowhere vanishing section of  $W$, then
$W \cong \underline{ \C} \psi \oplus \psi^{\perp}$, where
$\underline{ \C} \psi$ is the trivial complex line bundle generated
by $\psi$ and $ \psi^{\perp}$ is the orthogonal complement of $\psi$
under the Hermitian product. Moreover, 
$\sigma (\psi, \cdot)$ defines a bundle isomorphism between
$ \underline{ \R}  \psi \oplus  \psi^{\perp}$ and
$T^*Y \otimes i\R$. This last claim implies that a $\spinc$ structure
on $(Y, g)$ can be identified with Turaev's Euler structure
\cite{Turaev} on $(Y, g)$.
\end{enumerate}
\label{basic:cal}
\end{Lem}
\begin{proof}
Claims (a) and (b) follow from direct calculations.
We only prove claim (d), as from (d), we know that the
on $Y\backslash\psi^{-1}(0)$, the kernel of the map $\sigma(\psi, \cdot)$
consists of sections of $\underline{i\R}\psi$, and this implies claim (c). 

For a nowhere vanishing section $\psi$ of $W$, we write
$\psi =|\psi|\tau_1$ where $\tau_1$ is a unit-length spinor.
Choose a local basis $\{\tau_1, \tau_2\}$ for the $\spinc$ bundle,
so that  Clifford multiplication in the local orthonormal
coframe $\{e^1, e^2, e^3\}$ for $T^*Y$ is given by
\[
\rho(e_1) = \left(\begin{array}{lr}
i &0 \\ 0 & -i \end{array}\right),
\rho(e_2)= \left(\begin{array}{lr}
0 &-1\\ 1 &0 \end{array}\right),
\rho(e_3)= \left(\begin{array}{lr}
0 &  i\\  i &0  \end{array}\right).
\]
where $\{e^1, e^2, e^3\}$ can be expressed as
\[
e^1 = -2i \sigma (\tau_1, \tau_1), e^2= 2i \sigma (\tau_1, i\tau_1),
e^3 = -2i \sigma (\tau_1, \tau_2).
\]
Any section of $\underline {\R\psi} \oplus \psi^\perp$ can be
written as $u\tau_1 + v\tau_2$ for a real-valued
function $u$ and a complex-valued function $v$, then
\[
\sigma (\psi, \phi) = \frac i2 |\psi| (u e^1 -Im (v) e^2 + Re (v)e^3).
\]
Hence, $\sigma (\psi,\cdot)$ defines a bundle isomorphism 
between $T^*Y\otimes i\R$ and $\underline {\R\psi} \oplus \psi^\perp$.

In \cite{Turaev}, Turaev defined an Euler structure to be a homology
class of a nowhere vanishing vector field on $Y$, where 
two vector fields are called homologous if they are homotopic away from
a point. Identifying $TY$ with $T^*Y$ using the Riemannian metric, for
any Euler structure on $(Y, g)$, there is
a unit-length section $v$ of $T^*Y$ representing that Euler structure,
which gives rise to the following decomposition $T^*Y = \underline{\R.v}
\oplus L$. Then the Hodge star operator can endow $L$ with
a complex structure given by
$J= *(v\wedge \cdot)$. The corresponding $\spinc$ structure
is given by $W = \underline{\C.v}
\oplus L$ and Clifford multiplication by $v$
is given by $\left(\begin{array}{lr}
i &0 \\ 0 & -i \end{array}\right)$. 
\end{proof}

With  the Levi-Civita connection $\nabla$ on the cotangent bundle
$T^*Y$, a $U(1)$-connection $A$ on the determinant bundle $det (W)$
determines a connection 
$\nabla_A$ on $W$ such that $\rho$ is parallel, that is,
for $v$ and $\psi$ sections of $T^*Y$ and $W$ respectively,
$\nabla_A$ satisfies 
\[
\nabla_A (v.\psi) = (\nabla v) .\psi + v. \nabla_A (\psi).
\]
Then $\nabla_A$ is called a $\spinc$ connection on $W$.
Applying the Clifford multiplication, we can define a Dirac operator
$\dirac _A = \rho\circ \nabla_A: \Gamma (W) \to \Gamma (W) $.

Now we can introduce the Seiberg-Witten equations on $(Y, g, \s)$
 for a pair $(A, \psi)$ consisting of a  $U(1)$ connection $A$ on the
determinant line bundle of $\s$ and a spinor section
$\psi$ of $W$:
\cite{KM1} \cite{CMWZ} \cite{MT} \cite{Lim} \cite{Chen}: \ba
\left\{ \begin{array}{l}
\dirac_A \psi =0, \\
*F_A = \sigma (\psi, \psi) + \eta
\end{array}\right.
\label{SW:3d}
\na
where $\eta$ is a co-closed, imaginary-valued 1-form on $Y$.

Denote by $ \A_Y(\s)$  the
Seiberg-Witten configuration space on $(Y, \s)$
consisting of pairs $(A, \psi)$ as above
Then $\A_Y(\s)$ is an affine space modelled on
$\Omega^1 (X, i\R) \times \Gamma(W)$. For analytic reasons, we shall use
the $L^2_1$-integrable configuration space, denoted $\A_{L^2_1}$,
and we require$\eta$ to be $L^2_2$-integrable, hence $\|\eta\|_{C^0}$ is
bounded.

The automorphism group of $det(W)$ is the gauge group
$\G_Y=\G_{L^2_2}$ of $L^2_2$-maps of $Y$ to
$U(1)$, locally modelled on the  Lie algebra $\Omega_{L^2_2}^0( Y, i\R)$.
The gauge group $\G_{L^2_2}$ acts on $\A_{L^2_1}$ by gauge transformations:
$(A, \psi) \sta{u}{\longmapsto} (A- 2u^{-1}du, u\psi)$, which 
leave the  Seiberg-Witten equations invariant. It
acts freely on $\{ (A, \psi) | \psi \ne 0 \}$, 
those elements are said to be irreducible
while  $(A, \psi)$ with $\psi = 0$ is said to be reducible. Denote
by $\B = \A_{L^2_1}/\G_{L^2_2}$ the quotient of $\A_{L^2_1}$ by  
gauge transformations. We use the notation
$\A_{L^2_1}^*$ for the open subset of irreducible configurations,
and $\B^* = \A_{L^2_1}^*/\G_{L^2_2}$ for the
quotient space of $\A_{L^2_1}^*$
by the gauge group $\G_{L^2_2}$. Then $\B^*$ is a Hausdorff space
\cite{FU}\cite{Morgan}, in fact,  $\B^*$ is a smooth
infinite dimensional  Hilbert manifold. 

At $(A, \psi) \in \A_{L^2_1}$, the tangent space can be identified
with the $L^2_1$-completion of $\Omega^{1}(Y, i\R) \oplus \Gamma (W)$, we
denote this completion by 
$\Omega^{1}_{L^2_1}(Y, i\R) \oplus \Gamma_{L^2_1}(W)$. The infinitesimal
action of $\G_{L^2_2}$ at $(A, \psi)$ defines a map
\ba\begin{array}{c}
G_{(A, \psi)}: \qquad \Omega^0_{L^2_2}(Y, i\R)
\rightarrow \Omega^{1}_{L^2_1}(Y, i\R) \oplus \Gamma_{L^2_1}(W)\\[2mm]
G_{(A, \psi)} (f) = (-2df, f\psi).
\end{array}
\label{gauge:action}
\na
Local slices exist for the action of $\G_{L^2_2}$ on $\A_{L^2_1}$.
At $[A, \psi] \in \B^*$,
the tangent space is given by the $L^2_1$-completion of
\[
\{ (\alpha, \phi) \in \Omega^{1}(Y, i\R) \oplus \Gamma (W) |
d^*\alpha + iIm <\psi, \phi> =0\}.
\]
The $L^2_1$ metric on $T_{[A, \psi]}\B^*$ is defined by
the natural $L^2_1$-norm on one forms and the $L^2_1$-norm
on sections of $W$ using the connection $\nabla_A$. This
metric is independent of the choice of representative $(A, \psi)$.
At a  reducible point $[A, 0] \in \B$, there is a neighbourhood of 
$[A, 0] \in \B$, which is homeomorphic to the quotient of
the $L^2_1$-completion of
\[
Ker d^* \oplus  \Gamma (W)
\]
by the action of $U(1)$ on the $\Gamma (W)$ factor. 
Hence, the reducible point $[A, 0]$ is a singular point in $\B$
with the link of singularity  homotopy equivalent to
$\C P^\infty$.

\begin{Def} The moduli space $\M_Y(\s, \eta)$ 
of the  Seiberg-Witten monopoles on $(Y, g, \s)$
is the solution space of  the Seiberg-Witten
equations (\ref{SW:3d}) modulo the gauge transformation group $\G_{L^2_2}$.
The moduli space $\M_Y^*(\s, \eta)$ is the irreducible part of 
$\M_Y(\s, \eta)$.  \end{Def}

There is an a priori pointwise estimate to the spinor part of any solution to
the Seiberg-Witten equations
(\ref{SW:3d}), which we now state as follows.

\begin{Lem}\label{estimate:psi}
Every solution to the  Seiberg-Witten equations
(\ref{SW:3d}) is gauge equivalent to a $C^\infty$ solution. Let $(A, \psi)$
be a smooth solution to the  Seiberg-Witten equations
(\ref{SW:3d}), then 
\[
|\psi|^2 \leq max_{y\in Y} \{ 0, -s (y) + 2 \|\eta\|_{C^0}\}
\]
where $s (y)$ is the scalar curvature for $(Y, g)$.
\end{Lem}
\begin{proof} The first claim follows from standard elliptic regularity
theory.  The estimate on $|\psi|^2$ follows from the Weitzenb\"ock formula
(Cf. \cite{KM1}\cite{Morgan}) 
\[
\dirac^2_A = \nabla_A^* \nabla_A + \displaystyle{\frac s4} - 
\displaystyle{\frac 1 2} \rho (*F_A).
\]
\end{proof}

Let $(A, \psi) \in \A_{L^2_1}$ be a
solution to the Seiberg-Witten equations
(\ref{SW:3d}), then there is a natural associated elliptic complex which
incorporates the linearization of the gauge action (\ref{gauge:action})
and the linearization of equations (\ref{SW:3d}) as follows.
The  linearization of the Seiberg-Witten equations (\ref{SW:3d})
at a solution $(A, \psi)$
defines the following map:
\ba
\begin{array}{c}
L_{(A, \psi)}: \Omega^1(Y, i\R)
\oplus \Gamma (W) {\rightarrow} \Omega^1(Y, i\R)
\oplus \Gamma (W)\\[2mm]
L_{(A, \psi)}(\alpha, \phi) = (*d\alpha -\sigma (\psi, \phi), \dirac_A\phi 
+ \frac{1}{2} \alpha. \psi).
\end{array}
\na

The infinitesimal action of $\G_Y=\G_{L^2_2}$ and 
linearization of the Seiberg-Witten 
equations  at  a solution $(A, \psi)$ define the
{\it deformation complex} 
\ba
\Omega^0_{L^2_2}(Y, i\R) \stackrel{G_{(A, \psi)}}{\to}
\Omega^1_{L^2_1}(Y, i\R) \oplus\Gamma_{ L^2_1}(W) \stackrel{L_{(A, \psi)}}{\to}
\Omega^1_{L^2}(Y, i\R) \oplus \Gamma_{ L^2}(W).
\label{3d:complex}
\na

It is convenient to put the deformation complex (\ref{3d:complex})
 in the following form:
\ba\begin{array}{c}
\Omega^1_{L^2_1}(Y, i\R) \oplus \Gamma_{L^2_1}(W) 
\oplus \Omega^0_{L^2_2}(Y, i\R) {\longrightarrow}
 \Omega^1_{L^2}(Y, i\R)
\oplus \Gamma_{L^2}(W) \oplus \Omega^0_{L^2}(Y, i\R)
\\
 {\cal T}_{(A, \psi)}= \left( \begin{array}{cc}
L_{(A, \psi)} & G_{(A, \psi)} \\ G^*_{(A, \psi)} 
& 0 \end{array}\right) \end{array}
\label{complex:3d1}
\na
where $G^*_{(A, \psi)}$ is the dual operator of $G_{(A, \psi)}$ under
the $L^2$-inner product: 
$$ G^*_{(A, \psi)} (a, \phi )= -2(d^*a + iIm \la \psi, \phi \ra).$$
Note that
${\cal T}_{(A, \psi)}$ is a compact perturbation, courtesy of Lemma
\ref{estimate:psi}, of the sum of the
signature operator and the twisted Dirac operator,
\[
\left( \begin{array}{cc}
*d & -2d \\ -2d^* &0 \end{array}\right) \oplus \dirac_{A} 
\]
which has index $0$. 

\begin{Def} 
We say that $[A, \psi]$ in  $\M_Y (\s, \eta)$  is  non-degenerate
if the middle cohomology of (\ref{3d:complex})
is zero:
\[
Ker L_{(A, \psi)} / Im G_{(A, \psi)}  =0.
\]
\end{Def}

We summarize the main properties of the monopole moduli space
 $\M_Y (\s, \eta)$,
their proofs can be found in \cite{Fro}  
 \cite{Lim} \cite{KM2} \cite{Wang}.

\begin{Pro}\label{3Dmonopole}
Let $(Y, g, \s)$ be a closed, Riemannian 3-manifold with
a $Spin^c$ structure $\s$, then the Seiberg-Witten moduli space 
$\M_Y (\s, \eta)$ has the following properties. 
There exists a Baire set of co-closed 1-forms $\eta
\in \Omega^1_{L^2_2}(Y, i\R)$ such that all the points in
$\M_Y (\s, \eta)$ are non-degenerate. Moreover, if $b_1(Y) >0$,
$\M_Y (\s, \eta)$ consists of only finitely many irreducible
points;  if $Y$ is a rational homology 3-sphere, assume
that  a generic $\eta$ satisfies $Ker \dirac_\theta =0$
(where $\theta$ is the unique reducible point in $\M_Y (\s, \eta)$,
that is, $*F_\theta = \eta$) , then
$\M_Y^* (\s, \eta)=  \M_Y(\s, \eta) - \{\theta\}$
 consists of only finitely many irreducible
points.
\end{Pro}

\begin{Rem} Proposition \ref{3Dmonopole} still holds for a
Baire set of co-closed imaginary-valued 1-forms on a
fixed non-empty open set $U\subset Y$ (if $b_1(Y)>0$ and $c_1(det(\s))=0$,
$U$ is chosen such that the
map: $H^2_{cpt}(U, \R)\to H^2(Y, \Z)$ is non-trivial). This is
established in  \cite{CMW}.
\end{Rem}

\section{ Seiberg-Witten-Floer homology}

In this Section we will discuss the construction of the Seiberg-Witten-Floer
homology groups. The main technical parts have been established
in \cite{MW}. Here we give another account of analytical parts of
the construction of Seiberg-Witten-Floer
homology groups and introduce
 a different  perturbation of the Chern-Simons-Dirac
functional to achieve transversality for the moduli space
of flowlines.

\subsection{The basics of the CSD functional: critical points, relative indices}

Let $Y$ be a closed, oriented 3-manifold with a Riemannian metric $g$ and
a $\spinc$ structure $\s = (W, \rho)$ on $Y$. Let $\A$ be
the $L^2_1$-configuration space for the Seiberg-Witten equations, the tangent
space at any point is the space of $L^2_1$-sections: 
$\Omega^1_{L^2_1}(Y, i\R) \oplus
L^2_1(W)$.
Let $\G=\G_Y$ be
 the $L^2_2$ gauge transformation group on $Y$. Its tangent space
at the identity is
$\Omega^0_{L^2_2}(Y, i\R)$ and $\G$ acts on $\A$ with the
quotient space denoted by ${\cal B}$. For simplicity, we assume
that $c_1(det (\s)) \neq 0$ in $H^2(Y, \R)$ if $b_1(Y)>0$.

Fix a $C^\infty$ connection $A_0$ on $det (W)$, 
then there is a  Chern-Simons-Dirac function on $\A$ 
as in \cite{KM1}:
\ba
{\cal C}_\eta (A, \psi) = - \frac 12 \int_Y (A-A_0) \wedge 
(F_A + F_{A_0} -*2\eta )
+  \int_Y \la \psi, \dirac_A \psi \ra dvol_Y.
\label{CSD}
\na
Here $\eta$ is a co-closed, imaginary-valued 1-form on $Y$
which is $L^2_2$-integrable and we can assume that $*\eta$ represents
a trivial cohomology class in $H^1(Y, i\R)$.
When $\eta=0$, we denote ${\cal C}_\eta $ by $\cal C$.

For $g \in \G$, we have
 \ba
{\cal C}_\eta (g. ( A, \psi)) - {\cal C}_\eta( A, \psi) =
 ( c_1(det(W)) \wedge [g])([Y]).
\label{gC-C}
\na
Here $[g]$ is the cohomology class
$[g]= [-2\pi ig^{-1}dg]$ on $Y$ which defines
a  homomorphism: $\G \to H^1(Y, \Z)$, the  kernel of this map is the
identity component of  $\G$.
Therefore, ${\cal C}_\eta$ descends to a map $\B \to \R/ d(\s)\Z$ where
$d(\s)$ is the divisibility of $c_1(\s)=c_1(det(W))$ in
$H^2(Y, \Z)/Torsion$:
\[
d(\s) = g.c.d.\{ <c_1(\s), \sigma> : \ for \ \sigma \in
H_2(Y,\Z)\}.
\]
In the case of $c_1(det(W)) =0$, we see
from (\ref{gC-C}) that ${\cal C}_\eta$ descends to a map $\B \to \R$.
If $c_1(det(W)) $ is not a torsion element, we know that
${\cal C}_\eta$ descends to a map $\B \to S^1$. 

With the $L^2$-metric on sections
of the spinor bundle $W$ given by twice the real part of the 
$L^2$-hermitian product,
the gradient of ${\cal C}_\eta$ at $(A, \psi)\in \A$ is given by 
$$\nab{\cal C}_\eta (A, \psi) = (*F_A - \sigma(\psi, \psi)- \eta, 
\dirac_A \psi),$$ 
which is $\G$-equivariant. Hence, $\nab{\cal C}_\eta$ defines
a section of the $L^2$-tangent bundle on $\B^*$.
Note that the tangent space of $\B^*$ at $[A, \psi]$ is the $L^2_1$-completion
of
\[
Ker G^*_{[A, \psi]} = \{ (a, \phi) \in \Omega^1(Y, i\R) \oplus
\Gamma (W) | d^* a + iIm \la \psi, \phi \ra =0  \}.
\]
The Hessian operator of
${\cal C}_\eta$ at $[A, \psi] \in \B^*$ is given by
\[\begin{array}{c}
Q_{[A, \psi]} : Ker G^*_{[A, \psi]} \longrightarrow Ker G^*_{[A, \psi]}\\[2mm]
Q_{[A, \psi]} : (a, \phi)  \mapsto
(*da - 2\sigma (\psi, \phi) - 2 df, \dirac_A \phi
+ \displaystyle{\frac 12} a. \psi + f \psi)
\end{array}
\]
where  $f$ is the unique solution in $\Omega^0_{L^2_2}(Y, i\R)$ of
the equation
\[
(d^*d +  \displaystyle{\frac 12} |\psi|^2 ) f = iIm \la \dirac_A\psi, \phi \ra.
\]
Notice that if $[A, \psi]$ is a critical point of ${\cal C}_\eta$, then $f =0 $ by
the maximum principle. The following Lemma describes the
properties of the Hessian operator  $Q_{[A, \psi]}$
of ${\cal C}_\eta$ at $[A, \psi] \in \B^*$.

\begin{Lem} $Q_{[A, \psi]}$ defines a closed, unbounded,
 essentially self-adjoint,
Fredholm operator on the $L^2$-completion of $Ker G^*_{[A, \psi]}$, and its
eigenvectors form an $L^2$-complete orthonormal basis for $Ker G^*_{[A, \psi]}$.The domain of $Q_{[A, \psi]}$ is the $L^2_1$-tangent space of $\B^*$. The
eigenvalues form a discrete subset of the real line which has no accumulation
points, and which is unbounded in both directions. Each eigenvalue has
finite multiplicity.
\label{spectrum:Hessian}
\end{Lem}
\begin{proof} We follow the arguments in \cite{Tau0}. First $Q_{[A, \psi]}$
defines a bounded, essentially self-adjoint,
Fredholm operator on the $L^2_1$-completion of $Ker G^*_{[A, \psi]}$ to the
$L^2$-completion of $Ker G^*_{[A, \psi]}$. Then
the forgetful map from $L^2_1$ to $L^2$ is compact and  the resolvent of
$Q_{[A, \psi]}$ is compact so that $Q_{[A, \psi]}$
has no essential spectrum (accumulation points or isolated eigenvalue
with infinite multiplicity). Therefore $Q_{[A, \psi]}$ has  discrete spectrum
without  accumulation points and each eigenvalue has
finite multiplicity. The remaining statements are standard in the
elliptic theory on compact manifolds.
\end{proof}

Let $[A, \psi]$ be a critical point of ${\cal C}_\eta$ on $\B$, then
$[A, \psi]$ is a gauge equivalence class of a solution $(A, \psi) \in \A$
to the Seiberg-Witten equations (\ref{SW:3d}). Hence, the
set of  critical points of ${\cal C}_\eta$ on $\B$
is $\M_Y(\s, \eta)$.  A critical point 
$[A, \psi]\in \M_Y^*(\s, \eta)$ is called non-degenerate 
if its Hessian operator is non-degenerate. It is easy to see that
a non-degenerate  critical point is isolated in $\M_Y(\s, \eta)$.
A reducible critical point is of the form
 $[A, 0]$ with $*F_A = \eta$. For generic
$\eta$, $*F_A = \eta$ admits solutions only when $Y$ is a rational
homology 3-sphere. We call a reducible critical point
non-degenerate if it is isolated in $\M_Y(\s, \eta)$. From  
Proposition \ref{3Dmonopole}, we obtain the following property
of  this critical point set.

\begin{Lem}\label{non-deg:C}
For generic $\eta$,
the critical points of ${\cal C}_\eta$ on $\B$ are all non-degenerate. 
\end{Lem}
\begin{proof}  
For an irreducible critical point $c\in \M_Y(\s, \eta)$, 
choose a $C^\infty$ solution $(A, \psi)$ to (\ref{SW:3d}) which represents $c$.
The Hessian operator $Q_c$ has a $\G$-equivariant extension given by
(\ref{complex:3d1}), namely
\[
\begin{array}{c}
 {\cal T}_{(A, \psi)}= \left( \begin{array}{cc}
L_{(A, \psi)} & G_{(A, \psi)} \\ G^*_{(A, \psi)}
& 0 \end{array}\right) \end{array}.
\] 
By direct calculation, we obtain 
\[
Ker Q_c \cong Ker {\cal T}_{(A, \psi)}
\cong Ker L_{(A, \psi)}/ Im G_{(A, \psi)}.
\]
Then Proposition \ref{3Dmonopole} implies that 
$Ker Q_c = 0$ for generic $\eta$.

If $\theta=[A, 0]$ is a critical point of ${\cal C}_{\eta}$
on $\B$, which only happens when $Y$ is a rational homology 
3-sphere ($*F_A = \eta$), then one can check that under the 
condition $Ker \dirac_A=0$, 
$Ker L_{(A, 0)}/Im G_{(A, 0)}$ is zero for the
deformation complex (\ref{complex:3d1}) at $(A, 0)$, hence 
$\theta$ is also non-degenerate.  
\end{proof}

Henceforth $\M_Y(\s, \eta)$ denotes the set of critical points
of ${\cal C}_\eta$ which are assumed all non-degenerate.

There are two fundamental lemmata concerning
the gradient of the Chern-Simons-Dirac functional and the estimate
of  local distance near a critical point. The first of these shows
that ${\cal C}_{\eta}$ satisfies the Palais-Smale condition
which was proved for the unperturbed case ($\eta=0$) in \cite{MST} \cite{JPW}.
The proof here follows the arguments from \cite{MST}.

\begin{Lem} \label{Palais:Smale}
For any $\epsilon > 0$, there  is $\lambda > 0$ such that if $[A,\psi]
\in {\cal B}$ has $L^2_1$-distance at least $\epsilon $ from all
the critical points in ${\cal M}_Y(\s , \eta)$, then
\[ \| \nabla {\cal C}_{\eta, \nu}({[A, \psi]}) \|_{L^2} > \lambda. \]
\end{Lem}
\begin{proof}
Suppose that, on the contrary, there is a sequence $c_i$
in $\B^*$ whose $L^2_1$-distance from $c_i$ to the set of
critical points is at least $\epsilon $, for which
$\|\nabla{\cal C}_\eta (c_i) \|_{L^2} \to 0$ as $i\to \infty$.
Choose a sequence of $(A_i, \psi_i)$ representing  $c_i$,
then there is a constant $C$ such that
\[
\int_{Y} (|*F_{A_i} -\sigma(\psi_i, \psi_i) -\eta|^2
+2 |\dirac_{A_i} \psi_{A_i}|^2)dvol \leq
C.
\]
Using the Weitzenb\"ock formula and Lemma \ref{basic:cal}, we
have
\[\displaystyle{
\|F_{A_i}\|^2_{L^2} +\|\sigma(\psi_i, \psi_i)+\eta\|^2_{L^2}
-2<*F_{A_i}, \eta> + 2\|\nabla_{A_i} \psi\|^2_{L^2}
+\frac s 2 \|\psi_i\|^2_{L^2} \leq C,}
\]
where $s$ is the scalar curvature on $(Y, g)$. Applying the inequality
$2|ab| \leq \frac 12 |a|^2 + 2 |b|^2$ and the Cauchy-Schwartz inequality,
we can rewrite the above expression as
\[\begin{array}{lll}
&&\displaystyle{\frac 12 \|F_{A_i}\|^2_{L^2} + \frac 18 \|\psi_i\|^4_{L^4}
+ 2\|\nabla_{A_i} \psi\|^2_{L^2}}\\[2mm]
&\leq&\displaystyle{ C + 3\|\eta\|^2_{L^2}+
max_{y\in Y} \{-\frac {s(y)}{2}\}\sqrt{Vol (Y)} \|\psi_i\|^2_{L^4}}\\[2mm]
&\leq& \displaystyle{ C + 3\|\eta\|^2_{L^2} + \frac {1}{16} \|\psi_i\|^4_{L^4} 
+ max_{y\in Y} \{-\frac {|s(y)|^2}{2}\} Vol (Y)}.
\end{array}
\]
{}From this inequality, we see that
$\|F_{A_i}\|_{L^2}$, $\|\psi_i\|_{L^4}$ and $\|\nabla_{A_i} \psi\|_{L^2}$
are all bounded independent of $i$.
By standard elliptic regularity theory we know
that there exists a subsequence of $\{A_i, \psi_i\}$,
such that $(A_i, \psi_i)$ converges in $L^2_1$-topology to
a limit $(A_\infty, \psi_\infty)$ which 
is a solution to the Seiberg-Witten equations (\ref{SW:3d}) and hence
represents a critical point in $\M_Y(\s, \eta)$. This
contradicts the fact that $\{[A_i, \psi_i]\}$ has distance at
least $\epsilon$ from the critical points.
\end{proof}

\begin{Rem}\label{condition:eta}
 Note that $\lambda$ is independent of $\eta$ when
$\eta$ is sufficiently small, hence from now on, we choose 
a sufficiently small $\epsilon_0 <<1$ and
$\eta$ to have $L^2$-norm less
than $\epsilon_0\lambda$.
\end{Rem}

The second lemma is the following $L^2_1$-distance estimate near
the non-degenerate critical points of ${\cal C}_{\eta}$, which 
was also established in \cite{MST}. 

\begin{Lem} \label{near:critical}
Suppose $\alpha$ is a  non-degenerate critical point of ${\cal C}_{\eta}$.
There exists a positive constant $C_\alpha$ and an $\epsilon>0$, depending
on $\alpha$,  such that if the
$L^2_1$-distance from $[A, \psi]$ to $\alpha$ is less than $\epsilon$,
then the $L^2_1$-distance from $[A, \psi]$ to $\alpha$ (denoted by
$dist_{L^2_1}( [A, \psi], \alpha )$) is bounded by
\[ C_\alpha \|\nab {\cal C}_{\eta, \nu} (A, \psi) \|_{L^2}. \]
\end{Lem}
\begin{proof}
When $dist_{L^2_1}( [A, \psi], \alpha )$ is sufficiently small, there
are representatives $(A, \psi)$ and $(A_\alpha, \psi_\alpha)$
for $[A, \psi]$ and $\alpha$ respectively such that
$(A, \psi) = (A_\alpha, \psi_\alpha) + (b, \phi)$ with
$(b, \phi) \in T_\alpha \B^*$ with $L^2_1$-norm less than some $\epsilon$ (to be
determined). Then
\[
\nabla {\cal C}_\eta (A, \psi) = Q_\alpha (b, \phi) + N (b, \phi),
\]
where $N (b, \phi)$ is a quadratic term in $(b, \phi)$,
satisfying $\|N (b, \phi)\|_{L^2} \leq C \epsilon \|(b, \phi)\|_{L^2_1}$
where $C$ is the universal constant for the Sobolev multiplication
theorem. As $Q_\alpha$ is non-degenerate, let $\lambda_1$ be the
first absolute eigenvalue of $Q_\alpha$
(Cf. Lemma \ref{spectrum:Hessian}).
We then obtain
\[
\begin{array}{lll}
\lambda_1 \|(b, \phi)\|_{L^2_1} &\leq& \|Q_\alpha(b, \phi)\|_{L^2}\\[2mm]
&\leq& \|\nabla {\cal C}_\eta (A, \psi)\|_{L^2}+ \|N (b, \phi)\|_{L^2}\\[2mm]
&\leq& \|\nabla {\cal C}_\eta (A, \psi)\|_{L^2}+
C \epsilon \|(b, \phi)\|_{L^2_1}.
\end{array}\]
Choose $\epsilon < \frac {\lambda_1}{2C}$, then
\[
\|(b, \phi)\|_{L^2_1} \leq \displaystyle{\frac {1}{\lambda_1 - C\epsilon}}
 \|\nabla {\cal C}_\eta (A, \psi)\|_{L^2}
\leq \displaystyle{\frac {2}{\lambda_1}}
 \|\nabla {\cal C}_\eta (A, \psi)\|_{L^2}.
\]
Set $C_\alpha = \displaystyle{\frac {2}{\lambda_1}}$, we get
$dist_{L^2_1} ([A, \psi], \alpha) \leq C_\alpha
 \|\nabla {\cal C}_\eta (A, \psi)\|_{L^2}.$
\end{proof}

For any two irreducible critical points $[A_0, \psi_0]$ and
$[A_1, \psi_1]$,  
choose a path $[A_t, \psi_t]$ on $\B^*$ connecting 
$[A_0, \psi_0]$ to $[A_1, \psi_1]$. Then by Lemma \ref{spectrum:Hessian},
 $Q_{[A_t, \psi_t]}$
defines a path of self-adjoint Fredholm operators, each of which has 
eigenvalues forming a discrete subset of the real line
with no accumulation points.
The spectral flow of $Q_{[A_t, \psi_t]}$ along this path is defined to
be the intersection number of the spectrum graph
$$
\bigcup_{t\in [0, 1]} Spec Q_{[A_t, \psi_t]}$$
with the line $\epsilon =0$ (Cf. \cite{APS}). 
The spectral flow of the Hessian operator 
defines a locally constant function on the space
of continuous paths in $\B^*$ between $[A_0, \psi_0]$ and
$[A_1, \psi_1]$ depending only on the homotopy class
of the path between $[A_0, \psi_0]$ and
$[A_1, \psi_1]$. By the Atiyah-Patodi-Singer
index theorem,  its $mod\ (d(\s))$ reduction doesn't depend
on the homotopy class of the path, where $d(\s)$
 is the divisibility of the first Chern class $c_1(det (W))$
 in $H^2(Y, \Z)/Torsion$.

We know that
$Q_{[A, \psi]}$ has a $\G$-equivariant extension 
${\cal T}_{(A, \psi)}$ referred to as 
the extended Hessian operator at a smooth solution $(A, \psi)$ to the
Seiberg-Witten equations (\ref{SW:3d})
on $\A$ representing $[A, \psi]$. 
Let $(A_0, \psi_0)$ and $(A_1, \psi_1)$ represent two irreducible
critical points $[A_0, \psi_0]$ and $[A_1, \psi_1]$
in $\M_Y(\s, \eta)$ respectively.
By its construction, the kernel of
${\cal T}_{(A_i, \psi_i)}$ is trivial, 
and so the spectral flow of $T_{(A, \psi)}$ from
$(A_0, \psi_0)$ to $(A_1, \psi_1)$ is a well defined $\Z$-valued function and
\ba
SF_{[A_0, \psi_0]}^{[A_1, \psi_1]} Q =
SF_{(A_0, \psi_0)}^{(A_1, \psi_1)} {\cal T} (mod\  d(\s)).
\label{spec:T}
\na

 For a homology 3-sphere $Y$, we need to choose the perturbation
$\eta$ so that the unique reducible critical point $\theta$
is isolated, that is, $(g, \eta)$ is away from the codimension one
set (Cf. \cite{MW}) where the corresponding Dirac operator
 has non-trivial kernel. Then   we can
replace $[A_0, \psi_0]$ by $\theta$ in (\ref{spec:T}). Now 
${\cal T}$ has one dimensional kernel at $\theta$, therefore
the spectral flow of ${\cal T}$ from $\theta$ to an
irreducible critical point is defined to be the intersection number of the
eigenvalues with the line $\lambda = -\epsilon$ where $\epsilon$ is any
sufficiently small positive number (see \cite{APS}).

Therefore, using the spectral flow of  the Hessian
operator,  we can define a relative
index map on the critical point set
$\M_Y (\s, \eta)$.

\begin{Def} \label{index:rel}
Using the spectral flow of  the Hessian
operator of the Chern-Simons-Dirac functional, there is a relative
index map
on the critical point set $\M_Y(\s, \eta)$:
\ba
\ind : \qquad \M_Y(\s, \eta) \times \M_Y(\s, \eta) \longrightarrow \Z_{d(\s)}
\na
which is the spectral flow ($mod \ d(\s)$)
of the extended Hessian operator along
any path connecting  two smooth solutions to (\ref{SW:3d}) representing
two critical points in $\B$.
When $c_1(\s)$ is a torsion class 
(as occurs for example when $Y$ is a homology 3-sphere)
then $\ind$ is an integer-valued function.
\end{Def}

\subsection{Estimates on gradient flowlines}

The downward gradient flow equation for ${\cal C}_{\eta}$ on $\A$ is
by definition given by the
paths of pairs  $(A(t), \psi(t))$ that satisfy the equations
\ba
\left\{ \begin{array}{l}
\displaystyle{\frac {dA(t)}{dt} }= -* F_{A(t)} + \sigma (\psi, \psi) + \eta ,
\\[2mm]
\displaystyle{\frac {d\psi (t)}{dt}} = - \dirac_{A(t)} \psi(t).
\end{array}
\right.
\label{grad:flow}
\na

Using the projection $\pi: Y\times \R \to Y$,
we identify
$ \Omega^{2, +}(Y\times \R, i\R) \cong  \pi^*(\Omega^1(Y, i\R)) $
by sending $\rho(t)$ to $\frac 12 (*\rho(t) + \rho(t)\wedge dt)$, and identify
$\Omega^1(Y\times \R, i\R) \cong\pi^*(\Omega^0(Y, i\R)\oplus \Omega^1(Y, i\R))$
by sending $(f(t) , \rho(t))$ to
 $\rho(t)+ f(t)dt$.
Using
 Clifford multiplication
by  $dt$, we can identify
the positive and negative spinor bundles $W^{\pm}$
 both as $\pi^*(W)$ via $\rho^{\pm}$,
such that for
any 1-form $a \in \Omega^1(Y, i\R)$  and $\psi \in W$, we have
$\rho^+(a.\psi) = a.dt.\rho^+(\psi) =a.\rho^-(\psi).$

Let $\psi(t)$ be a section of
$\pi^*(W)$ and $A(t)$ be a family of connections on $det(W)$.
Then we can view
$A(t)$ as a connection on $det(\pi^*W)$, and the twisted Dirac operator
$\Dirac_A$ for $(Y\times \R, \pi^*(W))$
can be expressed as $\d_t + \dirac_{A(t)}$ in the sense of
\[
\Dirac_A (\rho^+ \psi) = \rho^- \Bigl( (\d_t + \dirac_{A(t)} )\psi(t)\Bigr).
\]

Then, as discovered by Kronheimer-Mrowka \cite{KM1},
 the downward gradient flow equations (\ref{grad:flow})
are  the following Seiberg-Witten
equations on $(Y\times \R, \pi^*(W))$ in temporal gauge
(the $dt$ component of $U(1)$-connection on the determinant bundle vanishes
identically). 
\ba
\left\{
\begin{array}{l}
F^+_A = q (\psi, \psi) + \eta^+ \\[2mm]
\Dirac_A \psi = 0
\end{array}
\right.  \label{SW:4d}
\na
where $F_A^+$ is the self-dual part of the curvature $F_A$,
$\eta^+ = \frac 12( *_3 \eta + \eta \wedge dt)$,
and  $ q (\psi, \psi)$ is the
self-dual 2-form given by
\[
q (\psi, \psi) =(dt \wedge \sigma (\psi, \psi))^+ 
= \frac 12 (dt \wedge \sigma (\psi, \psi) + *_3 \sigma (\psi, \psi)).
\]
Two solutions to the equations (\ref{SW:4d}) are gauge equivalent
if the paths in $\A$ they
determine in temporal gauge are gauge equivalent under  the gauge group
 $\G_Y = Map(Y, U(1))$.

For any solution $ S(t)= (A(t), \psi (t) ) $ of the
downward gradient flow equations (\ref{grad:flow}), we have, for 
any $t_1 < t_2$,
\ba\begin{array}{lll}
&& {\cal C}_{\eta} (S(t_1)) -  {\cal C}_{\eta} (S(t_2)) \\[2mm]
&=& -\displaystyle{\int_{t_1}^{t_2} <\nabla {\cal C}_\eta ((S(t)), 
\frac{\partial S(t)}{\partial t} >_{L^2} dt}\\[2mm]
&=& \displaystyle{\int_{t_1}^{t_2}dt \int_Y
(|\frac {\partial A(t)}{\partial t}|^2 + 2|\frac {\partial \psi (t)}{\partial t}|^2) dvol_Y}
 \\[2mm]
&=&\displaystyle{\int_{t_1}^{t_2}\int_Y } (|*F_A -\sigma(\psi, \psi)-\eta|^2
+ 2|\dirac_A\psi|^2 )dvol_Y \\[2mm]
&=& \displaystyle{\int_{t_1}^{t_2}dt \int_Y
(|F_A|^2 + |\sigma (\psi, \psi)+\eta|^2 
-2 <*F_A, \eta> +2 |\nabla_A\psi|^2 +
\displaystyle{\frac s 2}|\psi|^2 ) dvol_Y}
\\[2mm]
&=&\displaystyle{\int_{t_1}^{t_2}dt \int_Y
(|F_A|^2 + |\sigma (\psi, \psi)+\eta|^2 + 
2(F_A\wedge \eta)+ 2|\nabla_A\psi|^2 +
\displaystyle{\frac s 2}|\psi|^2 ) dvol_Y} \end{array}
\label{3d:energy}
\na
where $F_A$ is the curvature of $A(t)$ on a $t$-slice of $[t_1, t_2]\times Y$,
$*$ is the Hodge star operator on $(Y, g)$, $\nabla_A$ is the
covariant derivative on $W$ over the $t$-slice of $[t_1, t_2]\times Y$
and all the norms are on the $t$-slice.   

Note that $A(t)$ can be viewed as a connection on
the determinant line bundle of the pull-back $\spinc$ structure
on $[t_1, t_2]\times Y$. We will write $\mathbb{A}=A(t)$ for this  connection on
the determinant line bundle over $[t_1, t_2]\times Y$ to distinguish it
from $A(t)$  viewed as a connection on the determinant line bundle 
over the $t$-slice of $[t_1, t_2]\times Y$. 
Similarly, denote the corresponding section of the  pull-back $\spinc$
bundle over $[t_1, t_2]\times Y$ by $\Psi = \psi(t)$. If $(A(t), \psi(t))$
is a solution to the gradient flow equation (\ref{grad:flow}),   then
$(\mathbb{A}, \Psi)$ is a solution to the 4-dimensional Seiberg-Witten
equations (\ref{SW:4d}) on $[t_1, t_2]\times Y$. With this
notation understood, the last integration in (\ref{3d:energy})
has a 4-dimensional analogue, called the energy of the
4-dimensional monopole. We denote it by $E_{[t_1, t_2]\times Y}(\mathbb A,
\Psi)$ and define it as follows.
\ba\label{energy:4d}
\begin{array}{lll}&&
E_{[t_1, t_2]\times Y}(\mathbb A,
\Psi)\\[2mm] &=& 
\displaystyle{\int_{[t_1, t_2]\times Y} }  (|F_{\mathbb A}|^2 +
2|q(\Psi, \Psi) + \eta^+|^2 + 2|\nabla_{\mathbb A}\Psi|^2 +
\displaystyle{\frac s 2}|\Psi|^2)dvol  + 2 (F_{\mathbb A}\wedge dt \wedge \eta) 
\end{array}
\na
Since 
$$|F_{\mathbb A}|^2 = |F_A|^2 + |\displaystyle{\frac
{\partial A(t)}{\partial t}}|^2, |\nabla_{\mathbb A}\Psi|^2
= |\nabla_A \psi|^2 +| \displaystyle{\frac
{\partial \psi (t)}{\partial t}}|^2$$
 and 
$$2|q(\Psi, \Psi) + \eta^+|^2
= |\sigma (\psi, \psi)+\eta|^2,$$ we know that
\ba \label{E=2C}
 E_{[t_1, t_2]\times Y}(\mathbb A,
\Psi) = 2 \bigl({\cal C}_{\eta} (S(t_1)) -  {\cal C}_{\eta} (S(t_2))\bigr).
\na
Hence we will say that any gradient flow line
$S(t)= (A(t), \psi (t) ) $ ($t\in [t_1, t_2]$) on $\A$ has
 a finite variation of ${\cal C}_{\eta}$ if $S(t)$ satisfies  
\[\displaystyle{
{\cal C}_{\eta} (S(t_1)) -  {\cal C}_{\eta} (S(t_2))
=  \int_{t_1}^{t_2} (\|\frac{\partial A(t)}{\partial t} \|^2_{L^2(Y)}
+ 2\|\frac{\partial \psi (t)}{\partial t} \|^2_{L^2(Y)}) dt < \infty.
}\]
 That is, the corresponding 4-dimensional monopole $({\mathbb A}, \Psi)$
 has finite energy over $[t_1, t_2]\times Y$, or equivalently 
\ba\displaystyle{
\int_{t_1}^{t_2} \| \nab {\cal C}_{\eta} (S(t))\|^2_{L^2(Y)}dt
= \int_{t_1}^{t_2} (\|\frac{\partial A(t)}{\partial t} \|^2_{L^2(Y)} 
+ \|\frac{\partial \psi (t)}{\partial t} \|^2_{L^2(Y)})dt < \infty.
}\label{finite:energy}
\na

\begin{Lem}
\label{first:estimate} Suppose $({\mathbb A}, \Psi) $ to be a
solution to the Seiberg-Witten equations (\ref{SW:4d})
on $[t-2, t+2]\times Y$ with finite energy
\[
E_{[t-2, t+2]\times Y}({\mathbb A}, \Psi) < E_0.
\]
For any positive integer $k$, there exists a constant $C$ depending on
the metric on $Y$ and the perturbation term $\eta$ such that
all the $C^k$-norms of $F_{\mathbb A} $ and $\Psi$ on
$[t-\frac 32, t+\frac 32] \times Y$ are bounded by $C (\sqrt{E_0} +1)$.
\end{Lem}
\begin{proof}
Let $N_s = [t-s, t+s]\times Y$. From $E_{N_2}({\mathbb A}, \Psi) < E_0$,
we immediately have
\[
\begin{array}{l}
\|F_{\mathbb A}\|^2_{L^2(N_2)} \leq \frac {E_0}{4} + C_0,\\[2mm]
\|\Psi\|^4_{L^4(N_2)} \leq 4 E_0 + C_0,\\[2mm]
\|\nabla_{\mathbb A}\Psi\|^2_{L^2(N_2)} \leq \frac {E_0}{4} + C_0,
\end{array}\]
for some constant $C_0$ depending on
the metric on $Y$ and $\eta$. The standard bootstrapping argument 
in elliptic regularity theory (Cf. \cite{Morgan}) 
implies that the $C^k$-norms of $F_{\mathbb A} $ and $\Psi$ on
$[t-\frac 32, t+\frac 32] \times Y$ are bounded by $C (\sqrt{E_0} +1)$
for some constant $C_0$ depending on
the metric on $Y$ and $\eta$.
\end{proof}

\begin{Lem}\label{smallenergy:smalldistance}
Given a sufficiently small 
$\epsilon >0$, there exists a constant $E_0>0$ such that, 
for a solution $(A(t), \psi(t))$ to the Seiberg-Witten equations 
(\ref{SW:4d}) on $[T-2, T+2]\times Y$ in temporal gauge, with finite energy 
$E_{[T-2, T+2]\times Y} (A(t), \psi(t)) < E_0$, there exists a critical
point $\alpha$ in $\M_{Y}(\s, \eta)$ with 
$L_1^2$-distance  from $[A(t), \psi(t)]\in \B$ to
$\alpha$ less than $\epsilon$ for any $t\in [T-1, T+1]$.
\end{Lem}
\begin{proof} The critical points of $\M_{Y}(\s, \eta)$ are non-degenerate
and hence isolated. Given an $\epsilon$ sufficiently less than 
$min\{dist_{L^2_1}(\alpha, \beta): \alpha, \beta \in \M_{Y}(\s, \eta)\}$,
let $\lambda$ be the constant in Lemma \ref{Palais:Smale} such that
if $[A, \psi]$ has $L_1^2$-distance at least $\epsilon$ to any
critical point in $\M_{Y}(\s, \eta)$, then 
$\|\nabla {\cal C}_\eta([A, \psi])\|_{L^2(Y)} >\lambda$.
Now we only need to show that for any $t\in [T-1, T+1]$, 
$\|\nabla {\cal C}_\eta([A(t), \psi(t)])\|_{L^2(Y)} < \lambda$.

Let $N_s = [T-s, T+s]\times Y$ and $(\dot{A}, \dot{\psi})= 
(\displaystyle{\frac {\partial A(t)}{\partial t}, 
\frac {\partial \psi(t)}{\partial t}})$. Differentiating the
Seiberg-Witten equations (\ref{SW:4d}) and noting that 
$(\displaystyle{\frac {\partial A(t)}{\partial t}, 
\frac {\partial \psi(t)}{\partial t}})= -\nabla {\cal C}_\eta (A, \psi)$
is $L^2$-orthogonal to the gauge orbit at each point, we obtain
the following elliptic system of equations on $[T-2, T+2]\times Y$,
\ba
\left\{ \begin{array}{l}
P_+(d \dot{A}) -2 q(\dot{\psi}, \psi)=0\\[2mm]
d^*\dot{A} + iIm <\dot{\psi}, \psi>=0\\[2mm]
\Dirac_{A} \dot{\psi} + \frac 12 \dot{A}.\psi =0.
\end{array}\right.
\na
We deduce from this system of equations and the estimates on $\psi$ in
Lemma \ref{first:estimate} that
\[
\begin{array}{lll}
\|(\dot{A}, \dot{\psi})\|_{L^2_1(N_{\frac 32})} &\leq &
C_0\bigl( \|(\dot{A}, \dot{\psi})\|_{L^2_1(N_{\frac 32})} +
\|(P_+(d \dot{A}), d^*\dot{A}, \Dirac_{A} \dot{\psi})\|_{L^2_1(N_2)}\bigr)
\\[2mm]
&\leq& C_0\bigl( \|(\dot{A}, \dot{\psi})\|_{L^2_1(N_{\frac 32})} +
\|(q(\dot{\psi}, \psi), iIm <\dot{\psi}, \psi>, \frac 12 \dot{A}.\psi)
\|_{L^2_1(N_2)}\bigr)
\\[2mm]
&\leq & C_1\bigl( \sqrt{E_0} + 
\|\psi\|_{C^0} \|(\dot{A}, \dot{\psi})\|_{L^2_1(N_2)} \bigr)
\\[2mm]
&\leq& C_1\bigl( \sqrt{E_0} + C_2 (\sqrt{E_0} +1) \sqrt{E_0}\bigr)\\[2mm]
&\leq& C_3\sqrt{E_0}.
\end{array}
\]
Here the $C_i$'s are constants depending on the metric on $Y$ and
$\eta$, and $\sqrt{E_0} $ bounds the energy of $(A(t), \psi(t))$. Therefore,
from the standard elliptic regularity theory of Degiorgi-Nash-Moser 
(\cite{Taylor} Page 169-179), we get that for any $t\in [T-1, T+1]$,
\[\|\nabla {\cal C}_\eta([A(t), \psi(t)])\|_{L^2(Y)}=
\|(\dot{A}, \dot{\psi})\|_{L^2(\{t\}\times Y)} \leq C_4 
\|(\dot{A}, \dot{\psi})\|_{L^2_1(N_{\frac 32})} \leq C_5 \sqrt{E_0}.
\]
Now we can choose our solution $(A(t), \psi(t))$
so that its energy is bounded by $E_0$ with $C_5 \sqrt{E_0} \leq \epsilon$,
then the claim in the lemma follows. 
\end{proof}

We remark that standard elliptic regularity theory implies
that $[A(t), \psi(t)]$ has, for all $k$, $C^k$-distance
$\epsilon$ or less to a critical point in $\M_Y(\s, \eta)$ for
any $t\in [T-1, T+1]$. The claims in Lemma \ref{smallenergy:smalldistance}
still hold for a solution to
the Seiberg-Witten equations on any compact set times $Y$ with sufficiently
small energy.

Now we can establish the crucial lemma concerning the 
exponential decay estimate on a
flowline with small energy.

\begin{Pro} There is a constant $\delta>0$, such that a choice
of sufficiently small $\epsilon $ determines $E_0>0$ with the following
property: suppose that $(A(t), \psi(t))$ is a solution to
the Seiberg-Witten equations on $[R_1-2, R_2+2]\times Y$ ($R_1 < R_2$)
in temporal gauge with energy
\[
E_{[R_1-2, R_2+2]\times Y} (A(t), \psi(t)) < E_0,
\]
then there exists a critical point $\alpha \in \M_Y(\s, \eta)$
such that the $L^2_1$-distance from $[A(t), \psi(t)]$ (for any
$t \in [R_1,  R_2]$) to $\alpha$ is at most 
$4\epsilon\bigl(
e^{-\delta (t-R_1)} + e^{-\delta (R_2 -t)}\bigr).
$
Furthermore, there exists a sequence of constants 
$\{\lambda_k\}_{k=0, 1, 2, \cdots}$ such that there is a
point $(A_\alpha, \psi_\alpha)$ representing $\alpha$ for which
\ba\label{C^k:estimate}
|\nabla^k (A(t)-A_\alpha)| + | (\nabla_{A_\alpha})^k (\psi(t) -\psi_\alpha)|
\leq \lambda_k \bigl(e^{-\delta (t-R_1)} + e^{-\delta (R_2 -t)}\bigr),
\na
for any point $(t, y) \in [R_1,  R_2]\times Y$.\label{decay:expo}
\end{Pro}
\begin{proof}
We can choose $\epsilon>0$ and $E_0>0$ as in Lemma
\ref{smallenergy:smalldistance}, so
that $dist_{L^2_1}([A(t), \psi(t)], \alpha) < \epsilon$ for some
critical point $\alpha \in \M_Y(\s, \eta)$ and any $t\in [R_1-1, R_2+1]$.
Then, on $[R_1-1, R_2+1]\times Y$, we can write 
\[
(A(t), \psi(t)) = (A_\alpha, \psi_\alpha) + (b, \phi)
\]
for some representative $(A_\alpha, \psi_\alpha)$ of
$\alpha$, and $(b, \phi)$ has $C^1$-norm less than $\epsilon$ 
and satisfies $d^*b +iIm<\psi_\alpha, \phi> =0$. 
As $(A(t), \psi(t))$ is a solution to the Seiberg-Witten equations
on $[R_1-1, R_2+1]\times Y$,
the pair $(b, \phi)$ must satisfy
\ba
\partial_t (b, \phi) + Q_\alpha (b, \phi) + N(b, \phi)=0,
\label{local:equation}
\na
where $Q_\alpha (b, \phi) = (*db -2\sigma (\psi_\alpha, \phi),
\dirac_{A_\alpha} \phi + \frac 12 b. \psi_\alpha)$ and 
$ N(b, \phi) = (-\sigma (\phi, \phi),  \frac 12 b. \phi).$

Note that $Q_\alpha$ is the Hessian operator of ${\cal C}_\eta$
at $\alpha$, whose $L^2$ spectrum is discrete, real and without
accumulation points (Cf. Lemma \ref{spectrum:Hessian}). The non-degeneracy
of $\alpha$ implies that $0$ is not in the spectrum of $Q_\alpha$. 
Let $\mu_\alpha$ be the distance between $0$ and the spectrum of $Q_\alpha$.
On the $t$-slice of $[R_1-1, R_2+1]\times Y$, 
\[
\|N(b, \phi)\|_{L^2_1(\{T\}\times Y)} < 
\epsilon \|(b, \phi)\|_{L^2_1(\{T\}\times Y)}.
\]

Let $\lambda_{\pm}$ be the components of $(b, \phi)$ in the 
positive and negative spectral subspaces
of $Q_\alpha$ respectively. 
Let $\|\lambda_{\pm}\|$ be the functions on $[R_1-1, R_2+1]$
given by taking the $L^2(Y)$-norm on the $t$-slice of $[R_1-1, R_2+1]\times Y$.
Then the differential equation (\ref{local:equation}) gives rise to
the following two differential inequalities
\[\begin{array}{c}
\partial_t \|\lambda_+\| - (\mu_\alpha -\epsilon) \|\lambda_+\| + \epsilon
\lambda_-\| \geq 0;\\[2mm]
\partial_t \|\lambda_-\| + (\mu_\alpha -\epsilon) \|\lambda_-\|- \epsilon
\lambda_+\| \leq 0.
\end{array}
\]
When $\epsilon < < \mu_\alpha$, the comparison principle (Cf. Lemma 9.4 in
\cite{Tau2}) can be invoked to establish exponential decay with
decay rate $\delta \leq \frac {\mu_\alpha-\epsilon}{2}$, that is,
for any $t\in [R_1-1, R_2+1]$, the $L^2$-norm of $(b, \phi)$
on the $t$-slice of $[R_1-1, R_2+1] \times Y$ is 
bounded by $4\epsilon \bigl(
e^{-\delta (t-R_1)} + e^{-\delta (R_2 -t)}\bigr)$. The standard bootstrapping
argument of elliptic regularity theory implies that
the $L^2_1$-distance from $[A(t), \psi(t)]$ for any
$t\in [R_1, R_2]$ to the critical point $\alpha$ is bounded
by $\epsilon \bigl(
e^{-\delta (t-R_1)} + e^{-\delta (R_2 -t)}\bigr)$, which is
the claim (\ref{C^k:estimate}).
\end{proof}

\begin{Cor} Let $(A(t), \psi(t))$ be a solution to the Seiberg-Witten
equations (\ref{SW:4d}) in temporal gauge with finite energy on
$\R \times Y$, then there
exist two 3-dimensional monopoles $(A_\alpha, \psi_\alpha)$ and
$(A_\beta, \psi_\beta)$ representing two critical points $\alpha$ and $\beta$
in $\M_Y (\s, \eta)$, such that
\begin{enumerate}
\item $A-A_\alpha$ and $\psi - \psi_\alpha$ decay exponentially along
  with their first derivatives as $t \to -\infty$,
\item $A-A_\beta$ and $\psi - \psi_\beta$ decay exponentially along with their
first derivatives as $t \to\infty$.
\end{enumerate}\label{limit:exist}
\end{Cor}
\begin{proof} To establish claim (a), we apply Lemma \ref{decay:expo} to
$(A(t), \psi(t))$ as a solution to the Seiberg-Witten
equations (\ref{SW:4d}) in temporal gauge on  
$(-\infty, -R] \times Y$ with energy arbitrarily small for
sufficiently large $R$. This gives the estimate 
(\ref{C^k:estimate}) for $(A_\alpha, \psi_\alpha)$ representing a
critical point $\alpha$ at any point $(t, y) \in (-\infty, -R] \times Y$
where, in (\ref{C^k:estimate})
we take $R_1$ to infinity and $R_2= -R$. Claim (b) can be proved
in the same way.
\end{proof}

{}From this corollary, we know that if 
$(A(t), \psi(t))$ is a solution to the Seiberg-Witten
equations (\ref{SW:4d}) in temporal gauge with finite energy on
$\R \times Y$, then there exists a sufficiently large $T$, such that
for $t\in (-\infty, -T]$ or $t\in [T, \infty)$, the gauge equivalent class
of $[A(t), \psi(t)]$ lies in a small $\epsilon$-neighborhood of
some critical point $\alpha$ or $\beta$ respectively, and approaches 
the critical point exponentially fast.

\subsection{Moduli spaces of flowlines: transversality, compactification
and orientibility.}

Let $\alpha, \beta $ be two critical points in $\M_Y(\s, \eta)$.
Denote by $\M(\alpha, \beta)$ the moduli space of
solutions to the gradient flow
equations (\ref{grad:flow}) with the asymptotic limits $\alpha, \beta $
as $t \to - \infty, + \infty$. Note that there is an $\R$-action on 
$\M(\alpha, \beta)$ given by translation in the $t$ variable, $\R$ acts freely
on $\M(\alpha, \beta)$ when $\alpha$ and $ \beta$ are
two different critical points. Denote by $\hat\M(\alpha, \beta)$
the quotient by this $\R$-action on $\M(\alpha, \beta)$.

Let $(A(t), \psi(t))$ represent a flowline in $\M(\alpha, \beta)$, then
any flowline in the same component as $[A(t), \psi(t)]$ has
constant energy, given by
\[
{\cal C}_\eta (A(-\infty), \psi(-\infty) ) -
{\cal C}_\eta (A(\infty), \psi(\infty) ),
\]
where $(A(-\infty), \psi(-\infty) )$ and $(A(\infty), \psi(\infty) )$
are two asymptotic limits of $(A(t), \psi(t))$ as $t\to \pm\infty$, whose
existence follows from Corollary \ref{limit:exist}. There is a constant
$\delta >0$ depending on
the critical points,  a constant $C$ depending only
on the local geometry of $Y$ and perturbation
$\eta$  such that for  sufficiently large $T$
\[
\sum_{0\leq k\leq 2} (|\nabla^k (A(t) -A_{\pm \infty})| +
 |(\nabla_{A_{\pm \infty}})^k (\psi(t) -\psi(\pm \infty)| ) 
\leq C e^{- \delta (|t|-T)}
\]
for $t \in (-\infty, T]$ and $t \in[T, \infty)$ respectively.

In this subsection, we will prove that generically 
each component of  the moduli space $\M(\alpha, \beta)$ is a
smooth manifold with dimension prescribed by the index theorem.
We remind the reader that we always identify
a path $(A(t), \psi(t))$ in $\A$ with a
pair consisting  of a $U(1)$-connection on the determinant bundle of
the pull-back $\spinc$ structure and a spinor on  $\pi^*(W)$ in
temporal gauge.  

To achieve transversality for the moduli space $\M(\alpha, \beta)$, 
we need to perturb the gradient flow equations
in a non-local fashion:
\ba
\left\{ \begin{array}{l}
\frac{\partial A }{\partial t} =- *F_A +\sigma (\psi, \psi )+\eta (A, \psi)
 \\[2mm]
\frac{\partial \psi }{\partial t} = -  \dirac _A (\psi) - \nu (A,\psi).\psi
+ \phi (A, \psi)
\end{array}\right.
\label{grad:pert}
\na
where $\eta (A, \psi), \nu (A,\psi)$ and $\phi (A, \psi)$
are functions of $(A, \psi)$.
This general perturbation is chosen so as to
 preserve translation invariance of
the gradient flow equation  under the action of $\R$.

We now describe in detail how to construct such perturbations.
For any Riemannian 3-manifold $(Y, g)$ with a $\spinc$
structure $\s$, we choose a complete $L^2$-basis $\{\nu_j\}_{j=1}^{\infty}$
for the imaginary-valued 1-forms on $Y$, and
a  complete $L^2$-basis $\{ \mu_j\}_{j=1}^{\infty}$ for the co-closed 
imaginary-valued 1-forms on $Y$. Under the Hodge decomposition,
\[
\Omega^1_{L^2}(Y, i\R) = H^1(Y, i\R) \oplus Im (d^*) \oplus Im (d),
\]
$\{ \mu_j\}_{j=1}^{\infty}$ span the space $H^1(Y, i\R) \oplus Im (d^*)$
and $\{\nu_j\}_{j=1}^{\infty}$ span the space $\Omega^1_{L^2}(Y, i\R) $.
We also need to choose a complete $L^2$-basis $\{\psi_j\}_{j=1}^{\infty}$
for the sections of the spinor bundle $W$. They are 
eigenvectors of the fixed Dirac operator $\dirac_{A_0}$ for a
a fixed $U(1)$-connection $A_0$ on the determinant bundle $det (\s)$. For each
co-closed 1-form $ \mu_j$, we associate a function on the
configuration space ${\cal A}$ 
by
\[
\tau_j (A, \psi) = \tau_j (A-A_0) = \int_Y (A-A_0) \wedge * \mu_j.
\]
 For simplicity, we assume that $\{ \mu_j\}_{j=1}^{b_1}$ consist of a
basis of $ H^1(Y, i\R)$, and $[*\mu_j] =0$ for $j>b_1$. It is easy to
see that:\\
(1) $\tau_j$ is invariant under gauge transformations for $j> b_1$;\\
(2) for the map obtained from the first $b_1$ functions
\[
(\tau_1,\cdots,  \tau_{b_1}): {\cal A} \to \R^{b_1},
\]
the gauge transformation $g: Y \to U(1)$ on ${\cal A}$ commutes
with the $H^1(Y, \Z) \cong \Z^{b_1}$-action on $\R^{b_1}$ by
translation
by the vector
 $$<([g^{-1}dg] \cup [*\mu_1], [g^{-1}dg] \cup [*\mu_2], \cdots,
[g^{-1}dg] \cup [*\mu_{b_1}]), [Y]>.$$

To each imaginary-valued 1-form $\nu_j$ we associate a function $\zeta_j $
on ${\cal A}$ as follows
\[
\zeta_j (A, \psi) = \zeta_j (\psi, \psi)=
\int_Y \la \nu_j. \psi, \psi \ra dvol_Y,
\]
where $\la, \ra$ is the Hermitian metric on the space of spinors. It is easy to
see that $\zeta_j$ is gauge invariant and real-valued. We also define
$\zeta_j (\psi, \phi) = 2\int_Y Re  \la \nu_j.\psi, \phi\ra dvol_Y.$

To each $\psi_j$, we  associate a function $\xi_j$ on ${\cal A}$ as follows
\[
\xi_j (A, \psi) = \int_Y \la e^{-\frac 12 Gd^* (A-A_0)}\psi_j, \psi\ra dvol_Y,
\]
where $G$ is the Green's operator for the ordinary Laplacian on
$L^2(Y)$. Let $\psi_j^A = e^{-\frac 12 Gd^* (A-A_0)\psi_j}$. Then
$\xi_j$ is invariant under the action of the subgroup
\[
H= \{e^{if}| f: Y\to \R, \int_Y f dvol_Y=0\}
\]
of the identity component of the gauge group $\G_Y$ 
where $\G_Y/H \cong S^1\times H^1(Y, \Z)$. The action 
of the $S^1$-factor has
weight $-1$, that is, for $u\in S^1$,
\[
\xi_j(u(A, \psi)) = u^{-1} \xi_j(A, \psi).
\]

 For any $n \geq b_1 (Y), k\geq 0$ and $l\geq 0$, $(\tau, \zeta, \xi)
= (\tau_1, \cdots, \tau_n, \zeta_1, \cdots, \zeta_k, \xi_1, \cdots, \xi_l)$
defines a map $\A \to \R^n \times \R^k \times \C^{l}$
with the following key property.

\begin{Lem} For any two points $(A_1, \psi_1)$ and
$(A_2, \psi_2)$ in $\A$ representing two different points in $\B$,
there exist $n$, $k$ and $l$ such that
\[
(\tau, \zeta, \xi) (A_1, \psi_1) \neq (\tau, \zeta, \xi) (A_2, \psi_2).
\]
\label{injectivity}
\end{Lem}
\begin{proof} On the contrary, suppose that for any $n$, $k$ and $l$,
$(\tau_1, \cdots, \tau_n, \zeta_1, \cdots, \zeta_k, \xi_1, \cdots, \xi_l)$
takes the same values on $(A_1, \psi_1)$ and
$(A_2, \psi_2)$. Then from the definition of
$\tau_j$, we know that $A_1 = A_2 -idf $ for some real-valued
function $f$ on $Y$, and from the definition of $\xi_j$,
we know that 
\[
e^{-\frac 12 G d^*(A_1-A_0)}\psi_1 = e^{-\frac 12 G d^*(A_2-A_0)}\psi_2,
\]
which implies that $\psi_1 = e^{\frac i2 f}\psi$. Hence
$(A_1, \psi_1)$ and
$(A_2, \psi_2)$ are gauge equivalent. This  
contradicts the assumption of the lemma. 
\end{proof}

As a corollary we have injectivity of the induced map
\[
(\tau, \zeta, \xi): \B \to \bigcup (\R^n\times \R^k\times \C^l)/(H^1(Y, \Z)
\times S^1)s,
\]
which implies that there are sufficiently many perturbations to achieve
the transversality of the moduli spaces of flowlines.
This will be crucial in the proof of Proposition \ref{transversal:M}.

Now we choose any function $\p \in C^{\infty} (\R^n \times \R^k \times \C^{l}, \R)$ (for
$n\ge b_1$, $k, l>0$)  where $\p$ is invariant under the action of $H^1(Y, Z)$
on $\R^{b_1}$ (the first $b_1$-components of $\R^N$), and
invariant under the action of $S^1\times H^1(Y, \Z)$ on $\C^{l}$. We can define
a function on ${\cal A}/{\cal G}$ by composition:
\[
\p( \tau_1, \cdots, \tau_n, \zeta_1, \cdots, \zeta_k, \xi_1, \cdots, \xi_l).
\]
Let ${\cal C}_\eta$ be the Chern-Simons-Dirac function as in (\ref{CSD}). 
We can
perturb this functional to form 
\[
{\cal C}_{\eta, \p} =
{\cal C}_\eta + \p( \tau_1, \cdots, \tau_n, \zeta_1, \cdots, \zeta_k,
\xi_1, \cdots, \xi_l).\]  

Let ${\cal P}$ be the subspace of 
$ \bigcup_{n\ge b_1, k, l> 0} C^{\infty} (\R^n \times\R^k\times \C^l, \R)
$
which is invariant under the above action of 
$H^1(Y, \Z)\times S^1\times H^1(Y, \Z)$ and has finite Floer $\epsilon$ norm.
To be specific, 
choose $\underline{\epsilon} = (\epsilon_k)_{k\in \N}$
to
be a given sequence of positive real numbers,
then the Floer $\epsilon$-norm (cf. 
\cite{Floer2}): is given by
\[
\| \p \|_{\underline{\epsilon}} 
=\sum_{k\ge 0} \epsilon_k sup |\nabla ^k \p|.
\]
By arguments similar to those in Lemma 5.1 of \cite{Floer2}, the sequence 
$\underline{\epsilon}$ can be chosen  
such that ${\cal P}$ is a Banach space and 
${\cal P}$ is dense in the subspace of  $L^2_k$, whose elements are invariant
under the above action of $H^1(Y, \Z)\times S^1\times H^1(Y, \Z)$.

For each $\p \in {\cal P}$, the gradient of ${\cal C}_{\eta, \p}$
 at $[A, \psi]$  with respect to the natural $L^2$ 
inner product is given by
\ba\begin{array}{c}
\Bigl( *F_A -\sigma (\psi, \psi)-\eta-
 \sum_{j=1}^{n} \frac{\d \p}{\d \tau_j}\mu_j
+ \sum_{j=1}^{l} \frac{\d \p}{\d \xi_j} d\circ G (iIm <\psi_j^A, \psi>),
\\[2mm]
\ \dirac_A (\psi)+ \sum_{j=1}^{K} \frac{\d \p}{\d \zeta_j}\nu_j.\psi 
+ \sum_{j=1}^{l} \frac{\d \p}{\d \xi_j} \psi_j^A \Bigr).
\end{array}\label{grad:perturbed}
\na
Here $\nabla \p = (- \sum_{j=1}^{n} \frac{\d \p}{\d \tau_j}\mu_j
+ \sum_{j=1}^{l} \frac{\d \p}{\d \xi_j} d\circ G (iIm <\psi_j^A, \psi>), 
\sum_{j=1}^{K} \frac{\d \p}{\d \zeta_j}\nu_j.\psi+
\sum_{j=1}^{l} \frac{\d \p}{\d \xi_j} \psi_j^A)$.

In order not to change the critical points of ${\cal C}_\eta$ when
perturbed by $\p$, we need to impose some condition on $\p \in {\cal P}$.
We remind the reader that, 
as in Lemma \ref{Palais:Smale}, given a small $\epsilon >0$,
there is a constant $\lambda$, such that 
\[
\|\nabla {\cal C}_\eta ([A, \psi])\|_{L^2} = 
\|\bigl(*F_A -\sigma (\psi, \psi)-\eta,
\dirac_A (\psi)\bigr)\|_{L^2} > \lambda
\]
 for
$[A, \psi]$ at $L^2_1$-distance
at least $\epsilon$ away from the set of critical points.
 For each critical point $\alpha \in \M_Y(\s, \eta)$,
let $B(\alpha, \epsilon)$ be ball of radius $\epsilon$
centered at $\alpha$ in the $L^2_1$ topology. 

\begin{Con}
We require that the perturbation function  $\p \in {\cal P}$ 
satisfies the following conditions:
\begin{enumerate}
\item  $\p =0$ on $\bigcup_{\alpha\in \M_Y(\s, \eta)} (\tau, \zeta, \tau) 
\bigl( B(\alpha, \epsilon)\bigr)$.
\item  The corresponding perturbation term $\nabla \p$ has
$L^2$-norm less than $\epsilon_0\lambda$ for a sufficiently small
$\epsilon_0 <<1$.
\end{enumerate}
Let ${\cal P}_c$ be the set of perturbation functions  $\p \in {\cal P}$
satisfying (a) and (b).
\label{condition}
\end{Con}

The corresponding downward gradient flow equation of ${\cal C}_\p$
is given by
\ba
\left\{ \begin{array}{lll}
\displaystyle{\frac{\d A}{\d t} }&=& -*F_A + \sigma (\psi, \psi) + \eta +
\sum_{j=1}^{N} \displaystyle{\frac{\d \p}{\d \tau_j}} \mu_j
- \sum_{j=1}^{l} \frac{\d \p}{\d \xi_j} d\circ G (iIm <\psi_j^A, \psi>) \\[2mm]
\displaystyle{\frac{\d \psi}{\d t}}& =& -\dirac_A \psi
- \sum_{j=1}^{K} \displaystyle{\frac{\d \p}{\d \zeta_j}}\nu_j.\psi -
\sum_{j=1}^{l} \frac{\d \p}{\d \xi_j} \psi_j^A
\end{array}\right.
\label{grad:pert2}
\na
for a  pair $(A(t), \Phi(t))$ regarded as a family of connections and spinors.
These are the perturbed  4-dimensional Seiberg-Witten equations  
on $Y \times \R$ in temporal gauge.
The corresponding perturbation of the Seiberg-Witten equations
on $Y \times \R$ is given by
\ba
\left\{
\begin{array}{l}
F^+_A = q (\psi, \psi ) + \eta^+ +
      \sum_{j=1}^{N} \displaystyle{\frac{\d \p}{\d \tau_j}} \mu_j^+
- \sum_{j=1}^{l} \frac{\d \p}{\d \xi_j} (d\circ G (iIm <\psi_j^A, \psi>))^+
\\[2mm]
\Dirac_A (\psi) + \sum_{j=1}^{K} \displaystyle{\frac{\d \p}{\d \zeta_j}}
      \nu_j.\psi+ \sum_{j=1}^{l} \frac{\d \p}{\d \xi_j} \psi_j^A= 0
\end{array}
\right.
\label{SW:4dpert}
\na
where $\mu_j^+ =\frac 12 ( *_3 \mu_j + \mu_j \wedge dt)$ and
$(d\circ G (iIm <\psi_j^A, \psi>))^+$ is the self-dual part of
$d\circ G (iIm <\psi_j^A, \psi>)\wedge dt$, and the  function 
$\p$ is the corresponding function on the 4-dimensional configuration
space of $U(1)$-connections and spinors on $(\R\times Y dt^2+ g)$ with
the pull back $\spinc$ structure. Note that $\p$ is invariant under gauge 
transformations. It is easy to see that $\p$ is a well-defined function
over the 4-dimensional configuration space.
For any solution $(A(t), \psi(t))$ to the perturbed Seiberg-Witten equations
(\ref{SW:4dpert}) on $[R_1, R_2]\times Y$ in temporal gauge, we say that
$(A(t), \psi(t))$ has finite energy if and only if
\[ 
\displaystyle{\int_{R_1}^{R_2} (\|\frac {\partial A(t)}{\partial t}\|_{L^2(Y)}^2
+ \|\frac {\partial \psi(t)}{\partial t}\|_{L^2(Y)}^2 )dt }< \infty.
\] 
 
\begin{Lem}\label{energy}
Let $ \gamma (t) = (A(t), \psi(t))$ be a solution to 
the perturbed Seiberg-Witten equations
(\ref{SW:4dpert}) on $[R_1, R_2]\times Y$ in temporal gauge, then 
$(A(t), \psi(t))$ has finite energy if and only if one of the
following terms is finite: 
\begin{enumerate}
\item ${\cal C}_\eta(A(R_1), \psi(R_1)) -{\cal C}_\eta(A(R_2), \psi(R_2))$, 
\item ${\cal C}_{\eta, \p}(A(R_1), \psi(R_1)) -
 {\cal C}_{\eta, \p}(A(R_2), \psi(R_2))$,
\item ${\cal C}(A(R_1), \psi(R_1)) -{\cal C}(A(R_2), \psi(R_2))$,
\item $E_{[R_1, R_2]\times Y}({\mathbb A}, \Psi)$ 
(as defined in (\ref{energy:4d})),
where $({\mathbb A}, \Psi) = (A(t), \psi(t))$ is the
corresponding 4-dimensional monopole on $[R_1, R_2]\times Y$.
\end{enumerate}
\end{Lem}
\begin{proof} Note that $\gamma (t)=(A(t), \psi(t))$ 
satisfies the perturbed Seiberg-Witten 
equations
(\ref{SW:4dpert}) on $[R_1, R_2]\times Y$ in temporal gauge, that is,
$\displaystyle{\frac {\partial \gamma(t)}{\partial t}}
= - \nabla {\cal C}_{\eta, \p} (\gamma(t))$, hence we have
\ba
\label{identity:2}
\begin{array}{lll}
&&{\cal C}_\eta(A(R_1), \psi(R_1)) -{\cal C}_\eta(A(R_2), \psi(R_2))
\\[2mm]
&=& -  \displaystyle{\int_{R_1}^{R_2}} < \nabla {\cal C}_\eta (\gamma(t)), 
\displaystyle{\frac {\partial \gamma(t)}{\partial t}}> dt\\[2mm]
&=&\displaystyle{ \int_{R_1}^{R_2} }< \nabla {\cal C}_\eta (\gamma(t)), 
\nabla {\cal C}_{\eta, \p} (\gamma(t))>dt\\[2mm]
&=& \displaystyle{\int_{R_1}^{R_2}}(\|\nabla{\cal C}_\eta (\gamma(t))\|^2_{L^2}
+ <  \nabla {\cal C}_\eta (\gamma(t)), \nabla \p (\gamma(t))> )dt.
\end{array}
\na
and 
\ba
\label{identity:3}
\begin{array}{lll}
&&{\cal C}_{\eta, \p}(A(R_1), \psi(R_1)) -{\cal C}_{\eta, \p}(A(R_2), \psi(R_2))\\[2mm]
&=& -\displaystyle{ \int_{R_1}^{R_2}} < \nabla {\cal C}_{\eta, \p} (\gamma(t)),
\displaystyle{\frac {\partial \gamma(t)}{\partial t}}> dt\\[2mm]
&=&\displaystyle{
\int_{R_1}^{R_2} }\| \displaystyle{\frac {\partial \gamma(t)}{\partial t}}
\|^2_{L^2}dt\\[2mm]
&=& \displaystyle{ \int_{R_1}^{R_2} }
\|\nabla {\cal C}_{\eta, \p} (\gamma(t))\|^2_{L^2}dt.
\end{array}
\na
{}From Condition \ref{condition}, we get
\[ \|\nabla {\cal C}_\eta (\gamma(t)) - \nabla {\cal C}_{\eta, \p} (\gamma(t))
\|_{L^2} < \epsilon_0 \|{\cal C}_\eta (\gamma(t))\|_{L^2}.
\]
Putting all these together and setting $\epsilon = \epsilon_0/5$, we have 
\[
1-\epsilon < \displaystyle{\frac { {\cal C}_{\eta}(A(R_1), \psi(R_1)) -{\cal C}_{\eta}(A(R_2), \psi(R_2) }{
{\cal C}_{\eta, \p}(A(R_1), \psi(R_1)) -{\cal C}_{\eta, \p}(A(R_2),\psi(R_2))}}
< 1+\epsilon.
\]
Similarly, under the condition that $\eta$ is sufficiently
small as in Remark \ref{condition:eta}, we can
prove that
\[
1-\epsilon < \displaystyle{\frac { {\cal C}_{\eta}(A(R_1), \psi(R_1)) -{\cal C}_{\eta}(A(R_2), \psi(R_2) }{
{\cal C}(A(R_1), \psi(R_1)) -{\cal C} (A(R_2),\psi(R_2))}}
< 1+\epsilon,
\]

These last two inequalities imply that $\gamma(t)$ has finite energy if and
only if any one of 
${\cal C}_{\eta}(\gamma(R_1)) -{\cal C}_{\eta}(\gamma(R_2) )$, 
${\cal C}(A(R_1), \psi(R_1)) -{\cal C} (A(R_2),\psi(R_2))$ or
$ {\cal C}_{\eta, \p}(\gamma(R_1)) 
-{\cal C}_{\eta, \p}(\gamma(R_2))$ is finite.  

For the proof of part (4),
direct calculation leads to the following identities:
\ba\label{identity:1}
\begin{array}{lll}
&&{\cal C}_\eta(A(R_1), \psi(R_1)) -{\cal C}_\eta(A(R_2), \psi(R_2))
\\[2mm]
&=&\displaystyle{ \frac 12 \int_{[R_1, R_2]\times Y}}
 (F_{\mathbb A} -*\eta )\wedge (F_{\mathbb A} -*\eta )\\[2mm]
&& + \displaystyle{\int_{\{R_1\}\times Y} <\psi, \dirac_A \psi> dvol_Y
- \int_{\{R_2\}\times Y} <\psi, \dirac_A \psi> dvol_Y}.
\end{array}
\na
and
\ba\label{identity:0}
\begin{array}{lll}
&&{\cal C}(A(R_1), \psi(R_1)) -{\cal C} (A(R_2),\psi(R_2))\\[2mm]
&=&\displaystyle{ \frac 12 \int_{[R_1, R_2]\times Y}}
F_{\mathbb A} \wedge F_{\mathbb A}\\[2mm]
&& + \displaystyle{\int_{\{R_1\}\times Y} <\psi, \dirac_A \psi> dvol_Y
- \int_{\{R_2\}\times Y} <\psi, \dirac_A \psi> dvol_Y}.
\end{array}
\na

To prove the last claim, we write the perturbed Seiberg-Witten equations
(\ref{SW:4dpert}) as follows:
\ba
\left\{
\begin{array}{l}
F^+_{\mathbb A} = q (\Psi, \Psi ) + \mu^+ + \p_{\mathbb A}^+
\\[2mm]
\Dirac_{\mathbb A} (\Psi) +\p_\Psi =0. 
\end{array}
\right.
\label{rewrite:SW}
\na
where $q (\Psi, \Psi ) = (dt \wedge \sigma (\psi, \psi ))^+$
and $(\p_{\mathbb A}^+, \p_\Psi)$ is the perturbation term coming
from the perturbation $\p$ (see (\ref{SW:4dpert}) for exact expression),
notice that, by Condition \ref{condition}, the $L^2$-norm
of $(\p_A^+, \p_\psi)$ is less than $\epsilon_0\|\nabla {\cal C}_\eta\|_{L^2}$
for a sufficiently small $\epsilon_0 <<1$.

Now we calculate
\[
\displaystyle{
\int_{[R_1, R_2]\times Y} (|\Dirac_{\mathbb A} (\Psi)|^2
+ |F^+_{\mathbb A} - q (\Psi, \Psi ) -  \mu^+|^2) dt dvol_Y}.
\] using the
3- and 4-dimensional Weitzenb\"ock formulae

The 4-dimensional Weitzenb\"ock formula is
\[
\Dirac_{\mathbb A}^* \Dirac_{\mathbb A} (\Psi) 
= \nabla ^*_{\mathbb A}\nabla_{\mathbb A}\Psi +
\displaystyle{ \frac s 4 \Psi + \frac 12 F^+_{\mathbb A}.\Psi},
\]
where $s$ denotes the scalar curvature on
$(\R\times Y, g_Y + dt^2)$. 
Take  the inner product of both sides with $\Psi = \psi (t)$, and
then integrate them over $[R_1, R_2]\times Y$. Here we need
to understand the boundary contribution when we integrate by parts. 
To see what this contribution is,   we proceed as 
follows
\[\begin{array}{lll}
&& \displaystyle{
\int_{[R_1, R_2]\times Y}}
|\Dirac_{\mathbb A} \Psi |^2 dt dvol_Y\\[2mm]
&=&
\displaystyle{
\int_{[R_1, R_2]\times Y}}\bigl( | \displaystyle{ \frac
{\partial \psi(t)}{\partial t}} + \dirac_A \psi |^2\bigr) dt dvol_Y \\[2mm]
&=& \displaystyle{
\int_{[R_1, R_2]\times Y}}\bigl(| \displaystyle{ \frac
{\partial \psi(t)}{\partial t}}|^2 + |\dirac_A \psi |^2 
+ < \displaystyle{ \frac
{\partial \psi(t)}{\partial t}}, \dirac_A \psi> +
<\dirac_A \psi, \displaystyle{ \frac
{\partial \psi(t)}{\partial t}}\bigr) dt dvol_Y\\[2mm]
&=& \displaystyle{
\int_{[R_1, R_2]\times Y}}\bigl(| \displaystyle{ \frac
{\partial \psi(t)}{\partial t}}|^2 + |\dirac_A \psi |^2 
+  \displaystyle{ \frac
{\partial }{\partial t}}(< \psi, \dirac_A \psi>) 
- \displaystyle{ \frac 12 < \frac
{\partial A(t) }{\partial t}}.\psi, \psi>\bigr) dt dvol_Y.
\end{array}\]
Now apply the 3-dimensional Weitzenb\"ock formula,
\[
\dirac^2_A\psi = \nabla_A^* \nabla_A\psi + \displaystyle{\frac s4}\psi -
\displaystyle{\frac 1 2} (*F_A).\psi,
\]
 and note that 
$F^+_{\mathbb A}. \Psi = (\displaystyle{  \frac
{\partial A(t) }{\partial t}} + *_3 F_A).\psi$,
where $*_3$ is the Hodge star operator on $Y$.  
Using the  various expressions above, we obtain that
\ba\begin{array}{lll}
&&\displaystyle{
\int_{[R_1, R_2]\times Y} |\Dirac_{\mathbb A} (\Psi)|^2 dt dvol_Y}
\\[2mm]
&=&\displaystyle{ \int_{[R_1, R_2]\times Y}}\bigl(|\nabla_{\mathbb A}\Psi|^2 + 
\frac s 4 |\Psi|^2 -\frac 12 <F^+_{\mathbb A}.\Psi, \Psi>\bigr)dt dvol_Y 
\\[2mm]
&&+ \displaystyle{\int_{\{R_2\}\times Y} <\psi, \dirac_A\psi> dvol_Y
- \int_{\{R_1\}\times Y} <\psi, \dirac_A\psi> dvol_Y}.
\end{array}
\na
We also have 
\ba\begin{array}{lll}
&&\displaystyle{
\int_{[R_1, R_2]\times Y}} 2|F^+_{\mathbb A} -q(\Psi, \Psi) -\eta^+|^2 dt dvol_Y
\\[2mm]
&= &\displaystyle{ \int_{[R_1, R_2]\times Y}}
\bigl( ( |F_{\mathbb A}|^2 + |\sigma (\psi, \psi) +\eta|^2 
+ <F^+_{\mathbb A}.\Psi, \Psi>) dt dvol_Y \\[2mm]
&& - (F_{\mathbb A} -*_3 \eta) \wedge (F_{\mathbb A} -*_3 \eta) 
+ 2 dt \wedge F_{\mathbb A}\wedge \eta.
\end{array}
\na
Using the perturbed Seiberg-Witten equations (\ref{rewrite:SW}) and
the identity (\ref{identity:1}),  we have 
\[
\begin{array}{lll}
&&\displaystyle{\int_{[R_1, R_2]\times Y}}
2(|(\p_{\mathbb A}^+, \p_\Psi)|^2) dt dvol_Y 
\\[2mm]
&=& \displaystyle{\int_{[R_1, R_2]\times Y}} 
(2 |\Dirac_{\mathbb A} (\Psi)|^2 + 2|F^+_{\mathbb A} -q(\Psi, \Psi) -\eta^+|^2)
dt dvol_Y \\[2mm]
&=& \displaystyle{\int_{[t_1, t_2]\times Y} } 
 (|F_{\mathbb A}|^2 + 2|q(\Psi, \Psi) + \eta^+|^2 + 
2|\nabla_{\mathbb A}\Psi|^2 +
\displaystyle{\frac s 2}|\Psi|^2)dt dvol_Y \\[2mm]
&& +\displaystyle{\int_{[t_1, t_2]\times Y} }\bigl(   
2 (F_{\mathbb A}\wedge dt \wedge \eta) - 
(F_{\mathbb A} -*_3 \eta) \wedge (F_{\mathbb A} -*_3 \eta)\bigr)
\\[2mm]
&& + 2 (\displaystyle{\int_{\{R_2\}\times Y} <\psi, \dirac_A\psi> dvol_Y
- \int_{\{R_1\}\times Y} <\psi, \dirac_A\psi> dvol_Y})\\[2mm]
&=& E_{[R_1, R_2]\times Y}({\mathbb A}, \Psi) + 2\bigl( 
{\cal C}_\eta(A(R_2), \psi(R_2)) - {\cal C}_\eta(A(R_1), \psi(R_1))\bigr).
\end{array}
\]
Thus, we have the following inequalities ($\epsilon = \epsilon_0/5$)
\[
2 <\displaystyle{\frac {
E_{[R_1, R_2]\times Y}({\mathbb A}, \Psi)}{ 
{\cal C}_{\eta}(A(R_1), \psi(R_1)) -{\cal C}_{\eta}(A(R_2), \psi(R_2) )} }
< 2+ \epsilon, 
\]
which implies that ${\cal C}_{\eta}(A(R_1), \psi(R_1)) -
{\cal C}_{\eta}(A(R_2), \psi(R_2) )$ is finite if and only if 
$E_{[R_1, R_2]\times Y}({\mathbb A}, \Psi)$ is finite.
The completes the proof of the Lemma.
\end{proof}

The following lemma states the basic properties of the
perturbed CSD functional ${\cal C}_{\eta, \p}$. The proof follows
the same arguments as are used for the corresponding properties of
the CSD functional ${\cal C}_{\eta}$
and of the perturbation terms satisfying 
Condition \ref{condition} for $\p\in {\cal P}_c$.

\begin{Lem} For any $\p\in {\cal P}_c$, the perturbed Chern-Simons-Dirac
functional ${\cal C}_{\eta, \p}$ on $\B$ has the following properties:
\begin{enumerate}
\item  The critical point set of ${\cal C}_{\eta, \p}$ agrees with 
the  critical point set of ${\cal C}_{\eta}$, and all points are non-degenerate.
\item  The Hessian operator of ${\cal C}_{\eta, \p}$ on $\B$ is a compact
perturbation of the Hessian operator of ${\cal C}_\eta$, hence has the
same properties as the Hessian operator of ${\cal C}_\eta$
in Lemma \ref{spectrum:Hessian}. 
\item Let $(A(t), \psi(t))$ be a solution to  the
 perturbed Seiberg-Witten equations
(\ref{SW:4dpert}) on $\R\times Y$ in temporal gauge
with finite energy, then  there
exists a sufficiently large $T$ such that on $(-\infty, -T]\times Y$
and $[T, \infty) \times Y$ the perturbation term from $\p$ is zero, in other
words, $\p(\tau,\zeta, \xi) (A(t), \psi(t))$ has compact support in
$[-T, T]$. Hence, $(A(t), \psi(t))$ is a solution to the gradient flow
equation (\ref{grad:flow}) of ${\cal C}_\eta$ outside the
compact set $[-T, T]\times Y$.
\item There is a constant $\delta>0$, such that a choice
of sufficiently small $\epsilon $ determines  $E_0>0$ with the following
property: suppose that  $(A(t), \psi(t))$ is a solution to  the
 perturbed Seiberg-Witten equations
(\ref{SW:4dpert}) on $[R_1-2, R_2+2]\times Y$ in temporal gauge
with energy
\[
E_{[R_1-2, R_2+2]\times Y} (A(t), \psi(t)) < E_0,
\]
then there exists a critical point $\alpha \in \M_Y(\s, \eta)$
such that the $L^2_1$-distance from $[A(t), \psi(t)]$ (for any
$t \in [R_1,  R_2]$) to $\alpha$ is at most
$4\epsilon\bigl(
e^{-\delta (t-R_1)} + e^{-\delta (R_2 -t)}\bigr).
$ 
Moreover,  there exists a sequence of constants
$\{\lambda_k\}_{k=0, 1, 2, \cdots}$ such that there is a 3-dimensional
monopole $(A_\alpha, \psi_\alpha)$ representing $\alpha$ for which
\[|\nabla^k (A(t)-A_\alpha)| + | (\nabla_{A_\alpha})^k (\psi(t) -\psi_\alpha)|
\leq \lambda_k \bigl(e^{-\delta (t-R_1)} + e^{-\delta (R_2 -t)}\bigr),
\]
for any point $(t, y) \in [R_1,  R_2]\times Y$.  
\item  Let $(A(t), \psi(t))$ be a solution to  the
 perturbed Seiberg-Witten equations
(\ref{SW:4dpert}) on $\R\times Y$ in temporal gauge
with finite energy, then  there
exist two 3-dimensional monopoles $(A_\alpha, \psi_\alpha)$ and
$(A_\beta, \psi_\beta)$ representing two critical points $\alpha$ and $\beta$
in $\M_Y (\s, \eta)$, such that
$A-A_\alpha$ and $\psi - \psi_\alpha$ decay exponentially along with
their first derivatives as $t \to -\infty$, and $A-A_\beta$ and $\psi - \psi_\beta$ decay exponentially along with their
first derivatives as $t \to\infty$.
\end{enumerate}
\label{CSD:p}
\end{Lem}

\begin{Cor}\label{F:finite} If $(A(t), \psi(t))$ is a solution to  the
 perturbed Seiberg-Witten equations
(\ref{SW:4dpert}) on $[0, \infty )\times Y$ in temporal gauge
with finite energy, then there exists a sequence of constants
$\{\lambda_k\}_{k=0, 1, 2, \cdots}$ such that  there is 
a 3-dimensional monopole
$(A_\alpha, \psi_\alpha)$ representing a critical point $\alpha$ in 
$\M_Y (\s, \eta)$ for which 
\[|\nabla^k (A(t)-A_\alpha)| + | (\nabla_{A_\alpha})^k (\psi(t) -\psi_\alpha)|
\leq \lambda_k e^{-\delta t},
\]
for any point $(t, y) \in [0, \infty)\times Y$, and 
\[
|\displaystyle{ \int_{[0, \infty)\times Y}}
 F_{\mathbb A} \wedge F_{\mathbb A} |< 
\infty.
\]
\end{Cor}

We now  study the geometric structure of the moduli space
$\M_\p(\alpha, \beta)$ for the gradient flow lines of ${\cal C}_{\eta, \p}$
on $\B$
connecting two different critical points $\alpha, \beta$ of $\M_Y(\s, \eta)$.
The space $\M_\p(\alpha, \beta)$ is also the moduli space of
solutions to the perturbed Seiberg-Witten equations on $(Y\times \R,
g + dt^2, \pi^*(\s))$, whose asymptotic limits at $\pm \infty$
represent $\alpha$ and $\beta$ respectively.
Let $\delta$ be  the  decay rate of the gradient flow line
in $\M_\p(\alpha, \beta)$
(see Proposition \ref{decay:expo} and Lemma \ref{CSD:p}).

We will consider the space of pairs consisting of $U(1)$ connections
on the
determinant bundle  and  sections of the spinor bundle
on \[\Bigl(Y\times \R, g+dt^2, \pi^*(\s)\Bigr),\]
  topologized with
weighted Sobolev norms of weight $\delta$ as in \cite{LMc}. Here the
weight function is $e_\delta(t)=e^{\tilde\delta t}$,
where $\tilde\delta$ is a smooth
function with bounded derivatives, $\tilde\delta : \R\to [-\delta,\delta]$
such that $\tilde\delta (t)\equiv -\delta$ for $t\leq -1$ and $\tilde\delta
(t)\equiv \delta$ for $t\geq 1$.
The $L^2_{k,\delta}$ norm is defined as 
\[ \| f \|_{L^2_{k,\delta}} =\Bigl(\int_{Y\times \R}e_\delta (t) 
(|f|^2 + |\nab f|^2 + |\nab ^2 f|^2 + \cdots +  |\nab ^k f|^2)
 dt dvol_Y \Bigr)^{\frac 12} 
\]
The weight $e_\delta$ imposes an exponential decay
as an  asymptotic condition along the cylinder.

An element $\gamma_0$ of $\M(\alpha, \beta )$ determines a path in $\A/\G$,
also denoted $\gamma_0$. We lift $\gamma_0$ to
a path in $\A$
which we denote by $\tilde\gamma_0 =  (A_0(t), \psi_0(t))$. From 
Lemma \ref{CSD:p}, we know that $(A_\alpha, \psi_\alpha)$ and
$(A_\beta, \psi_\beta)$ represent $\alpha$ and $\beta$ respectively, such that
$(A_0(t), \psi_0(t))$ approaches to $(A_\alpha, \psi_\alpha)$ and
$(A_\beta, \psi_\beta)$  exponentially fast in the $L^2_1$-topology
as $t\to \pm\infty$ respectively.

The operators $T(t)$ have a spectral flow denoted by
$\ind (\gamma_0 )$
along the path $\tilde\gamma_0$ (see (\ref{spec:T}) for the relation with
the spectral flow of the Hessian operator $Q$). We see that
\[
\ind (\gamma_0) = SF_{\gamma_0} (Q) = \ind (\alpha, \beta ) ( mod\  d(\s))
\]
where $\ind (\alpha, \beta )$ is the relative index
between $\alpha$ and $\beta $, defined in (\ref{index:rel}). 
Note that $\ind (\gamma_0)$ is $\Z$-valued.

Choose a $U(1)$-connection $\mathbb A_0$ on $det (W^+)$ and
a spinor section $\Psi_0$ of $W^+$ such that $(A_0, \Psi_0)$
is in temporal gauge outside a compact set and satisfies 
$({\mathbb A}_0, \Psi_0)= (A_\alpha, \psi_\alpha)$ for $t<< 0$ and 
$({\mathbb A}_0, \Psi_0)= (A_\beta, \psi_\beta)$ for $t>>0$.
 
For $k\ge 2$, let ${\cal A}_{k,\delta}(\alpha, \beta )$ be the space of pairs
of connections and spinor sections $(A, \psi)$
on $Y \times \R$ satisfying
\[ (A, \psi) \in ({\mathbb A}_0, \Psi_0) + (\Omega^1_{L^2_{k,\delta}}
(Y\times\R, i\R) \oplus L^2_{k,\delta} (W^+) ). \]
The gauge transformation group ${\cal G}_{k+1,\delta}(\alpha, \beta)$ is
locally
modelled on \[ L^2_{k+1,\delta}(\Omega^0(Y\times\R, i\R)) \]
and approaches elements
in the stabilizers $G_\alpha$ and $G_\beta$ of $\alpha$ and $\beta$
as $t\to \pm\infty$. This gauge
group acts on ${\cal A}_{k,\delta}(\alpha,\beta)$  freely,
so we can form the quotient ${\cal B}_{k,\delta}(\alpha,\beta)$,
a smooth Banach manifold, whose tangent space at $[A, \psi]$ (represented
by $(A, \psi)\in {\cal A}_{k,\delta}(\alpha,\beta)$)
is given by
\[
\{ (a, \phi)\in  \Omega^1_{L^2_{k, \delta}} \oplus L^2_{k, \delta}(W^+)
\ |G^*_{[A, \psi], \delta} (a, \phi) = 
(e_{-\delta}d^*  e_{\delta} ) a + iIm \la \psi , \phi \ra = 0  \}.
\]

Using Proposition \ref{decay:expo}, we know that
the component of ${\cal M}_\p(\alpha, \beta )$ containing the given
gradient flow line  $\gamma_0$, denoted by
${\cal M}_\p(\alpha, \beta )_{\gamma_0}$,
can be identified with
\[
\left\{ (A, \psi ) \in {\cal A}_{k, \delta} (\alpha, \beta)
\left| \begin{array}{l}
(A, \psi )  \hbox{ satisfies the }\\[2mm]
 \hbox{monopole equations (\ref{SW:4dpert})}.
\end{array}\right.
\right\} \slash  {\cal G}_{k+1, \delta}(\alpha, \beta).
\]

\begin{Pro} \label{transversal:M}
There is a Baire set ${\cal P}_0$ of 
smooth functions  $\p\in {\cal P}_c$, satisfying
Condition \ref{condition}, such that  for $\p \in {\cal P}_0$
${\cal M}_\p(\alpha, \beta )_{\gamma_0}$, the component
of ${\cal M}(\alpha, \beta )$ containing $\gamma_0$, if non-empty,
is a smooth oriented manifold of dimension $\ind (\gamma_0) - d_\alpha$ 
where $d_\alpha= 1$ if $\alpha$ is reducible and $d_\alpha=0$ otherwise.
\end{Pro}
\begin{proof}
Suppose $(A(t), \psi (t), \p)$, representing a point 
in ${\cal M}_\p(\alpha, \beta )_{\gamma_0}$, is a solution 
to the perturbed gradient
flow equation (\ref{grad:pert}) or the perturbed Seiberg-Witten 
equations (\ref{SW:4dpert}) in temporal gauge with the perturbation
$\p \in {\cal P}_c$. Using the linearization
of (\ref{SW:4dpert}) and gauge transformations
at $(A(t), \psi (t), \p)$ we may
define the following operator
\[\begin{array}{c}
{\cal D}_{A, \psi,  \p}: \quad {\cal P}_c \oplus
L^2_{k, \delta}( i\Omega^1 \oplus W^+) \to
 L^2_{k-1, \delta}(i\Omega^0 \oplus i \Omega^{2,+}\oplus  W^-),
\\[2mm]
{\cal D}_{A, \psi,  \p} (\tilde p, a, \phi) 
= \bigl( G^*_{(A, \psi), \delta}, d^+ a - (dt\wedge \sigma(\phi, \psi))^+
+ \p^+, \Dirac_A\phi +\frac 12 a.\psi+ \p_\psi\bigr),
\end{array}
\]
where $(\p^+, \p_\psi)$ is given by the linearization of the perturbation
term at $(A(t), \psi (t), \p)$.

Let $\pi: Y \times \R \to Y$ be the projection map, then as bundles over $Y \times \R$,
\[
\Lambda^1 (Y \times \R, i \R) \cong \pi^* \Lambda^1 (Y, i\R) \oplus
\pi^*(\Lambda^0 (Y, i\R)
\]
\[
\Lambda^{2, +} (Y \times \R, i \R) \cong \pi^* \Lambda^1 (Y, i\R).
\]
Also, using Clifford multiplication by $dt$, we identify
$W^+ \cong W^- \cong \pi^* W$.
Under these identifications,
fix $\p$, then the operator ${\cal D}_{A, \psi,  \p}$, restricted to
$L^2_{k, \delta}( i\Omega^1 \oplus W^+)$,  has
the form 
\[
\displaystyle {\frac {\partial }{\partial t} }  
+ \left( \begin{array}{cc}
L_{(A, \psi)} & G_{(A, \psi)} \\ G^*_{(A, \psi)}
& 0 \end{array}\right)  +\delta + o(1)
\]
where $L_{(A, \psi)}, G_{(A, \psi)}$ and $G^*_{(A, \psi)}$
are given by (\ref{complex:3d1}). By the result of Lockhart-McOwen
\cite{LMc}, as $\p$ varies through ${\cal P}_c$, 
${\cal D}_{A, \psi,  \p}$,  restricted to
$L^2_{k, \delta}( i\Omega^1 \oplus W^+)$, forms
 a continuous family of Fredholm operators
with index equal to $\ind (\gamma_0) -d_\alpha$.

Now we prove that ${\cal D}_{A(t), \psi(t),  \p}$ 
is surjective.
Note that any solution $(A(t), \psi(t))$ which represents a point 
in ${\cal M}_\p(\alpha, \beta )_{\gamma_0}$, is a solution
to the perturbed gradient
flow equation (\ref{grad:pert}) with finite energy 
${\cal C}_\eta (A_\alpha, \psi_\alpha) - {\cal C}_\eta (A_\beta,
\psi_\beta)$. We can use Condition \ref{condition}, to assert that there is
a compact subset $[-T, T]\times Y$ in $\R\times Y$, such
that outside this compact set, the perturbation function $\p\in 
{\cal P}_c$ is zero.

Suppose that $(f, b, \phi) \in L^2_{-\delta}(i\Omega^0 \oplus i\Omega^{2, +}
\oplus W^-)$ is $L^2$-orthogonal to the image of 
${\cal D}_{A(t), \psi(t),  \p}$, then using Lemma \ref{injectivity}
and  a direct calculation shows
that $(f, b, \phi)$ is zero on $[-T, T]\times Y$.
As an element in the cokernel of ${\cal D}_{A(t), \psi(t),  \p}$,
$(f, b, \phi)$ satisfies the weak unique continuation
principle \cite{BMW} and hence $(f, b, \phi)$ is zero everywhere. This proves
that ${\cal D}_{A(t), \psi(t),  \p}$
is surjective.

The Sard-Smale theorem implies that there exists a Baire set ${\cal P}_0$
 of functions $\p\in {\cal P}_c$, so that
the moduli space $\M_\p (\alpha, \beta )_{\gamma_0}$,
for $\p\in {\cal P}_0$, is a smooth
manifold of dimension $\ind (\gamma_0) - d_\alpha$ if
$\ind (\gamma_0) - d_\alpha > 0 $, and is empty otherwise.

The orientation of $\M_\p(\alpha, \beta)_{\gamma_0}$ is determined
by the determinant line bundle of ${\cal D}_{A, \psi,  \p}$ over
$\A_{k, \delta}(\alpha, \beta)$ for a fixed perturbation $\p$.
This is established in Proposition 2.15 of \cite{MW} where it is
shown that $\M_\p(\alpha, \beta)_{\gamma_0}$ is orientable, 
the orientation is determined  by choosing an orientation of
\[
\wedge^{top} H_{\delta}^1(Y\times \R, i\R) \otimes
\wedge^{top} H_{\delta}^{2, +}(Y\times \R, i\R).
 \]
\end{proof}

Proposition \ref{transversal:M} has an immediate corollary.

\begin{Cor}\label{moduli:least} Let  $(Y, g, \s)$
be a closed, oriented 3-manifold with
a Riemannian metric $g$ and a $\spinc$ structure $ \s$. 
Let $\M(\alpha, \beta)$ be the
moduli space of gradient flows of ${\cal C}_{\eta, \p}$
which connect two different critical points $\alpha$ and $\beta$.  Then,
for a generic perturbation 
$\p$, we have the following results. 
\begin{enumerate}
      \item \label{M:Z}  When $Y$ is a homology 3-sphere,
$\M(\alpha, \beta)$ is empty if $\ind (\alpha) - \ind (\beta) - d_\alpha \le 0$,
otherwise, $\M(\alpha, \beta)$ is an oriented, smooth manifold of dimension
$\ind (\alpha) - \ind (\beta) - d_\alpha$.
       \item \label{M:Zd} If $c_1(\s)$ is non-zero in $H^2(Y, \R)$, 
$\M(\alpha, \beta)$
has infinitely  many components. 
Each of them is an oriented, smooth manifold with
dimension given by $\ind (\gamma) - d_\alpha > 0$, where $\gamma$ is
a chosen element in that component. The dimensions of two non-empty
components differ by a multiple of $d(\s)$, where $d(\s)$
is the divisibility of $c_1(\s)$ in $H^2(Y, \Z)/Torsion$.
\end{enumerate}
\end{Cor}

Now we discuss  the compactness of the trajectory moduli space. For 
the case of a 3-manifold
$Y$ with $b_1(Y) > 0$ and $c_1(\s) \neq 0$ in $H^2(Y, \R)$, we know
that $ \hat{\cal M}(\alpha, \beta)$, the quotient of 
${\cal M}(\alpha, \beta)$ by $\R$, has infinitely many components. Assume 
$\ind (\alpha) - \ind (\beta) = k +1 (mod \  d(\s))$ with $k\ge0$, then
\[
\hat{\cal M}(\alpha, \beta) = \cup_{n\in \N} \hat\M^{k+ nd(\s)}
(\alpha, \beta)
\]
where $ \hat\M^{k+ nd(\s)} (\alpha, \beta)$ is the union of components of
dimension $k+ nd(\s)$ in $\hat{\cal M}(\alpha, \beta)$.
   We now summarize the compactness 
results and refer to page 501-549 in \cite{MW} for the proofs.

\begin{Pro} \label{compactification} (Theorem 4.23 \cite{MW})
For any 3-manifold $(Y, g, \s)$ with a Riemannian metric $g$ and a $\spinc$
structure $\s$, assume $c_1(\s) \neq 0$ if $b_1(Y)>0$, then 
$\hat \M^{k+nd(\s)}(\alpha, \beta)$ ($ \hat\M (\alpha, \beta)$ 
for $b_1(Y)=0$) can be compactified by adding lower dimensional boundary 
strata of broken trajectory moduli spaces.
Namely, for $b_1(Y)=0$, the  boundary strata are of 
the form 
\[
\bigcup_{\alpha_1, \cdots \alpha_j}
\hat \M(\alpha, \alpha_1) \times \hat \M(\alpha_1, \alpha_2)
\times \cdots \times \hat \M(\alpha_j, \beta).
\]
Here the union is over all possible sequences of critical points
$\alpha, \alpha_1, \cdots \alpha_j, \beta$ with decreasing indices. 
For $b_1(Y)>0$ and
$c_1(Y) \neq 0$,  $\hat \M^{k+nd(\s)}(\alpha, \beta)$ can be
compactified by adding boundary strata of 
the form
\[
\bigcup_{\alpha_1, \cdots \alpha_j}
\hat \M^{k_0} (\alpha, \alpha_1) \times \hat \M ^{k_1}(\alpha_1, \alpha_2)
\times \cdots \times \hat \M ^{k_j}(\alpha_j, \beta).
\]
Here the union is over all possible sequences of critical points
$\alpha, \alpha_1, \cdots \alpha_j, \beta$ with $k_0 + k_1 +\cdots + k_j
\leq k+ n d(\s) -j$ and $k_i = \ind (\alpha_i, \alpha_{i+1}) (mod \ d (\s))$.
\end{Pro}

With these compactness results 
for the moduli space, we have the following
corollary which will be crucial for the construction of Seiberg-Witten-Floer
homology.

\begin{Cor}\begin{enumerate}
\item  Let $(Y, \s)$ be a homology 3-sphere or $b_1(Y)>0$ with 
$c_1(\s)\neq 0$, then the zero dimensional components 
 $\hat \M^0 (\alpha, \beta)$ for two critical points $\alpha$
and $\beta$ of relative index $1$  are compact, and consist of
finitely many oriented points.
\item  Let $(Y, \s)$ be a homology 3-sphere or $b_1(Y)>0$ with 
$c_1(\s)\neq 0$, 
let $\alpha, \gamma$ be two irreducible critical points
in $\M_Y(\s, \eta)$ with relative  index
$ \ind (\alpha, \gamma) = 2$,  or 2 $mod \ d(\s)$ for $b_1(Y)>0$.
Then the boundary of $ \hat\M^1(\alpha, \gamma)$ (an oriented,
compact 1-manifold) consists of the union
\[
\bigcup_{\beta \in \M_{\s, \p}^*}
\hat\M^0(\alpha, \beta) \times \hat\M^0(\beta, \gamma)
\]
where $\beta$ runs over those critical points $\beta$ with
relative index $\ind(\alpha, \beta) =1$ or $1$ $ mod\ d(\s)$ for $b_1(Y)>0$. 
\end{enumerate}
\label{boundary:1}
\end{Cor}

\begin{Rem} 
For a 3-manifold $(Y, \s)$ with $b_1(Y) >0$ and $c_1(\s) =0$, we
need to perturb by a co-closed 2-form $\eta$ such that $[*\eta]$ representing a
non-trivial cohomology class. Then ${\cal C}_\eta$
has finitely many critical points $\M_Y(\s, \eta)$, which are
all irreducible and have $\Z$-valued indices. Note that in this case
${\cal C}_\eta (g. (A, \psi)) - {\cal C}_\eta (A, \psi)
= \la [g] \cup [\eta], [Y] \ra$, 
which is a
multiple of a positive number $e_\eta$ ($e_\eta$ depends on
$\eta$). For two critical points $\alpha$ and $ \beta$ with relative
index $k+1>0$, and for generic $\p$,
  $\hat \M(\alpha, \beta)$
is an oriented and smooth manifold of dimension $k$, 
$\hat \M(\alpha, \beta) = \bigcup_{n\in \N} \hat \M_{(n)}(\alpha, \beta)$
where $\hat \M_{(0)}(\alpha, \beta)$ is the union of components
in $\hat \M (\alpha, \beta)$ with minimal energy $e_0(\alpha, \beta)$,
and for $n>0$, $\hat \M_{(n)}(\alpha, \beta)$ 
is the the union of components
in $\hat \M (\alpha, \beta)$ with energy 
$e_0(\alpha, \beta) +n e_\eta$. Moreover $\hat \M_{(n)}(\alpha, \beta)$
can be compactified by adding boundary strata of
the form
\[
\bigcup_{\alpha_1, \cdots, \alpha_l} \hat \M_{(i_0)}(\alpha, \alpha_1)
\times \hat \M_{(i_1)} (\alpha_1, \alpha_2) \times
\cdots \times \hat \M_{(i_l)}(\alpha_l, \beta).
\]
Here the union is over all possible sequences of
$\alpha_0=\alpha, \alpha_1, \cdots, \alpha_l,  \alpha_{l+1} =\beta$
with decreasing indices and $\sum_{j=0}^{l} e_0(\alpha_j, \alpha_{j+1})
+\sum_{j=0}^{l} i_j e_\eta \leq e_0(\alpha, \beta) + n e_\eta$.
\end{Rem}

\subsection{Seiberg-Witten-Floer homology: definition and properties}
\label{SWF:def}
In this section, we will only construct the Seiberg-Witten-Floer
homology for $(Y, \s)$ when  $(Y, \s)$ 
has $b_1 (Y) >0$ and $c_1(\s) \neq 0$. When $Y$ is a homology sphere,
there is an equivariant Seiberg-Witten-Floer
homology, developed in \cite{MW}, which is a topological invariant. 

Let $(Y, g, \s)$ be an oriented, closed 3-manifold with a
Riemannian metric $g$ and a $\spinc$ structure $\s$. 
For a generic $\eta$, the critical points of ${\cal C}_\p$ consist
of finitely many, non-degenerate points, which we 
denote by $\M_Y(\s, \eta)$, all irreducible.

The Floer complex $C_*^{SW}(Y, \s)$ 
is generated freely by the critical points in $\M_Y(\s, \eta)$
with a relative $d(\s)$-grading
\[
C_{*}(Y, \s) = \oplus_{\alpha \in \M_Y(\s, \eta)} \Z. <\alpha>.
\]
The boundary operator $\d$ on $C_*(Y, \s)$ is defined  by
\[ \partial (<\alpha>) = \sum_{\beta\in \M_Y(\s, \eta)}
n_{\alpha\beta}<\beta>, \]
where $n_{\alpha\beta}$ is given by counting the points
in $\hat\M^0(\alpha, \beta)$ (an oriented, smooth, compact 0-manifold
from Corollary \ref{boundary:1})
with sign.

\begin{Lem} $\partial\circ\partial =0$. \end{Lem}
\begin{proof} By definition,
\[
\begin{array}{rcl}
\partial^2 (<\alpha>)&=&
\displaystyle{\sum_{
\beta\in \M_Y(\s, \eta)}}
 n_{\alpha \beta} \partial(<\beta>)\\[2mm]
&=& \displaystyle{\sum_{
\beta \in \M_Y(\s, \eta)}
\sum_{ \gamma\in \M_Y(\s, \eta) }}
n_{\alpha \beta}n_{\beta\gamma}<\gamma>.
\end{array}
\]
where $\beta$ runs over the critical points with 
 relative
index $\ind(\alpha, \beta) = 1 mod (\d(\s))$,
$\gamma$ runs over the critical points with 
 relative index $\ind(\alpha, \gamma)= 2 mod (\d(\s))$. We want to show that
\[  \displaystyle{\sum_{
\{\beta | \ind (\alpha, \beta) =1 mod(d(\s))\} } }
n_{\alpha \beta}n_{\beta\gamma} = 0\]
for any $\gamma \in \M_Y(\s, \eta) $ with
$\ind (\alpha, \gamma)= 2 mod  (d(\s))$.
We know that  the number
\[
\displaystyle{\sum_{
\{\beta |\ind (\alpha, \beta) =1 mod(d(\s))\} }}
 n_{\alpha \beta}n_{\beta\gamma}
\]
is the number of oriented boundary points of the 1-manifold 
$ \hat \M^1(\alpha, \gamma)$
by Corollary \ref{boundary:1},
and hence is zero.
\end{proof}

Now we can define the Seiberg-Witten-Floer homology.
\begin{Def}\label{non-equ:SWF} For any closed oriented 3-manifold
$Y$ with a $\spinc$ structure,
the Seiberg-Witten-Floer homology is defined to be
\ba
HF^{SW}_*(Y,\s) = H_*(C^{SW}_*(Y,\s), \partial_*),
\na
which is a $\Z_{d(\s)}$-graded and finitely generated Abelian group.
\end{Def}

\begin{Rem} 
The Seiberg-Witten-Floer homology can be thought as
 a refinement of
the Seiberg-Witten invariant for a 3-manifold $Y$ with a $\spinc$
structure $\s$, 
in the sense that
the Seiberg-Witten invariant for $(Y, \s)$
is  the Euler characteristic of the Seiberg-Witten-Floer
homology $HF^{SW}_*(Y, \s)$.\end{Rem}

The definition of $HF^{SW}_*(Y,\s)$
involves  both the metric and perturbations.
The following proposition shows
that for  a 3-manifold $Y$ with $b_1(Y) > 0$ and $c_1(\s)\neq 0$, 
the Seiberg-Witten-Floer homology
is a topological invariant, independent of the
metric and generic perturbations $\eta$ and $\p$. 
We then show that there is an the action of
$\AA(Y) = Sym^*(H_0(Y, \Z))\otimes \Lambda ^*(H_1(Y, \Z)/torsion)$
on $HF^{SW}_*(Y, \s)$.

\begin{Pro}\label{U-action}
For any closed oriented 3-manifold $Y$ with 
a $\spinc$ structure $\s$ such that $c_1(\s)$ is non-torsion, 
then $HF^{SW}_*(Y, \s)$ is a topological invariant and 
there is an action of $\Z[U]$ on $HF^{SW}_*(Y, \s)$
where the $U$-action decreases the relative grading by $2$.
\end{Pro}
\begin{proof} The proof of topological invariance can be
found in \cite{MW}. The $\Z[U]$-module structure can also be extracted
from the isomorphism between $HF^{SW}_*(Y, \s)$ and the
corresponding equivariant Seiberg-Witten-Floer homology
$HF^{SW}_{*, U(1)}(Y, \s)$. Here
we sketch a more direct proof of the existence of a $\Z[U]$-module structure
on $HF^{SW}_*(Y, \s)$. In subsection 5.3 of \cite{MW}, there is
an associated integer $m_{\alpha\gamma}$ for each pair of
critical points $\alpha$ and $\gamma$ in $\M_Y(\s, \eta)$ of
relative index 2 $mod (d(\s))$, defined
as follows.

Let $\M^2(\alpha, \beta)$ be the 2-dimensional components of
the moduli space of parametrized flowlines between $\alpha$ and $\gamma$. There
is a $U(1)$-principal bundle over $\M^2(\alpha, \beta)$,
denoted by $\tilde \M^2(\alpha, \beta)$, which
is the based moduli space corresponding to $\M^2(\alpha, \beta)$.
That is, fix a base point $(t_0, y_0) \in \R \times Y$ and
a complex line $L_{y_0}$ in the fiber $W_{y_0}$ of the
spinor bundle $W$ over $Y$ such that 
$L_{y_0}$ doesn't contain the spinor part of any critical
point in $\M_Y(\s, \eta)$. Such
$L_{y_0}$ exists due to the fact that
$\M_Y(\s, \eta)$ consists of finitely many irreducible
critical points. Note also that the space of such $L_{y_0}$
is path-connected. Then there is a complex line bundle
${\cal L}_{\alpha\gamma}$ over $\M^2(\alpha, \beta)$
\[
{\cal L}_{\alpha\gamma} 
= \tilde \M^2(\alpha, \beta) \times _{U(1)} (W_{y_0}/L_{y_0})
\]
with a canonical  section 
\[
s_{\alpha\gamma}([{\mathbb A}, \Psi]) = ([{\mathbb A}, \Psi], \Psi (t_0, y_0)).
\]
This canonical section determines a trivialization of
${\cal L}_{\alpha\gamma}$ away from a compact set (see Lemma 5.7
in \cite{MW}) by the choice of $L_{y_0}$. 
Then the relative Euler number of $({\cal L}_{\alpha\gamma},
s_{\alpha\gamma})$ defines an integer $m_{\alpha\gamma}$ for each pair of
$\alpha$ and $\gamma$ with relative index 2 $mod (d(\s))$.
The system of integers $\{m_{\alpha\gamma}: \ind (\alpha, \gamma)
= 2 mod (d(\s))\}$ satisfies the following properties (Cf. Lemma 5.7
and Remark 5.8 in \cite{MW}): $m_{\alpha\gamma}$ is independent of the 
choices of base points and $L_{y_0}$, and for
$\alpha$ and $\delta$ two critical points with relative
index 3 $mod (d(\s))$, then 
\ba
\label{property:m}
\sum_{\{\beta | \ind (\alpha, \beta) =1 mod (d(\s))\}} n_{\alpha\beta}
m_{\beta\delta} - \sum_{\{\gamma| \ind (\alpha \gamma) =2 mod( d(\s))\}}
m_{\alpha\gamma}n_{\gamma \delta} =0.
\na
Now we can define a map 
\ba\label{U}
\begin{array}{cccc}
U: &\qquad C_*^{SW} (Y, \s)& \longrightarrow&
 C_{*-2}^{SW} (Y, \s) \\
&<\alpha>& \mapsto& 
\displaystyle{\sum_{
\{\gamma |
\ind (\alpha, \gamma) =1 mod(d(\s))\}}
} m_{\alpha\gamma} <\gamma>.
\end{array}
\na
Identity (\ref{property:m}) implies that $U$ induces
a map on the Seiberg-Witten-Floer homology which decreases
the relative grading by $2$. This is the $\Z[U]$ module structure
on $HF_*^{SW}(Y, \s)$.
\end{proof}

Now we define the action of $\AA(Y) = Sym^*(H_0(Y, \Z))\otimes
\Lambda^* (H_1(Y, \Z)/Torsion)$ on $HF_*^{SW}(Y, \s)$. 
Note that $Sym^*(H_0(Y, \Z))
\cong \Z[U]$, so we only need to define an action of
$\Lambda^* (H_1(Y, \Z)/Torsion)$ on $HF_*^{SW}(Y, \s)$.

Let $\pi$ be an embedded loop in $Y$ which represents a non-torsion
element in $H_1(Y, \Z)$. For any $\alpha$ and $\beta$  two critical points
in $\M_Y(\s, \eta)$ with relative index 1 $mod (d(\s))$, let
$\M^1(\alpha, \beta)$ denote the 1-dimensional components of
the moduli space of parametrized flowlines between 
$\alpha$ and $\beta$.  Then there is 
a universal line bundle over $\pi \times \M^1(\alpha, \beta)$
constructed as follows: for any point $y\in \pi$, the based
moduli space with base point $(t_0, y)$ defines a $U(1)$-bundle
over $\M^1(\alpha, \beta)$. Varying $y$ in $\pi$ gives rise
to a universal line $\widetilde {\M^1_\pi(\alpha, \beta)}$
 over $\pi \times \M^1(\alpha, \beta)$.

As in the proof of Proposition \ref{U-action}, choose
a fixed base point $y_0$ and a complex line $L_{y_0} \subset W_{y_0}$,
then there
is an associated complex line bundle 
\[
 {\cal L}^{\pi}_{\alpha\beta}
= \widetilde {\M^1_\pi(\alpha, \beta)}\times_{U(1)} (W_{y_0}/L_{y_0})
\]
over $\pi \times \M^1(\alpha, \beta)$, 
endowed with a canonical section $s_{\alpha\beta}^\pi$ defined
in a fashion similar to that in the proof of Proposition \ref{U-action}.
This canonical section $s_{\alpha\beta}^\pi$  also defines
a trivialization away from a compact set in $\pi \times \M^1(\alpha, \beta)$.
The corresponding relative Euler number is denoted
by $\pi_{\alpha\beta}$: it counts the
zeros of a generic section which agrees with the trivialization
determined by $s_{\alpha\beta}^\pi$. The
standard cobordism argument implies that $\pi_{\alpha\beta}$
is independent of the choice of embedded loop $\pi$, hence
we can denote this relative Euler number
as $[\pi]_{\alpha\beta}$ for any non-torsion element 
$[\pi] \in H_1(Y, \Z)$. This defines a map of degree $-1$ on
the Seiberg-Witten-Floer chain complex: 
\ba\label{pi-action}
\begin{array}{ccccc}
[\pi]&\qquad& C_*^{SW}(Y, \s)&\longrightarrow& C_{*-1}^{SW}(Y, \s)\\[2mm]
&&<\alpha> &\mapsto & \displaystyle{ \sum_{
\{\beta | \ind (\alpha, \beta) =1 mod(d(\s)) \}}}
[\pi]_{\alpha\beta} <\beta>.
\end{array}
\na

\begin{Pro} $[\pi]$ is a chain map and satisfies $[\pi]\circ [\pi] =0$,
hence it defines an action of $\Lambda^*(H_1(Y, \Z)/Torsion)$ on the
Seiberg-Witten-Floer homology $HF^{SW}_*(Y, \s)$.
\end{Pro}
\begin{proof}
First, we show that $[\pi]$ is a chain map. Let $\alpha$ and $\gamma$
be two critical points of relative index $2$ $ mod (d(\s))$. Let
$\M^2(\alpha, \gamma)^*$ be the induced partial compactification
of the 2-dimensional components of $\M (\alpha, \gamma)$ obtained
by adding the broken flowlines. Note that, $\M^2(\alpha, \gamma)^*$ 
consists of finitely many cylinders (homeomorphic
to $S^1\times \R$), and finitely many bands (homeomorphic
to $[0, 1]\times \R$). The boundary ends in the bands correspond
to the unparametrised broken flowlines
\[
\bigcup_{\{\beta | \ind (\alpha, \beta) =1 mod (d(\s))\}}
\hat \M^0(\alpha, \beta)\times \hat \M^0(\beta, \gamma),
\]
where $\hat \M^0(\alpha, \beta)$ and $\hat\M^0(\beta, \gamma)$ are
the quotients of $\M^1(\alpha, \beta)$ and $\hat\M^1(\beta, \gamma)$
by the $\R$-action respectively.
There is a complex line bundle ${\cal L}_{\alpha\gamma}^\pi$ 
over $\pi\times \M^2(\alpha, \gamma)^*$ with a
canonical section $s_{\alpha\gamma}^\pi$, constructed in a fashion
similar to that of
$({\cal L}_{\alpha\beta}^\pi, s_{\alpha\beta}^\pi)$
in the proof of Proposition 3.24. Away from
a compact set in $\pi\times \M^2(\alpha, \gamma)^*$, $s_{\alpha\gamma}^\pi$
is nowhere vanishing, hence it determines a trivialization 
of ${\cal L}_{\alpha\gamma}^\pi$ away from a compact set. 
The zeros of a section which agrees with 
$s_{\alpha\gamma}^\pi$ away from a compact set form an
oriented smooth 1-manifold. Such sections we refer to as 
generic sections of $({\cal L}_{\alpha\gamma}^\pi, s_{\alpha\gamma}^\pi)$.
Looking at the boundary of this zero set
we see that it consists of 
\[
\bigcup_{\{\beta | \ind (\alpha, \beta) =1 mod (d(\s))\}}
\bigl( \hat s_{\alpha\beta}^{-1} (0) \times \hat \M^0 (\beta, \gamma)
\cup \M^0 (\alpha, \beta) \times \hat s_{\beta\gamma}^{-1} (0)\bigr).
\]
Here $\hat s_{\alpha\beta}$ and $\hat s_{\beta\gamma}$
are generic sections of
 $({\cal L}_{\alpha\beta}^\pi, s_{\alpha\beta}^\pi)$
and $({\cal L}_{\beta\gamma}^\pi, s_{\beta\gamma}^\pi)$ respectively.
The counting of these boundary points implies that
\[
\sum_{\{\beta | \ind (\alpha, \beta) =1 mod (d(\s))\}}
[\pi]_{\alpha\beta} n_{\beta\gamma} 
- \sum_{\{\beta'| \ind (\alpha, \beta') =1 mod (d(\s))\}}
n_{\alpha\beta'}[\pi]_{\beta'\gamma} =0.
\]
Hence $[\pi]$ is a chain map so we have defined an action of
$H_1(Y, \Z)/Torsion$ on $HF^{SW}_*(Y, \s)$.

To get $[\pi]\circ [\pi] =0$ on the chain level, we need to give another
way of calculating the relative Euler number of $({\cal L}_{\alpha\beta}^\pi,
s_{\alpha\beta}^\pi)$ in terms of the following holonomy map:
\[
h_{\alpha\beta}^\pi: \qquad \M^1( \alpha, \beta) \longrightarrow
U(1),
\]
where $h_{\alpha\beta}^\pi ([{\mathbb A}, \Psi])$
is given by the holonomy of ${\mathbb A}$ along the loop $\{t_0\}\times \pi$
for any $[{\mathbb A}, \Psi] \in \M^1( \alpha, \beta)$.
Let $\omega$ be the standard closed 1-form on $U(1)$ with integral one
over $U(1)$.
Then pulling back:
\[
[\pi]_{\alpha\beta} = \displaystyle{ \int_{\M^1( \alpha, \beta)}}
(h_{\alpha\beta}^\pi)^* \omega 
\]
which is the winding number of the holonomy map $h_{\alpha\beta}^\pi$,
that is, the counting of a generic fiber of $h_{\alpha\beta}^\pi$. 

Now for two critical points $\alpha$ and $\gamma$ with relative
index $2$ $ mod (d(\s))$, the corresponding holonomy map
can be extended over the partial compactification of
$\M^2(\alpha, \gamma)^*$. We denote this holonomy map
by $h_{\alpha\gamma}^\pi$. Then the fiber of $h_{\alpha\gamma}^\pi$
at a generic point $p\in U(1)$ is an oriented smooth 1-manifold
whose boundary consists of
\[
\bigcup_{\{\beta | \ind (\alpha, \beta) =1 mod (d(\s))\}}
(h_{\alpha\beta}^\pi)^{-1} (p) \times (h_{\beta\gamma}^\pi)^{-1} (p). 
\]
Counting these boundary points with sign, we obtain 
\[
\sum_{\{\beta | \ind (\alpha, \beta) =1 mod (d(\s))\}}
[\pi]_{\alpha\beta} [\pi]_{\beta\gamma} =0.
\]
This shows that $[\pi]\circ [\pi] =0$ on the chain level. Hence
there is an action of $\Lambda^*(H_1(Y, \Z)/Torsion)$ on
the Seiberg-Witten-Floer homology $HF^{SW}_*(Y, \s)$.
\end{proof}

Let $(-Y, -\s)$ be the 3-manifold $Y$ with the reversed orientation. There
is a one-to-one correspondence between the 3-dimensional monopoles
on $(Y, \s)$ and $(-Y, -\s)$
\[
\M_{Y} (\s, \eta) \cong \M_{-Y}(-\s, -\eta)
\]
which sends $\alpha \mapsto \hat\alpha$,  such that 
\[ \ind_Y (\alpha, \beta) = \ind_{-Y} (\hat\beta, \hat\alpha).
\]
 There is also an orientation
preserving diffeomorphism between the moduli space of gradient 
flow lines for the Chern-Simons-Dirac functionals on
$(Y, \s)$ and $(-Y, -\s)$:
\[
\M_{\R\times Y}(\alpha, \beta) \cong \M_{\R\times -Y}
(\hat \beta, \hat \alpha).
\]
 Hence from the construction of the Seiberg-Witten-Floer homology, 
we see that $C_{-*}(-Y, -\s)$ can be identified  to   the
dual complex of $C_*(Y, \s)$. This 
gives rise to a natural pairing 
\ba 
\la \ , \ \ra : \qquad  HF^{SW}_*(Y, \s) \times HF^{SW}_{-*}(-Y, -\s) 
\longrightarrow \Z,
\label{pair:SWF}
\na
which satisfies the following properties whose
proofs are immediate from the definition:
\begin{enumerate}
\item For two relative degree 2 cycles $\Xi_1$ and $\Xi_2$ in
$HF^{SW}_*(Y, \s) $ and $HF^{SW}_{-*}(-Y, -\s)$ respectively
$< U (\Xi_1), \Xi_2 > = <\Xi_1, U( \Xi_2)>$. 
\item For two relative degree 1 cycles $\Xi_1$ and $\Xi_2$ in
$HF^{SW}_*(Y, \s) $ and $HF^{SW}_{-*}(-Y, -\s)$ respectively, and
any $[\pi]\in H_1(Y, \Z)/Torsion \cong H_1(-Y, \Z)/Torsion$,
$<[\pi](\Xi_1), \Xi_2 > = <\Xi_1, [\pi](\Xi_2)>$. 
\end{enumerate}
Thus, $<z.\Xi_1, \Xi_2> = <\Xi_1, z.\Xi_2>$
 for any $z\in \AA (Y)\cong \AA(-Y)$, for any cycles
$\Xi_1 \in HF^{SW}_*(Y, \s)$ and $\Xi_2\in HF^{SW}_{-*}(-Y, -\s)$
respectively.

\subsection{Variants of the Seiberg-Witten-Floer homology}

The Seiberg-Witten-Floer homology $HF_*^{SW}(Y, \s)$, 
as defined in the previous subsection,
is a finitely generated $\Z_{d(\s)}$-graded Abelian group. 
There are various ways of defining $\Z$-graded Seiberg-Witten-Floer
homology groups
for any closed oriented 3-manifold $Y$ with a non-torsion $\spinc$ structure
$\s$. These were briefly discussed at the end of section 4.1 of \cite{MW}.
In this subsection, we propose a way of defining $\Z$-graded 
Seiberg-Witten-Floer homology, which will be convenient for our later
description of
the relative Seiberg-Witten invariant for a 4-manifold with cylindrical end
modelled on $([-2, \infty)\times Y, \s)$.

Notice that the first Chern class of the non-torsion $\spinc$ structure
$\s$ defines a homomorphism:
$c_1(\s): \quad H^1(Y, \Z) \longrightarrow \Z$
by $c_1(\s) ([u]) = < [u]\cup c_1(\s), [Y]>$ for
any element $[u] \in H^1(Y, \Z)$. 

For any subgroup $K\subseteq Ker (c_1(\s))$, there is a subgroup $\G^K_Y$ of
the full gauge transformation group $\G_Y$, whose elements lie
in the connected components determined by $K$ (note that the group
of connected components of $\G_Y$ is $H^1(Y, \Z)$). 
Consider the Seiberg-Witten-Floer
homology theory for the Chern-Simons-Dirac functional on the
configuration space $\A_Y$ modulo the gauge group $\G^K_Y$. The
critical point set, denoted by $\M_{Y, K}(\s, \eta)$, is a covering space 
\[
\pi_K: \qquad \M_{Y, K}(\s, \eta) \longrightarrow \M_Y (\s, \eta)
\]
whose fiber is an $H^1(Y, \Z)/K$-homogeneous space, hence, there is 
a natural action of $H^1(Y, \Z)/K$ on $\M_{Y, K}(\s, \eta)$. 

The generators
of this variant of the Seiberg-Witten-Floer chain complex are elements
in $\M_{Y, K}(\s, \eta)$, with relative $\Z$-graded indices. 
Denote this chain complex by
\[
C^{SW}_{*, [K]}(Y, \s) = \sum_{\Gamma\in \M_{Y, K}(\s, \eta)} \Z <\Gamma>.
\]
The boundary operator $\partial^K$ is given by counting the gradient flowlines
of the perturted CSD functional on $\A_Y/\G^K_Y$ connecting two critical points
with relative index 1. Fix an element $\Gamma_\alpha \in \pi_K^{-1}(\alpha)$,
for any critical point $\beta \in \M_Y (\s, \eta)$, there 
is an orientation presevering diffeomorphism
\[
\M (\alpha, \beta) \cong 
\bigcup_{\Gamma_\beta \in \pi_K^{-1}(\beta)}
\M (\Gamma_\alpha, \Gamma_\beta),
\]
where $\M (\Gamma_\alpha, \Gamma_\beta)$ is the moduli space of flowlines
on  $\A_Y/\G^K_Y$ connecting $\Gamma_\alpha$ and $\Gamma_\beta$.
Hence,   it is easy to see that there
is a well-defined topological invariant, given by the
homology of this new version of the Seiberg-Witten-Floer chain complex,
though not finitely generated, still admitting an action of
$\AA (Y)$.

\begin{Def} For any subgroup $K\subseteq Ker (c_1(\s))$, there is
a variant of the Seiberg-Witten-Floer chain complex, whose generators
are elements from the covering space $\M_{Y, K}(\s, \eta)$ of 
$\M_Y (\s, \eta)$, and the boundary operator is given by 
counting the gradient flowlines
of the perturbed CSD functional on $\A_Y/\G^K_Y$ connecting two critical points
with relative index 1.  
 We denote the Seiberg-Witten-Floer homology in this setting
by $HF^{SW}_{*, [K]}(Y, \s)$, the actions of  elements in
$H_0(Y, \Z)$ and $H_1(Y, \Z)/Torsion$ decreasing degree
in $HF^{SW}_{*, [K]}(Y, \s)$ by $2$ and $1$ respectively. 
\end{Def}

From the definition, it is easy to get the following
periodicity property for $HF^{SW}_{*, [Ker(c_1(\s))]}(Y, \s)$:
\ba\label{isomorphism}
HF^{SW}_{m, [Ker(c_1(\s))]}(Y, \s) \cong HF^{SW}_{m \ \mod (d(\s))} (Y, \s),
\na
for any $m\in \Z$.

\begin{Rem}\label{u-action}
From the natural action of $H_1(Y, \Z)/K$ on the generators of
the chain complex $C^{SW}_{*, [K]}(Y, \s)$,
$HF^{SW}_{*, [K]}(Y, \s)$ can be thought as an $H_1(Y, \Z)/K$-equivariant
Seiberg-Witten-Floer homology of $(Y, \s)$. For any $[u]\in H_1(Y, \Z)/K$,
the action of $[u]$ on $HF^{SW}_{*, [K]}(Y, \s)$ is the following
$\AA(Y)$-equivariant homomorphism:
\[
[u]:\qquad HF^{SW}_{*, [K]}(Y, \s)
\longrightarrow HF^{SW}_{*-n, [K]}(Y, \s)
\]
with $n= <[u]\wedge c_1(\s), [Y]>.$
\end{Rem}

For $(-Y, -\s)$, where $-Y$ is $Y$ with the reversed orientation and
 $-\s$ is the induced $\spinc$ structure,
the corresponding Seiberg-Witten-Floer complex $C^{SW}_{*, [K]}(-Y, -\s)$
is the dual complex of $C^{SW}_{*, [K]}(Y, \s)$. This gives a natural pairing
\[
\la \ , \ \ra : \qquad  HF^{SW}_{*, [K]}(Y, \s) \times
HF^{SW}_{-*, [K]}(-Y, -\s)\longrightarrow \Z
\]
satisfying $<z.\Xi_1, \Xi_2> = <\Xi_1, z.\Xi_2>$
 for any $z\in \AA (Y)\cong \AA(-Y)$ and any cycles
$\Xi_1 \in HF^{SW}_{*, [K]}(Y, \s) $ and 
$\Xi_2\in HF^{SW}_{-*, [K]}(-Y, -\s)$ respectively.

 For any subgroup $K$ in $Ker(c_1(\s))$, there is a chain map
$C^{SW}_{*, [K]}(Y, \s) \to C^{SW}_*(Y, \s)$ induced from the 
map $\pi_K$, which descends to an $\AA (Y)$-equivariant homomorphism
\ba\label{SWF}
\pi_K:  \qquad  HF^{SW}_{*, [K]} (Y, \s) \longrightarrow
HF^{SW}_* (Y, \s).
\na

Let $K_1 \subset K_2$ be two subgroups in
$Ker(c_1(\s))$. There is a covering map 
\[
\pi: \qquad  \M_{Y, K_1}(\s, \eta)   \longrightarrow \M_{Y, K_2}(\s, \eta)
\]
which induces  a map:
\[
\pi_*: \qquad C^{SW}_{*, [K_1]} (Y, \s) \longrightarrow  
C^{SW}_{*, [K_2]} (Y, \s)
\]
assigning the generator $<\Gamma>$ of $C_{*, [K_1]} (Y, \s)$ to $<\pi(\Gamma)>$
in $C_{*, [K_2]} (Y, \s)$. 
For any $\alpha, \beta$ in $\M_{Y, K_2}(\s, \eta)$ of relative
index 1, and for any fixed $\Gamma_\alpha$ in $\pi^{-1}(\alpha)$, there 
is an orientation preserving diffeomorphism:
\[
\M_{\R\times Y} (\alpha, \beta) \cong 
\bigcup_{\Gamma_\beta\in \pi^{-1}(\beta)} 
\M_{\R\times Y} (\Gamma_\alpha, \Gamma_\beta).
\]
which implies that $\pi_*$ is a chain map and 
induces an  $\AA (Y)$-equivariant homomorphism 
\ba
\pi:  \qquad  HF^{SW}_{*, [K_1]} (Y, \s) \longrightarrow  
HF^{SW}_{*, [K_2]} (Y, \s).
\label{equ:homo}
\na

For any two subgroups $K_1$ and  $K_2$  in
$Ker(c_1(\s))$, we have the following commutative diagram which relates
various covering spaces of $\M_Y(\s, \eta)$:
\ba
\begin{array}{ccccc}
&& M_{Y, K_1}(\s, \eta)&& \\[-2mm]
& \stackrel{\pi^1\quad}{\nearrow} && \stackrel{\quad\pi_1}{\searrow}&
\\[-2mm]
 M_{Y, K_1\cap K_2}(\s, \eta) && \stackrel{\pi_0}{\longrightarrow} &&
\M_Y(\s, \eta)\\[-2mm]
& \stackrel{\quad\pi^2}{\searrow} && \stackrel{\pi_2\quad}
{\nearrow}&\\[-2mm]
&&\M_{Y, K_2}(\s, \eta)&& 
\end{array}\label{diag}
\na
where the fibers of $\pi_0, \pi^1, \pi_1, \pi^2$ and $\pi_2$
are the homogeneous spaces of $H^1(Y, \Z)/(K_1\cap K_2)$,
$K_1/(K_1\cap K_2)$, $H^1(Y, \Z)/K_1$, $K_2/(K_1\cap K_2)$ and
$H^1(Y, \Z)/K_2$ respectively. 

Fix an element in $\M_{Y, K_1\cap K_2}(\s, \eta)$, then there is
a well-defined $\Z$-graded index on $\M_{Y, K_1\cap K_2}(\s, \eta)$,
as $K_1$ and $K_2$ are two subgroups of $Ker (c_1(\s))$. Hence
there are the induced $\Z$-graded indices on
$\M_{Y, K_1}(\s, \eta)$ and $\M_{Y, K_2}(\s, \eta)$.
 The $\A(Y)$-equivariant homomorphisms 
as provided in (\ref{equ:homo}) and  
(\ref{SWF}) induce the following commutative
diagram of $\AA(Y)$-equivariant homomorphisms.
\ba
\begin{array}{ccccc}
&&HF^{SW}_{*, [K_1]}(Y, \s) && \\[-2mm]
& \stackrel{\pi^1\quad}{\nearrow} && \stackrel{\quad\pi_1}{\searrow}&
\\[-2mm]
HF^{SW}_{*, [K_1\cap K_2]}(\s, \eta) && \stackrel{\pi_0}{\longrightarrow} &&
HF^{SW}_*(\s, \eta) \\[-2mm]
& \stackrel{\quad\pi^2}{\searrow} && \stackrel{\pi_2\quad}
{\nearrow}&\\[-2mm]
&&HF^{SW}_{*, [K_2]}(\s, \eta) &&\end{array}\label{diag:1}
\na

\section{Gluing formulae for 4-dimensional Seiberg-Witten invariants}

In this Section, we will study the Seiberg-Witten invariant of a closed
smooth 4-manifold $X$ with $ b_2^+ (X)>1$ and 
$\spinc$ structure $\s$ by stretching
its metric along a smooth embedded separating 3-manifold $Y$ with the
condition that $\s|_Y$ is non-torsion.  Consider a 1-parameter family
of metrics $\{g_R\}_{R>0}$ on
$X$ such that for each $(X, g_R)$, there is an isometrical embedding
\ba\label{isometry}
([-R-2, R+2]\times Y, dt^2 + g_Y) \longrightarrow 
(X, g_R).
\na
Write for brevity, $X(R) = (X, g_R)$.

As $Y$ separates $X = X(0)$, we write $X(0) = X_+\cup X_-$ with
$Y = \{0\}\times Y = X_+\cap  X_-$. There is a canonical orientation
preserving isometry
\[
\Theta_R: \qquad X(R) - \{0\}\times Y \longrightarrow X_+(R) \cup X_-(R),
\]
where $X_+(R) = X_+ \cup [0, R] \times Y, X_-(R) = X_- \cup [-R, 0]\times Y$.
As $R\to \infty$, $X_\pm (R)$ become 4-manifolds with cylindrical ends
modelled on $[-2, \infty) \times Y$ and $(-\infty, 2]\times Y$.
We denote these 4-manifolds by $X_\pm (\infty)$.

Now, it is generally believed that, in some favorable cases, the moduli space 
of the Seiberg-Witten monopoles on $X(R)$, for a sufficiently large $R$, can be
obtained by gluing together the moduli spaces for $X_\pm(\infty)$ over
the ends. In this section, we will establish this gluing theorem and
give an expression for the Seiberg-Witten invariant of $X$ in terms of certain
relative Seiberg-Witten invariants for $X_\pm(\infty)$.

First we need to study the moduli space of Seiberg-Witten monopoles on 
$X_\pm(\infty)$ with $\spinc$ structures $\s_\pm$ which
are the pull-back $\spinc$ structures when restricted to the
ends. We assume that $\s_\pm|_Y = \pm \t$ is a non-torsion $\spinc$
structure. Then we will define relative Seiberg-Witten invariants
for $(X_\pm(\infty), \s_\pm)$ which take values in certain variants
of the Seiberg-Witten-Floer homology groups $HF_*^{SW}(Y, \t)$ and
$HF_*^{SW}(-Y, -\t)$ respectively.  As an application 
of the Seiberg-Witten-Floer homology, we will prove the gluing 
formula for  the Seiberg-Witten invariants on $X$ as claimed in Theorem 
\ref{Theorem:2}.

\subsection{Moduli space for  4-manifolds with cylindrical ends} 

 Let $X_+(\infty)$ be a Riemannian 4-manifold with a cylindrical 
end modelled on $[-2, \infty) \times Y$ and endowed with a $\spinc$ structure
$\s_+$ such that $\s_+|_{[-2, \infty) \times Y}$ is the pull-back $\spinc$
structure of a non-torsion $\spinc$
structure $\t$ on $Y$. Associated to $\s_+$, there is a pair of rank 2
Hermitian bundles $W^\pm$ on $X_+(\infty)$, a complex line bundle
$det(\s) = det (W^\pm)$, and a Clifford multiplication endomorphism
\[
T^*X_+(\infty) \otimes W^\pm \longrightarrow W^\mp.
\]
Using the Clifford multiplication of $dt$ over $[-2, \infty) \times Y$,
$W^\pm$ can be identified with a rank $2$ Hermitian bundle $W$ on $Y$, and
there is an induced Clifford multiplication homomorphism
\[
\rho: \qquad T^*Y \times W \longrightarrow W.
\]
Then $\t = (W, \rho)$ is the $\spinc$ structure on $Y$ as in the Definition
\ref{spinc:3d}.

The perturbed Seiberg-Witten equations on $(X_+(\infty), \s_+)$ are the
following equations, for a pair consisting of an $L^2_{2, loc}$-connection $\AA$
on $det(\s)$ and an $L^2_{2, loc}$-section $\Psi$ in $W^+$,  given by
\ba
\left\{
\begin{array}{l}
F^+_\AA = q (\Psi, \Psi) + \omega + \chi (\eta^+ + \p^+_\AA)\\[2mm]
\Dirac_\AA \Psi = \chi \p_\Psi.
\end{array}
\right.\label{SW:cyl}
\na
Here $\omega$ is an imaginary valued self-dual 2-form supported in a
non-empty open set in $X_+(-2) = X_+(0) - [-2, 0]\times Y$,
$\chi$ is a cut-off function on $X_+(\infty)$ which is
1 on $[0, \infty)\times Y$ and 0 on $X_+(-2)$, and 
$(\eta^+ + \p^+_\AA, \p_\Psi)$ is the perturbation of the
Seiberg-Witten equations on $[-2, \infty)\times Y$ as written
in (\ref{SW:4dpert}) (where $\eta$ and $\p$ satisfy Remark (\ref{condition:eta})
and Condition  (\ref{condition}) in their Baire sets). The perturbed
Seiberg-Witten equations (\ref{SW:cyl})  are invariant under the
action of the
gauge group, which is the group of $L^2_{3, loc}$-maps from $X_+(\infty)$
to $U(1)$.

Let $(\AA, \Psi)$ be a solution to the 
perturbed Seiberg-Witten equations (\ref{SW:cyl}) on $(X_+(\infty), \s_+)$,
in temporal gauge when restricted to the end $[-2, \infty) \times Y$. The
corresponding path, which is a gradient flowline of ${\cal C}_{\eta, \p}$,
is denoted by $\gamma (t) = (A(t), \psi(t)): [-2, \infty) \to \A_Y$.
We say that $(\AA, \Psi)$ has finite energy if
\[
\lim_{t\to\infty} \bigl( {\cal C}_{\eta, \p}(\gamma (-2)) - 
{\cal C}_{\eta, \p}(\gamma (t))\bigr) 
\]
exists and is finite. The finite energy condition is equivalent to
the following condition (see Lemma \ref{energy}):
\ba\label{finite:energy:X+}
\displaystyle{
\int_{-2}^\infty \bigl( \|\frac {\partial A}{\partial t}\|^2_{L^2(Y)}
+ \|\frac {\partial \psi}{\partial t}\|^2_{L^2(Y)} \bigr) dt }< \infty.
\na 

Any solution to the 
perturbed Seiberg-Witten equations (\ref{SW:cyl}) on $(X_+(\infty), \s_+)$
is smooth. We put a topology on the set of solutions 
by its embedding to 
\[
L^2_{2, loc}(\Omega^1(X_+(\infty), i\R) \oplus \Gamma (W^+)) \oplus \R
\]
where the last component is given by 
$\displaystyle{
\int_{-2}^\infty \bigl( \|\frac {\partial A}{\partial t}\|^2_{L^2(Y)}
+ \|\frac {\partial \psi}{\partial t}\|^2_{L^2(Y)} \bigr) dt}$
for the associated path $(A(t), \psi(t))$ determined by restricting
a solution $(\AA, \Psi)$ to $[-2, \infty)\times Y$ (which is 
in temporal gauge). Under this topology, 
the action of the $L^2_{3, loc}$ gauge group
is smooth on the set of finite energy solutions. 

The moduli space
of 4-dimensional monopoles on $(X_+(\infty), \s_+)$ is defined to be the
quotient of the set of finite energy solutions to (\ref{SW:cyl})
by the action of the $L^2_{3, loc}$ gauge group with the quotient topology.
Denote by $\M_{X_+}(\s_+)$ the moduli space of monopoles on 
$(X_+(\infty), \s_+)$.

By Lemma \ref{energy} and Lemma \ref{CSD:p}, we know that any finite energy
solutions to (\ref{SW:cyl}) decay exponentially to a 3-dimensional 
Seiberg-Witten monopole representing a critical point in $\M_Y( \s, \eta)$.
The proof of the following proposition is straightforward.

\begin{Pro} \label{decay:cyl}
Let $(X_+(\infty), \s^+)$ be a complete Riemannian 4-manifold with a 
cylindrical end $[-2, \infty)\times Y$ and a $\spinc$ structure $\s$
such that $\s|_{[-2, \infty)\times Y}$ is the pull-back $\spinc$ structure
of a non-torsion $\spinc$ structure $\t$ on $Y$.  
Let $(\AA, \Psi)$ be a finite energy solution to the
perturbed Seiberg-Witten equations (\ref{SW:cyl})) on $(X_+(\infty), \s_+)$,
in temporal gauge when restricted to the end $[-2, \infty) \times Y$ and
denoted by $(A(t), \psi(t))$. Then there is a constant $\delta >0$, and
a constant $C>0$, such that there is a 3-dimensional monopole 
$(A_\alpha, \psi_\alpha)$ representing a critical point 
$\alpha \in \M_Y(\t, \eta)$ for which
\[
\sum_{0\leq k\leq 2} \Bigl( |\nabla ^k(A(t)-A_\alpha)| 
+ |(\nabla _{A_\alpha})^k (\psi (t) - \psi_\alpha)|\Bigr) \leq
C e^{-\delta t}
\]
for any point $(t, y) \in [0, \infty)\times Y$.
\end{Pro}

This proposition  ensures that there is a boundary asymptotic limit map
\ba
\partial_\infty: \qquad \M_{X_+} (\s_+) \longrightarrow 
\M_{Y, X_+} (\t, \eta)
\label{boundary:limit}
\na
where $\M_{Y, X_+} (\t, \eta)$ is the quotient of solutions to the
perturbed Seiberg-Witten equations (\ref{SW:3d}) on $(Y, \t, \eta)$
by the action of those gauge transformations on $Y$ which can be extended to
$X_+(0)$. Note that $\M_{Y, X_+} (\t, \eta)$ is a covering space
of $\M_Y (\t, \eta)$, whose fiber is an $H^1(Y, \Z)/Im(i^*_+)$-homogeneous
space, where $Im(i^*_+)$ is the range of $i_+^*:
H^1(X_+, \Z)\to H^1(Y, \Z)$. Notice that $Im(i^*_+)$ is a subgroup
of $Ker (c_1(\t))$ where $c_1(\t)$ is the map from 
$H^1(Y, \Z)\to \Z$ defined in Subsection 3.5.

Over the end $[-2, \infty) \times Y$,
if $(X_+, \s_+)$ is modelled on $(Y, \t= (W, \rho))$ up to an isomorphism $u\in
C^\infty (Y, U(1))$, the equivalent classes of such modelling
is given by $[u] \in H^1(Y, \Z)/Im (i_+^*)$ determined by the connected
component of $C^\infty (Y, U(1))$ which $u$ belongs to. 
Then the corresponding 
asymptotic value map  is given by $[u]\circ \partial_\infty$
with $[u] \in H^1(Y, \Z)/Im (i_+^*)$ acting on $\M_{Y, X_+} (\t, \eta)$.

Denote by $\pi_+$ the covering map $\M_{Y, X_+} (\t, \eta)\to \M_Y (\t, \eta)$.
Let  $\M_{X_+} (\s_+, \alpha)$ be the subspace of 
$\M_{X_+} (\s_+)$ whose elements decay in the $L^2_{2,\delta}$-topology to
a 3-dimensional monopole representing $\alpha$ in $\M_Y (\t, \eta)$. 
Then we know that
\[
\M_{X_+} (\s_+, \alpha) = \bigcup_{\Gamma_\alpha \in \pi^{-1}_+(\alpha)}
\M_{X_+} (\s_+, \Gamma_\alpha)
\]
where $\M_{X_+} (\s_+, \Gamma_\alpha) = \partial_\infty^{-1} (\Gamma_\alpha)$
is the fiber of the boundary limit map (\ref{boundary:limit})
over $\Gamma_\alpha$. Then Proposition \ref{decay:cyl} tells us that
$\M_{X_+} (\s_+, \Gamma_\alpha)$
is the moduli space of solutions to the perturbed Seiberg-Witten equations
 (\ref{SW:cyl}) with $(\AA, \Psi)$ belonging to
\ba
(\AA_\alpha, \Psi_\alpha) + 
L^2_{2, \delta} (\Omega^1(X_+(\infty), i\R) \oplus \Gamma (W^+))
\label{config:delta}
\na
where $(\AA_\alpha, \Psi_\alpha)$ is a smooth pair consisting of 
a $U(1)$-connection
on $det(\s)$  and a section of $W^+$ on $X_+(\infty)$, 
which extends $(A_\alpha, \psi_\alpha)$
on $[-2, \infty)\times Y$ (as a constant pair representing $\Gamma_\alpha
\in \M_{Y, X_+}(Y, \s_+)$). The corresponding gauge group is the 
$L^2_{3, \delta}$-gauge transformation group. 

The following proposition 
describes the structure of $\M_{X_+} (\s_+, \Gamma_\alpha)$.
It is
a smooth oriented manifold with dimension given by the Fredholm index
of an operator ${\cal D}_x$ labelled by an element 
$x \in \M_{X_+} (\s_+, \Gamma_\alpha)$ representing  a solution
$(\AA, \Psi)$ to the  perturbed Seiberg-Witten equations
 (\ref{SW:cyl}). The definition of 
$${\cal D}_x: L^2_{2, \delta} \bigl(\Omega^1(X_+(\infty), i\R) 
\oplus \Gamma (W^+)\bigr)
 \longrightarrow  L^2_{1, \delta} \bigl(\Omega^0(X_+(\infty), i\R) 
\oplus \Omega^{2,+} (X_+(\infty), i\R) \oplus \Gamma (W^-)\bigr)$$
is
\ba
\begin{array}{c}
 {\cal D}_x (b, \Phi) =
\left\{ \begin{array}{l}
G^*_{(\AA, \Psi)} (b, \Phi)\\[1mm]
d^+ b -2 q (\Psi, \Phi) - \p^+_\AA(b, \Phi)\\[1mm]
\Dirac _\AA \Phi + \frac 12 b.\Psi + \p_\Psi (b, \Phi).
\end{array}\right.
\end{array}\label{linearization:cyl}
\na
where $G^*_{(\AA, \Psi)}$ is the adjoint of the linearization of the
gauge group action in $L^2_{\delta}$-norm, 
the terms involving $\p^+_\AA$ and $\p_\Psi $
come from the linearization of the perturbation term in 
(\ref{SW:cyl}) supported in $[-2, \infty )\times Y$. Note that
${\cal D}_x $ is indeed a Fredholm operator, whose  index is given by
the Atiyah-Patodi-Singer index formula:
\ba
\ind_{X_+} (\Gamma_\alpha) = 
\displaystyle{
\frac 14 \bigl(\frac{-1}{4\pi^2} \int_{X_+(\infty)} F_\AA\wedge
 F_\AA - 2 \chi(X_+) -3\sigma (X_+) \bigr)
+\frac 12 \rho (\Gamma_\alpha)},
\label{index:cyl}
\na
where $\rho (\Gamma_\alpha)$ is the rho-invariant of 
the boundary operator ${\cal T}_{(A_\alpha, \psi_\alpha)}$ (Cf. the
extended Hessian operator (\ref{complex:3d1})). Note that
$\frac{-1}{4\pi^2} \int_{X_+(\infty)} F_\AA\wedge
 F_\AA$ is finite for any finite energy solution 
(Cf. Corollary \ref{F:finite})), and $\sigma (X_+) $ is the signature of
the intersection form on $H^2(X_+(0), Y; \R)$ and
$\chi(X_+)$ is $2-2b_1+b_2$ with $b_1= dim H^1 (X_+(0), Y; \R)$ and
$b_2 = dim H^2(X_+(0), Y; \R)$. 

\begin{Pro} \label{smooth:cyl}
Fix an open set $U \subset X_+(-2)$ away from the boundary of $X_+(-2)$.
There is a Baire set of smooth imaginary-valued 2-forms with compact
support in $U$ such that for any $\Gamma_\alpha
\in \M_{Y, X_+} (\t, \eta)$: if 
$\ind_{X_+} (\Gamma_\alpha)< 0$, then $\M_{X_+}(\s_+, \Gamma_\alpha)$
is empty while if $\ind_{X_+} (\Gamma_\alpha)\geq 0$ and 
$\M_{X_+}(\s_+, \Gamma_\alpha)$ is non-empty, 
then $\M_{X_+}(\s_+, \Gamma_\alpha)$ is a smooth manifold of dimension
$\ind_{X_+} (\Gamma_\alpha)$ oriented with a choice of
orientation for the line $\Lambda^{top} H^1(X_+(0), Y; \R)
\otimes \Lambda^{top} H^{2, +}(X_+(0), Y; \R)$. Moreover, the cokernel
of ${\cal D}_x$ is trivial for any $x \in \M_{X_+}(\s_+, \Gamma_\alpha)$.
\end{Pro}
The proof of this proposition is fairly standard, we omit the details here,
see Proposition 2.14 \cite{MW}, Proposition 8.5 and Corollary 9.2 in \cite{MST},
Proposition 2.2 and Proposition 2.3 \cite{Tau3}. 
 
\begin{Rem} From the expression for the
dimension formula (\ref{index:cyl}), we know
that if $\Gamma_\alpha$ and $\Gamma_\alpha'$ are two different points in 
$\pi^{-1} (\alpha )$, there  is a gauge transformation $u$ on $det(\t)$ 
whose associated cohomology class $[u]$ is non-zero in
$H^1(Y, \Z)/Im (i^*_+)$, such that $u(\Gamma_\alpha)= \Gamma_\alpha'$
in $\pi^{-1}_+(\alpha)$. Then 
\[
\ind_{X_+} (\Gamma_\alpha) - \ind_{X_+} (\Gamma_\alpha') 
= < [u]\wedge c_1(\t), [Y]>.
\]
Note that if $[u] \in Ker (c_1(\t))$, then $< [u]\wedge c_1(\t), [Y]>=0$. The
relative grading on $\M_Y(\t, \eta)$ defined in Definition \ref{index:rel} is
related to $\ind_{X_+} $ in the following way:
\[
\ind( \alpha, \beta) = \ind_{X_+} (\Gamma_\beta) - 
\ind_{X_+}( \Gamma_\alpha) (mod \ d(\t)),
\]
where $\Gamma_\alpha \in \pi^{-1}_+(\alpha)$ and
$\Gamma_\beta\in \pi^{-1}_+(\beta)$.
\end{Rem}

Under the smoothness condition in Proposition \ref{smooth:cyl}, 
we now study the compactness of various components of the moduli space
$\M_{X_+}(\s_+)$. It is convenient to assemble the components in
 $\M_{X_+}(\s_+)$ with fixed dimension $d\geq 0$ as follows
\ba\label{moduli:d}
\M_{X_+}^d(\s_+) = \bigcup_{\alpha\in \M_Y(\t, \eta)} 
\Bigl( \bigcup_{\Gamma_\alpha\in \pi^{-1}_+(\alpha), \ind_{X_+} (\Gamma_\alpha)
=d} \M_{X_+}(\s_+, \Gamma_\alpha)\Bigr).
\na
Then $\M_{X_+}(\s_+) = \cup_{d\geq 0} \M_{X_+}^d(\s_+)$, each
$\M_{X_+}^d(\s_+)$ is called a stratum of dimension $d$ in $\M_{X_+}(\s_+)$.

The next proposition establishes the required compactification for each
stratum of dimension $d$ in $\M_{X_+}(\s_+)$.

\begin{Pro} \label{compact:cyl}
For any $d\geq 0$, and each $ \alpha \in \M_Y(\t, \eta)$, there are
only finitely many $\Gamma_\alpha$ in $\pi_+^{-1} (\alpha)$ with non-empty
moduli space $\M_{X_+}(\s_+, \Gamma_\alpha)$ in $\M^d_{X_+}(\s_+)$. Moreover,
there exists a compactification of $\M^d_{X_+}(\s_+)$, denoted by
$\overline{\M^d_{X_+}(\s_+)}$, which is a smooth manifold with corners. The
codimension k boundary faces can be expressed in the following form
\[
\bigcup_{\{\alpha_1, \cdots, \alpha_k\}}
\Bigl( \bigcup_{\{\Gamma_{\alpha_1}, \cdots, \Gamma_{\alpha_k}\}}
\M_{X_+}(\s_+, \Gamma_{\alpha_0}) \times \hat \M_{\R\times Y}
(\Gamma_{\alpha_0}, \Gamma_{\alpha_1}) \times \cdots \times
\hat \M_{\R\times Y} (\Gamma_{\alpha_{k-1}}, \Gamma_{\alpha_{k}})\Bigr).
\]
Here $\{\alpha_1, \cdots, \alpha_k\} \subset \M_Y(\t, \eta)$ and
$\Gamma_{\alpha_i} \in \pi^{-1}_+ (\alpha_i)$ such
that $\M_{X_+}(\s_+, \Gamma_{\alpha_i})$
is non-empty for each $ 0\leq i \leq k$, 
\[
d= \ind_{X_+} (\Gamma_{\alpha_{k}}) > \ind_{X_+} (\Gamma_{\alpha_{k-1}})
> \cdots >\ind_{X_+} (\Gamma_{\alpha_1}) > \ind_{X_+} (\Gamma_{\alpha_0}) \geq
0,\]
and $\hat \M_{\R\times Y} (\Gamma_{\alpha_i}, \Gamma_{\alpha_{i-1}})$ 
is the moduli space of
unparametrized flowlines for ${\cal C}_{\eta, \p}$, connecting
$\Gamma_{\alpha_i}$ and $\Gamma_{\alpha_{i-1}}$ on the configuration space
$\A_Y$ modulo those gauge transformations which can be extended to 
$X_+$.
\end{Pro}
\begin{proof}
The proof is standard in Floer gauge theory, we only give a sketch here, 
following the techniques developed in section 4 of  \cite{MW}.

There is a constant $E$ depending only on the Riemannian metric 
on $X_+(\infty)$ such that for any finite energy solution 
$(\AA, \Psi)$ in temporal gauge when restricted to 
$[-2, \infty) \times Y$, representing a point in $\M_{X_+}^d(\s_+)$, 
we have 
$|\Psi(x)| \leq E$ for any $x\in X_+(\infty)$. This follows from the 
Weitzenb\"ock formula for the Dirac operator $\Dirac_\AA$ and the fact
that, over 
the end, $(\AA, \Psi)$ satisfies the asymptotic estimate of Proposition
\ref{decay:cyl}. 

With this pointwise uniform bound on $\Psi$ and 
 the asymptotic estimate in Proposition
\ref{decay:cyl}, there is a  uniform bound on the $L^2$-norm 
of $F^+_\AA$ so that
\[
\displaystyle{-\frac {1}{4\pi^2} \int_{X_+(\infty)} F_\AA \wedge
F_\AA = \frac {1}{4\pi^2} \int_{X_+(\infty)} 
(|F_\AA^+|^2 - |F_\AA^-|^2) dvol_{X_+}}
\]
is uniformly bounded. This implies that there is a uniform upper bound on
the $L^2$-norm of $F_\AA$. 

Now let $(\AA_i, \Psi_i)$ represent a sequence of elements in 
$\M_{X_+}^d(\s_+)$. Then the standard bootstrapping arguments in elliptic
regularity theory show that on any compact set of $X_+(\infty)$, there
exists a subsequence of $(\AA_i, \Psi_i)$ which, after suitable 
gauge transformations, converges in $C^\infty$-topology to a finite
energy solution to the perturbed
Seiberg-Witten equations (\ref{SW:cyl}) on the corresponding
compact set of $X_+(\infty)$. These procedures, as illustrated in Theorem 4.1
\cite{MW},
will produce an element in 
\ba\label{faces}
\bigcup_{\{\alpha_1, \cdots, \alpha_k\}}
\Bigl( \bigcup_{\{\Gamma_{\alpha_1}, \cdots, \Gamma_{\alpha_k}\}}
\M_{X_+}(\s_+, \Gamma_{\alpha_0}) \times \hat \M_{\R\times Y}
(\Gamma_{\alpha_0}, \Gamma_{\alpha_1}) \times \cdots \times
\hat \M_{\R\times Y} (\Gamma_{\alpha_{k-1}}, (\Gamma_{\alpha_{k}})\Bigr)
\na
with $\{\alpha_1, \cdots, \alpha_k\} \subset \M_Y(\t, \eta)$ and
$\Gamma_{\alpha_i} \in \pi^{-1}_+ (\alpha_i)$ such that 
for each $ 0\leq i \leq k$ 
\[
d= \ind_{X_+} (\Gamma_{\alpha_{k}}) > \ind_{X_+} (\Gamma_{\alpha_{k-1}})
> \cdots >\ind_{X_+} (\Gamma_{\alpha_1}) > \ind_{X_+} (\Gamma_{\alpha_0}) \geq
0.\]

An easy consequence of this convergence argument is that there
are only finitely many $\Gamma_\alpha \in \pi_+^{-1}(\alpha)$ 
with non-empty  $\M_{X_+}(\s_+, \Gamma_{\alpha})$ in 
$\M_{X_+}^d(\s_+, \Gamma_{\alpha})$.

Next we need to show that any element in (\ref{faces})
can glued together to produce an element in
$\M_{X_+}(\s_+, \Gamma_{\alpha_{k}})$. This is the multiple
gluing theorem  in Theorem 4.23 and Proposition 4.25 of \cite{MW}. This implies
that for any $\Gamma_{\alpha_i}$ appearing in  (\ref{faces}),
$\M_{X_+}(\s_+, \Gamma_{\alpha_i})$ is non-empty. Note that by the
generic choices of various perturbations, all the
moduli spaces in our discussion are smooth and the cokernels of
the corresponding Fredholm operators are trivial, hence there is
no obstruction in the multiple
gluing theorem in Theorem 4.9 and Proposition 4.25 of \cite{MW}.

The smooth structure with corners on the compactification can be
obtained with the help of Lemma 4.26 in \cite{MW}, we omit the details.
\end{proof}

We should remark that there is a family version of Proposition
\ref{smooth:cyl} and Proposition \ref{compact:cyl}, which will
be crucial in the proof of the topological invariance of the 
various relative invariants introduced in the next Subsection. 

For simplicity, we assume that the Riemannian
metric on $Y$, the perturbation $\eta$ for the non-degeneracy
of $\M_Y(\s, \eta)$ and the perturbation $\p$ for the smoothness
of the moduli spaces of flowlines for ${\cal C}_{\eta, \p}$ are fixed.

Let $g_0$ and $g_1$ be two Riemannian metrics on $X_+(\infty)$
which agree with the product metric $dt^2 + g_Y$ over the
end $[-2, \infty)\times Y$. Let $\omega _0$ and $\omega _1$
be two compactly supported self-dual forms, with respect to $g_0$ and
$g_1$ respectively, in the Baire set in Proposition \ref{smooth:cyl}
such that the claims in Proposition \ref{smooth:cyl} and
Proposition \ref{compact:cyl} hold for the moduli spaces
$\M_{X_+}(\s_+, g_0, \omega_0)$ and $\M_{X_+}(\s_+, g_1, \omega_1)$.

\begin{Pro}\label{family:cyl}
There is a continuous interpolating family of metrics $g_t$ 
which agree with the product metric $dt^2 + g_Y$ over the
end $[-2, \infty)\times Y$ and
compactly supported self-dual forms $\omega_t$ such that, for
any $\Gamma_\alpha \in \M_{Y, X_+}(\t, \eta)$ with
$\ind_{X_+}(\Gamma_\alpha) =d\geq 0$,
the family of moduli spaces 
\[
\M_d^F (\s_+, \Gamma_\alpha) = \bigcup_{t\in [0, t]}
 \M_{X_+}^d (\s_+, \Gamma_\alpha, g_t, \omega_t)
\]
satisfies the following properties:
\begin{enumerate}
\item $\M_d^F(\Gamma_\alpha)$ 
is a (d+1)-dimensional, oriented, smooth manifold and the
parameter projection $\M_d^F (\Gamma_\alpha) \to [0, 1]$ is smooth and admits
a product structure near the ends of $[0, 1]$. 
The orientation on $\M_d^F (\Gamma_\alpha)$ is compatible with 
the orientations of  $\M_{X_+}(\s_+,\Gamma_\alpha, g_0, \omega_0)$ and 
$\M_{X_+}(\s_+,\Gamma_\alpha, g_1, \omega_1)$ 
as the boundary. 
\item There is a compactification of $\M_d^F(\Gamma_\alpha)$ which is compatible
with the compactifications of 
$\M_{X_+}(\s_+,\Gamma_\alpha, g_0, \omega_0)$ and 
$\M_{X_+}(\s_+, \Gamma_\alpha, g_1, \omega_1)$ described in Proposition
\ref{compact:cyl}. In particular, the codimension one boundary of the 
compactification of $\M_d^F(\Gamma_\alpha)$ is given by 
\[
\begin{array}{c}
\M_{X_+}(\s_+, \Gamma_\alpha, g_0, \omega_0) \cup 
(-\M_{X_+}(\s_+,\Gamma_\alpha, g_1, \omega_1))
\\
\cup \bigcup_{t\in [0, 1]}
\Bigl( \bigcup_{\alpha_0\in \M_Y (\t, \eta)}
\bigcup_{\Gamma_{\alpha_0} \in \pi_+^{-1}(\alpha_0)}
 M_{X_+}(\s_+, \Gamma_{\alpha_0}, g_t, \omega_t) \times \hat \M_{\R\times Y}
(\Gamma_{\alpha_0}, \Gamma_{\alpha}) \Bigr). 
\end{array}
\]
Here $\Gamma_{\alpha_0} \in \pi_+^{-1}(\alpha_0)$ and 
$ 0\leq  \ind_{X_+} (\Gamma_{\alpha_0}) < d$.
\end{enumerate}
\end{Pro}
The proof of this proposition is just the family version of arguments in the
proofs of  Proposition \ref{smooth:cyl} and  Proposition
\ref{compact:cyl}. We will not repeat the arguments here and
leave the details to the reader.

\subsection{Relative invariants for 4-manifolds with boundary}
\label{Relative}

Let $X_+(0)$ be an oriented 4-manifold with one
boundary $Y$. We endow $X_+(0)$ with a complete metric, denoted
by $ X_+(\infty)$, so that $X_+(\infty)$ has a cylindrical
end modelled on  $[-2, \infty )\times Y$. 
Let $\s_+$ be a $\spinc$ structure
on $X_+(\infty)$ whose induced $\spinc$ structure on $Y$ is $\t$ such that 
$c_1(\t)$ is non-torsion. 
There is a subgroup $Im (i_+^*)$ of $Ker (c_1(\t))$ such that
the covering map $\M_{Y, X_+}(\t, \eta) \to \M_Y(\t, \eta)$ has
as its fiber an $H^1(Y, \Z)/Im (i_+^*)$-homogeneous space. 

Fix an orientation for the line $\Lambda^{top} H^1(X_+(0), Y; \R)
\otimes \Lambda^{top} H^{2, +}(X_+(0), Y; \R)$.  
As studied in the previous Section, the finite energy moduli space has a 
well-defined asymptotic limit map:
\[ \M_{X_+}(\s_+) \to \M_{Y, X_+}(\t, \eta) .\]
The structure of $\M_{X_+}(\s_+)$ was established in       
Proposition \ref{smooth:cyl} and Proposition \ref{compact:cyl}.

Denote by $\AA(X_+) = Sym^*(H_0(X_+)) \otimes \Lambda^*(H_1(X_1)/Torsion)$
the free graded algebra generated by the class
of the element in $H_0(X_+)$ and the 1-cycle $\gamma \in
H_1(X_+)$ with degrees 2 and 
1 respectively. Then the relative Seiberg-Witten 
invariant for $(X_+(\infty), \s_+)$ will be defined as a map
\[
SW_{X_+} (\s, \cdot): \qquad \AA (X_+) \mapsto HF^{SW}_{*, [Im(i^*_+)]}(Y, \t).
\]

For any monomial $z = U^k \gamma_1 \wedge \gamma_2 \wedge \cdots
\wedge \gamma_l$ in $\AA(X_+)$ 
with $2k+l =d\geq 0$, we need to consider all the 
components of dimension $d$, $\M^d_{X_+} (\s_+)$ in $\M_{X_+} (\s_+)$.  
Recall that  
\[
\M^d_{X_+} (\s_+) = 
\bigcup_{\alpha\in \M_Y(\t, \eta)}
\Bigl( \bigcup_{\Gamma_\alpha\in \pi^{-1}_+(\alpha), \ind_{X_+} (\Gamma_\alpha)
=d} \M_{X_+}(\s_+, \Gamma_\alpha)\Bigr).
\]

Choose $k$-points $\{x_1, \cdots, x_k\}$ in $X_+(\infty)$, and for each
$x_i$, choose a complex line $L_{x_i}$ in $W_{x_i}^+$ (the fiber
of $W^+$ over $x_i$). Also choose $l$-loops $\Xi =
\{\pi_1, \cdots, \pi_l\}$ representing 1-cycles $\{ \gamma_1, \cdots,
\gamma_l\}$. For a generic choice of
$\Lambda = \{ (x_1, L_{x_1}), \cdots, (x_k,  L_{x_k})\}$, as a
subset of $\M_{X_+}(\s_+, \Gamma_\alpha)$ (for $\Gamma_\alpha$ with
$ind_{X_+} (\Gamma_\alpha) =d$),
\[
\M_{X_+}^\Lambda (\s_+, \Gamma_\alpha)
= \{ [\AA, \Psi] \in \M_{X_+}(\s_+, \Gamma_\alpha) |
\Psi (x_i) \in L_{x_i}, \text{for each $x_i$}, 1\leq i \leq k\}
\]
is a $l$-dimensional, oriented and smooth
 submanifold of $\M_{X_+}(\s_+, \Gamma_\alpha)$. Moreover,
$M_{X_+}^\Lambda (\s_+, \Gamma_\alpha)$ meets with each boundary
face of the compactification $\overline{\M_{X_+}(\s_+, \Gamma_\alpha)}$
transversally. This is the consequence of Proposition \ref{smooth:cyl}
and Proposition \ref{compact:cyl}  and the standard application of
Sard theorem.

Similarly, for any $l$-loops $\Xi = \{\pi_1, \cdots, \pi_l\}$, 
there is a holonomy map
\[
h_\Xi: \qquad \M_{X_+}(\s_+, \Gamma_\alpha) \longrightarrow 
U(1)^l
\]
given by the holonomies of $[\AA, \Psi] \in \M_{X_+}(\s_+, \Gamma_\alpha)$
along loops in $\Xi$. Then for a generic point $q \in U(1)^l$, there is
a $2k$-dimensional, oriented and smooth submanifold
$\M_{X_+}^\Xi (\s_+, \Gamma_\alpha) = h_\Xi^{-1}(q)$ of 
$\M_{X_+}(\s_+, \Gamma_\alpha)$.
Moreover, $\M_{X_+}^\Xi (\s_+, \Gamma_\alpha)$ meets with each boundary
face of the compactification $\overline{\M_{X_+}(\s_+, \Gamma_\alpha)}$
transversally.

Then we obtain that
\[
 \M_{X_+} ^{\Lambda, \Xi} (\s_+, \Gamma_\alpha) 
= \M_{X_+} ^\Lambda(\s_+, \Gamma_\alpha) \cap \M_{X_+}^\Xi 
(\s_+, \Gamma_\alpha) 
\]
is a compact, 0-dimensional, oriented and smooth manifold.  
We write

\[
SW_{X_+}(\s_+, z, \Gamma_\alpha) = \# \bigl(
\M_{X_+} ^{\Lambda, \Xi} (\s_+, \Gamma_\alpha) \bigr),
\] 
an integer which is independent of the choices
of $\Lambda$ and $\Xi$, as the space of $\{\Lambda, \Xi\}$ is
path-connected.

For those $\Gamma_\alpha$ with $\ind_{X_+}(\Gamma_\alpha) \neq d$, we set
$ SW_{X_+}(\s_+, z, \Gamma_\alpha) =0$. Now we can define the following
chain element in $C_{d, [Im(i^*_+)]}^{SW}(Y, \t)$ (for $d =deg (z)$):
\ba
SW_{X_+}(\s_+, z) = 
\sum_{\alpha \in \M_Y(\t, \eta)} \sum_{\Gamma_\alpha} 
SW_{X_+}(\s_+, z, \Gamma_\alpha) <\Gamma_\alpha>
\label{relative:SW}
\na
which is a finite sum according to Proposition
\ref{compact:cyl}. 

\begin{Pro} $SW_{X_+}(\s_+, z)$ is a cycle in 
$C_{d, [Im(i^*_+)]}^{SW}(Y, \t)$
whose image in $HF^{SW}_{d, [Im(i^*_+)]}(Y, \t)$ is independent
of the metric and perturbation in the perturbed
Seiberg-Witten equations (\ref{SW:cyl}). 
\end{Pro}
\begin{proof}
To see that $SW_{X_+}(\s_+, z)$ is a cycle we need to show that
$\partial^{Im(i_+^*)}(SW_{X_+}(\s_+, z))$ is zero. Consider
 the $(d+1)$-dimensional
components, $\M_{X_+}^{d+1}(\s_+)$, of $\M_{X_+}(\s_+)$. Express
\[
M_{X_+}^{d+1}(\s_+) = \cup_{\beta\in \M_Y (\t, \eta)}
\bigcup_{\Gamma_\beta\in \pi^{-1}_+(\beta), \ind_{X_+}(\Gamma_\beta) = d+1}
\M_{X_+}(\s_+, \Gamma_\beta),
\]
whose compactification is described in Proposition \ref{compact:cyl}. Notice
that the codimension one boundary face of $\overline{M_{X_+}^{d+1}(\s_+)}$
consists of
\[
\cup_{\alpha, \beta \in \M_Y (\t, \eta)}
\Bigl( \bigcup_{\Gamma_\alpha, \Gamma_\beta}
\M_{X_+}(\s_+, \Gamma_\alpha) \times \hat\M_{\R\times Y} (\Gamma_\alpha,
\Gamma_\beta)\Bigr)
\]
where $\Gamma_\alpha\in \pi_+^{-1}(\alpha)$ and $\Gamma_\beta
\in \pi_+^{-1}(\beta)$ with $\ind_{X_+}(\Gamma_\beta) =d+1$ and
the unparametrized moduli space $\hat\M_{\R\times Y} (\Gamma_\alpha,
\Gamma_\beta)$ is non-empty. 

Choose $\Lambda$ and $\Xi$ away from the cylindrical end 
$[-2, \infty) \times Y$, then the submanifold 
\[\M_{X_+} ^{\Lambda, \Xi} (\s_+, \Gamma_\beta)
= \M_{X_+} ^\Lambda(\s_+, \Gamma_\beta)\cap \M_{X_+}^\Xi
(\s_+, \Gamma_\beta)
\]
is a smooth oriented 1-manifold whose
boundary consists of
\[
\cup_{\alpha\in \M_Y (\t, \eta)}
\bigcup_{\Gamma_\alpha \in \pi_+^{-1}(\alpha), \ind_{X_+}(\Gamma_\alpha)
= d} \M_{X_+} ^{\Lambda, \Xi} (\s_+, \Gamma_\alpha)
\times \hat\M_{\R\times Y} (\Gamma_\alpha,
\Gamma_\beta).
\]
The counting of boundary points with sign yields
\[
\sum_{\alpha\in \M_Y (\t, \eta)}
\sum_{\Gamma_\alpha \in \pi_+^{-1}(\alpha), \ind_{X_+}(\Gamma_\alpha)
= d}  SW_{X_+} (\s_+, z, \Gamma_\alpha) \#\bigl( 
\hat\M_{\R\times Y} (\Gamma_\alpha,
\Gamma_\beta)\bigr) =0.
\]
Note that $\#\bigl( 
\hat\M_{\R\times Y} (\Gamma_\alpha,
\Gamma_\beta)\bigr)$ is the coefficient of $\Gamma_\beta$ in
$\partial^{Im(i^*_+)} (\Gamma_\alpha)$, this shows
that $ \partial^{Im(i^*_+)} (SW_{X_+}(\s_+, z) ) =0$. 

To establish the independence of $[SW_{X_+}(\s_+, z)]$,
as an element in $HF^{SW}_{d, [Im(i^*_+)]} (Y, \t)$ on
the choice of metric and perturbation on $X_+(\infty)$, we need to study
the family version of moduli spaces interpolating between
$\M_{X_+}(\s_+, g_0, \omega_0)$ and $\M_{X_+}(\s_+, g_1, \omega_1)$.

Choose $\Lambda$ and $\Xi$ away from the cylindrical end
$[-2, \infty) \times Y$, then
Proposition \ref{family:cyl} implies that, 
\[\begin{array}{lll}
&&SW_{X_+} (\s_+, z, g_0, \omega_0) - SW_{X_+} (\s_+, z, g_1, \omega_1)
\\[2mm]
&=& \sum_{\alpha_0, \alpha\in \M_Y (\t, \eta)}
\sum_{\Gamma_{\alpha_0}, \Gamma_\alpha} 
H_{\Gamma_{\alpha_0}} n_{\Gamma_{\alpha_0}\Gamma_\alpha } <\Gamma_\alpha>,
\end{array}\]
where $\Gamma_{\alpha_0}  
\in \pi_+^{-1}(\alpha_0)$ with $\ind_{X_+} (\Gamma_{\alpha_0})
=d-1$ and $\Gamma_\alpha\in \pi_+^{-1}(\alpha)$ 
with $\ind_{X_+} (\Gamma_{\alpha})$, $n_{\Gamma_\alpha\Gamma_\beta} = \#\bigl( 
\hat\M_{\R\times Y} (\Gamma_{\alpha_0},
\Gamma_\alpha)\bigr)$,
 and $H_{\Gamma_{\alpha_0}}$ is the counting of the points
$([\AA, \Psi], t)$ in $\M^F_{d-1}(\Gamma_{\alpha_0}) =
\cup_{t\in [0, 1]} \M_{X_+}(\s_+, 
\Gamma_{\alpha_0}, g_t, \omega_t)$ satisfying
\begin{enumerate}
\item $[\AA, \Psi] \in \M_{X_+}(\s_+, \Gamma_{\alpha_0}, g_t, \omega_t)$ 
such that
$\Psi(x_i) \in L_{x_i}$ for each $(x_i, L_{x_i}) \in \Lambda$,
\item $h_{\Xi} ([\AA, \Psi]) =q$ for a generic $q \in U(1)^l$.
\end{enumerate} 
Proposition \ref{family:cyl} implies that $H_{\Gamma_{\alpha_0}}$
is well-defined, and 
\[
\sum_{\alpha_0, \alpha\in \M_Y (\t, \eta)} 
\sum_{\Gamma_{\alpha_0}, \Gamma_\alpha} 
H_{\Gamma_{\alpha_0}} 
n_{\Gamma_{\alpha_0}\Gamma_\alpha} <\Gamma_\alpha>\]
 is a boundary chain element,
hence, $SW_{X_+} (\s_+, z, g_0, \omega_0) $ and 
$ SW_{X_+} (\s_+, z, g_1, \omega_1)$
define the same element in $HF^{SW}_{*, [Im(i_+^*)]} (Y, \t)$.
\end{proof}

To summarize, we have defined a relative Seiberg-Witten invariant as the 
following linear functional:
\ba \label{relative:4d}\begin{array}{cccc}
SW_{X_+}( \s_+, \cdot): \qquad & 
\AA (X_+)& \longrightarrow & HF_{*, [Im (i_+^*)]}^{SW} (Y, \t)\\[2mm]
& z & \mapsto & [SW_{X_+}( \s_+, z)].
\end{array}
\na
$[SW_{X_+}( \s_+, z)]$ has an representative written as in (\ref{relative:SW}).

Notice that there is a homomorphism $(i_+)_*: \AA (Y) \to \AA(X)$,
which induces an action of $\AA (Y)$ on $\AA(X)$. Then it is evident
that from the definition of the action of $\AA (Y)$ on 
$HF_{*, [Im (i_+^*)]}^{SW} (Y, \t)$ and the definition of the relative
invariant, $SW_{X_+}( \s_+, \cdot)$ in (\ref{relative:4d}) is
$\AA (Y)$-equivariant, in the sense that, 
\[
[ SW_{X_+} (\s_+, (i_+)_*(\gamma) z)] = \gamma [SW_{X_+} (\s_+, z)] ,
\]
for any $\gamma\in \AA (Y)$ and $z\in \AA(X)$.

\begin{Rem} The definition of the  relative Seiberg-Witten invariant can be
generalized to the case that $X_+$ has several boundary components.
For example, if $\partial (X_+) = (-Y_1) \cup Y_2$, let $\s_+$ be a
$\spinc$ structure on $X_+$ whose restrictions to the boundary manifolds
are  non-torsion $\spinc$ structures $\t_1$ and $\t_2$ respectively. 
Let $i_1$ and
$i_2$ be the boundary embedding map, $Im (i_1^*)$ and $Im (i_2^*)$
be the range of the induced maps on $H^1(X_+, \Z)$. Endowing 
$X_+$ with a complete metric such that $X_+(\infty) $ is a cylindrical-end
manifold modelled on $(-\infty, 2] \times Y_1 \cup [2, \infty) \times Y_2$,
then we can define a relative Seiberg-Witten invariant for $(X_+, s_+)$:
\[
SW_{X_+} (\s_+, \cdot): \qquad \AA (X_+) \longrightarrow
HF^{SW}_{*, [Im (i_1^*)]}(-Y_1, -\t_1)\otimes 
HF^{SW}_{*, [Im (i_2^*)]}(Y_2, \t_2).
\]
In other words, for any $z\in \AA (X_+) $, $SW_{X_+} (\s_+, z)$
defines a homomorphism from $HF^{SW}_{*, [Im (i_1^*)]}(Y_1, \t_1)$
to $HF^{SW}_{*, [Im (i_2^*)]}(Y_2, \t_2)$ using the natural pairing
between $HF^{SW}_{*, [Im (i_1^*)]}(Y_1, \t_1)$ and
$HF^{SW}_{*, [Im (i_1^*)]}(-Y_1, -\t_1)$. 
\end{Rem}

\subsection{Gluing formulae for the Seiberg-Witten invariants}

Consider a closed 4-manifold $(X, \s)$  with $b^+_2 \geq 1$ and a $\spinc$
structure $\s$. Suppose that  $X$ has a decomposition $X_+ \cup_Y X_+$
 along a 3-dimensional submanifold $Y$ 
where $b_1(Y)>0$ and $\t =\s|_Y$ is non-torsion. 
In this subsection we will develop a gluing theorem for the
moduli space of the Seiberg-Witten monopoles on  $(X, \s)$ in terms of the the
moduli spaces for  the Seiberg-Witten monopoles on  $(X_\pm, \s_\pm)$
with $\s_\pm = \s|_{X_\pm}$. Let $i_\pm$ be the boundary embedding maps
from $X_\pm \to X$, and let $Im(i_\pm^*)$ be the range of 
the map $i_\pm^*: H^1(X_\pm, \Z) \to H^1(Y, \Z)$. 

The set of $\spinc$ structures on $X$ which agree with $\s_\pm$ when
restricted to $X_\pm$ is an affine space over $H^1(Y, \Z)/(Im(i_+^*)
+ Im(i_-^*))$. Denote this set of $\spinc$ structures on $X$ by
$\spinc (X, \s_\pm)$. Note that $\spinc (X, \s_\pm)$ can be
obtained by gluing the $\spinc$ structures $\s_\pm$ 
on $X_\pm$ along $Y$ using gauge transformations on $(Y, \t)$. The
equivalent classes of $\spinc$ structures, obtained in this way, are
classified by the connected components of $H^1(Y, \Z)/(Im(i_+^*)
+ Im(i_-^*))$. Denote by $\s_+\#_u\s_-$ the resulting
$\spinc$ structure on $X$ obtained 
by gluing $\s_\pm$ using the gauge transformation
$u$, then 
\[
\spinc (X, \s_\pm) =\{ \s_+\#_u\s_- | u\in \G_Y, [u] \in H^1(Y, \Z)/(Im(i_+^*)
+ Im(i_-^*))\}.
\]

Consider a family of metrics $\{g_R\}_{R>0}$ on
$X$ such that for each $X(R) = (X, g_R)$, there is an isometrically embedded
submanifold $([-R-2, R+2]\times Y, dt^2 + g_Y)$. There are two 4-manifolds
$X_\pm (R)$ obtained by $X(R) = X_+(R)\cup  X_-(R)$. As $R\to \infty$,
$X(R)$ has a geometric limit, which are  4-manifolds
with cylindrical ends, denoted by $X_\pm(\infty)$. 

In the previous Subsection, we defined relative
Seiberg-Witten invariants for $(X_\pm(\infty), \s_\pm)$ using
the finite energy moduli spaces of the perturbed Seiberg-Witten
equations on $(X_\pm(\infty), \s_\pm)$. These relative
Seiberg-Witten invariants are $\AA(Y)$-equivariant linear
functionals:
\[
SW_{X_\pm}(\s_\pm, \cdot): \qquad \AA(X_\pm)\longrightarrow 
HF^{SW}_{*, [Im (i_\pm^*)]} (\pm Y, \pm \t).
\]

Fix an orientation on the line $\Lambda^{top} H^1(X, \R)
\otimes \Lambda^{top} H^{2, +}(X, \R)$ which is induced
from tensoring the orientations on 
$\Lambda^{top} H^1(X_+(0), Y; \R)
\otimes \Lambda^{top} H^{2, +}(X_+(0), Y; \R)$ and
$\Lambda^{top} H^1(X_-(0), Y; \R)
\otimes \Lambda^{top} H^{2, +}(X_-(0), Y; \R)$. 
When $b^+_2 =1$, we also need to fix an  orientation on the line
$H^{2, +}(X, \R)$ such that $c_1(\s)\cdot \omega^+ >0$ for an oriented
generator $\omega^+$ of $H^{2, +}(X, \R)$.
Then there is a Seiberg-Witten invariant as defined by Taubes in \cite{Tau4},
which is a linear functional:
\[
SW_X(\s, \cdot): \qquad \AA(X)\longrightarrow \Z,
\]
where $\AA(X) =  Sym^*(H_0(X, \Z)) \otimes \Lambda^*(H_1(X, \Z)/Torsion)$
is the free graded algebra generated by the class
of elements in $H_0(X, \Z)$ and $H_1(X, \Z)$ 
with degree 2 and 1 respectively.

On $(X(R), \s)$ for $\s \in \spinc (X, \s_\pm)$, consider the following
perturbed Seiberg-Witten equations for a pair consisting of a $U(1)$-connection
$\AA$ on $det (\s)$ and a spinor section $\Psi$ in $W^+$:
\ba
\left\{
\begin{array}{l}
F^+_\AA = q (\Psi, \Psi) + \omega^+ + \omega^-  
+ \chi (\eta^+ + \p^+_\AA)\\[2mm]
\Dirac_\AA \Psi = \chi \p_\Psi
\end{array}
\right.\label{SW:4dX}
\na
 where we use the following notation.

\noindent$\bullet$ $\omega^\pm$ are  imaginary valued self-dual 2-forms
on $X_\pm (0)$ with compact support in a non-empty open set such that
Proposition \ref{smooth:cyl} and Proposition \ref{compact:cyl} hold
for the finite energy moduli spaces for $(X_\pm(\infty),\s_\pm)$. 

\noindent $\bullet$
$(\eta^+ + \p^+_\AA, \p_\Psi)$ is the perturbation of the
Seiberg-Witten equations on $[-2, \infty)\times Y$ as written
in (\ref{SW:4dpert}) (where $\eta$ and $\p$ satisfy Remark 
(\ref{condition:eta})and Condition  (\ref{condition})), and 

\noindent $\bullet$
$\chi$ is a cut-off function on $X(R)$ which is 1 on $[-R, R]\times Y$
and $0$ on $X_\pm(-2)$. 

The moduli space of the Seiberg-Witten monopoles on $(X(R), \s)$, denoted
by $\M_{X(R)}(\s)$, is the quotient of the set of smooth solutions to
the above perturbed Seiberg-Witten equations (\ref{SW:4dX}) by the
action of gauge group $C^\infty(X(R), U(1))$.
Then (Cf.\cite{KM1}\cite{Morgan})
for generic $\omega^\pm$, $\M_{X(R)}(\s)$, if non-empty,
is a compact, oriented and smooth manifold of dimension given by
\[
d_X(\s)= \frac 14 (c_1(\s)^2 - (2 \chi (X) + 3\sigma (X) )) \geq 0.
\]

Any solution $(\AA, \Psi)$ to the perturbed Seiberg-Witten 
equations (\ref{SW:4dX}), when restricted to 
$[-R, R]\times Y$, is gauge equivalent to a path of
gradient flowlines of ${\cal C}_{\eta, \p}$ on $\A_Y$. Denote
this associated path of flowlines by $\gamma (t) = (A(t), \psi(t))$, then
the next lemma ensures that the variation of ${\cal C}_{\eta, \p}$
along the associated path $\gamma (t)$ is uniformly bounded.

\begin{Lem} There is a constant $E$ independent of $R$ such that for any
solution $(\AA, \Psi)$ to the perturbed Seiberg-Witten 
equations (\ref{SW:4dX}) on $(X(R), \s)$, the variation of 
${\cal C}_{\eta, \p}$ along the associated path of gradient flowlines
$\gamma (t)= (A(t), \psi(t))$ is uniformly bounded by $E$, that is
\[
|{\cal C}_{\eta, \p}(\gamma (-R)) - {\cal C}_{\eta, \p}(\gamma (R)) | 
\leq E.
\]
\end{Lem}
\begin{proof} From Lemma \ref{CSD:p}, we only need to show that 
$|{\cal C}(\gamma (-R)) - {\cal C}(\gamma (R)) |$
is uniformly bounded. Note that
\[\begin{array}{lll}
0&\leq &{\cal C}(\gamma (-R)) - {\cal C}(\gamma (R))\\[2mm]
&=& \displaystyle{ \frac 12 \int_{[R_1, R_2]\times Y}}
F_{\mathbb A} \wedge F_{\mathbb A}
+ \displaystyle{\int_{\partial X_+(0)} <\psi, \dirac_A \psi> dvol_Y
- \int_{\partial X_-(0)} <\psi, \dirac_A \psi> dvol_Y}.
\end{array}
\]
A standard application of the Weitzenb\"ock formulae
implies that, there is a uniform pointwise bound for
$|\Psi|$. Hence there is a uniform pointwise bound for
$|F^+_\AA|$ independent of $R$ and $(\AA, \Psi)$, from which, with the aid
of the perturbed Seiberg-Witten equations, we obtain 
a uniform bound for 
\[
\displaystyle{\int_{\partial X_+(0)} <\psi, \dirac_A \psi> dvol_Y
\qquad \text{and}
\int_{\partial X_-(0)} <\psi, \dirac_A \psi> dvol_Y}.
\]
Hence, there is a constant $C>0$ such that
\ba
\displaystyle{ \frac 12 \int_{[-R, R]\times Y}}
F_{\mathbb A} \wedge F_{\mathbb A}
 \geq -C.
\label{finite:F}
\na
Notice that there is an identity
\ba\begin{array}{lll}
&&\displaystyle{ \frac{1}{4\pi^2}} <c_1(\s)^2, [X]>\\[2mm]
&=& - \displaystyle{ \int_{X(R)} }F_{\mathbb A} \wedge F_{\mathbb A}\\[2mm]
&=& \|F_\AA^+\|^2_{L^2(X_+(0)\cup X_-(0))} - 
 \|F_\AA^-\|^2_{L^2(X_+(0)\cup X_-(0))} -
\displaystyle{ \int_{[-R, R]\times Y}}
F_{\mathbb A} \wedge F_{\mathbb A}.
\end{array}\label{identity}
\na
{}From (\ref{finite:F}) and (\ref{identity}), we know that
$\|F_\AA^-\|^2_{L^2(X_+(0)\cup X_-(0))}$ is uniformly bounded, which
in turn implies that 
\[
| \displaystyle{ \frac 12 \int_{[-R, R]\times Y}}
F_{\mathbb A} \wedge F_{\mathbb A}|
\]
is also uniformly bounded. This completes the proof of the lemma.
\end{proof}

An immediate consequence of the above finite energy lemma concerning
any solution to the perturbed Seiberg-Witten equations (\ref{SW:4dX})
is the following geometric limit theorem as $R\to \infty$.

\begin{The}\label{geo:limit}
Let $(\AA_R, \Psi_R)$ be a solution to the perturbed 
Seiberg-Witten equations (\ref{SW:4dX})
on $(X(R), \s)$ for $\s\in \spinc (X, \s_\pm)$, then there is a finite
integer $k\geq 0$ such that as $R\to \infty$, there exists
a subsequence of $(\AA_R, \Psi_R)$ which converges in $L^2_{2,loc}$-topology
to 
\[
\{(\AA_+, \Psi_+), (\AA_1, \Psi_1), (\AA_2, \Psi_2),
\cdots (\AA_k, \Psi_k),(\AA_-, \Psi_-)\}
\]
with compatible boundary asymptotic limits in the sense
of section 6.2 in \cite{MMR}. Here $(\AA_\pm, \Psi_\pm)$ are solutions to 
the perturbed Seiberg-Witten equations (\ref{SW:cyl})
on $(X_\pm(\infty), \s_\pm)$, and $(\AA_i, \Psi_i)'s$ are solutions
to the  perturbed Seiberg-Witten equations (\ref{SW:4dpert}) 
on $(\R\times Y, \t)$.
\end{The}
\begin{proof} When restricted to $[-R, R]\times Y \subset X(R)$, 
$(\AA_R, \Psi_R)$ is gauge equivalent to a path $\gamma_R (t) = (A(t), \psi(t))$
$(t\in [-R, R]$ which is a gradient flowline of ${\cal C}_{\eta, \p}$ on 
$\A_Y$ with uniformly bounded energy.
Introduce the following function on $[-R, R]$
\[
f_R(t) = \|\nabla {\cal C}_{\eta, \p}(\gamma_R (t)\|^2_{L^2(Y)}
\]
Then $\int_{-R}^{R} f_R(t) dt$ is uniformly bounded.
The argument in \cite{KM1} can be employed here to show that
as $R\to \infty$, there is at least one point $t_0\in \R$ such
that the geometric limit of $f_R$ vanishes at $t_0$, hence the 
geometric limit of $(\AA_R, \Psi_R)$ splits into two
solutions to the perturbed Seiberg-Witten equations (\ref{SW:cyl})
on $(X_\pm(\infty), \s_\pm)$. Following the procedures in the proof
of Theorem 4.1 in \cite{MW}, there exists a subsequence of 
$\{(\AA_R, \Psi_R)\}_R$ which converges in $L^2_{2,loc}$-topology to
\[
\{(\AA_+, \Psi_+), (\AA_1, \Psi_1), (\AA_2, \Psi_2),
\cdots (\AA_k, \Psi_k),(\AA_-, \Psi_-)\}
\]
with the properties as claimed in the theorem.
The upper bound on $k$ follows from the uniformly bounded energy for 
$\{(\AA_R, \Psi_R)\}_R$.
\end{proof}

The standard gluing argument shows that any sequence of 
solutions with compatible asymptotic limits can be glued to
produce a solution to the perturbed Seiberg-Witten equations 
(\ref{SW:4dX}) on $(X(R), \s)$ for a sufficiently large $R$ and
some $\s \in \spinc (X, \s_\pm)$. In order to formulate  this
gluing theorem, we briefly discuss various gluing models by
studying the asymptotic limit maps for the moduli spaces
on $(X_\pm(\infty), \s_\pm)$. 

The moduli spaces $\M_{X_\pm}(\s_\pm)$ are described
in Proposition \ref{smooth:cyl} and Proposition \ref{compact:cyl}, in 
particular, there are two asymptotic limit maps:
\[
\partial_\infty^\pm: \qquad \M_{X_\pm}(\s_\pm) \longrightarrow
\M_{Y, X_\pm}(\t, \eta)
\]
where $\M_{Y, X_\pm}(\t, \eta)$ are covering spaces of $\M_Y (\t, \eta)$
with fibers consisting of $H^1(Y, \Z)/Im(i^*_\pm)$-homogeneous spaces.
Recall that $\M_{Y, X_\pm}(\t, \eta)$ and $\M_{Y, X}$ are the moduli spaces of 
3-dimensional monopoles on $(Y, \t)$ modulo those
gauge transformations which can be extended to $X_\pm$ and $X$ respectively.

Denote by $\M_{Y, +/-}$ the equivalent classes of 3-dimensional 
monopoles on $(Y, \t)$ under the action of gauge transformations on
$det (\t)$ which can be extended to either $X_+$ or $X_-$ 
(not necessarily both). Then 
$\M_{Y, +/-}$ is a covering space of $\M_Y (\t, \eta)$
with fiber a $H^1(Y, \Z)/(Im(i^*_+) + Im(i^*_-))$
homogeneous space. Hence, there is an action of
$H^1(Y, \Z)/(Im(i^*_+) + Im(i^*_-))$ on $\M_{Y, +/-}$, this
action defines an action  of
$H^1(Y, \Z)/(Im(i^*_+) + Im(i^*_-))$ on
$HF^{SW}_{*, [Im(i^*_+) + Im(i^*_-)]}(Y, \eta)$, that is, 
for any $[u]\in H^1(Y, \Z)/(Im(i^*_+) + Im(i^*_-))$, there
is an $\AA(Y)$-equivariant homomorphism on 
$HF^{SW}_{*, [Im(i^*_+) + Im(i^*_-)]}(Y, \eta)$ which decreases 
degree in $HF^{SW}_{*, [Im(i^*_+) + Im(i^*_-)]}(Y, \eta)$ by
$<[u]\wedge c_1(\eta), [Y]>$.

The following  commutative diagram 
illustrates the relationship between these covering spaces of
$\M_Y (\t, \eta)$ (the structure groups of covering spaces
over $\M_Y (\t, \eta)$
are given by Diagram (\ref{diag:subgroup}) in the Introduction):
\ba
\begin{array}{ccccc}
&& M_{Y,X_+}(\t, \eta)&& \\[-2mm]
& \nearrow&& \stackrel{\quad\pi_+}{\searrow}&
\\[-2mm]
 M_{Y, X}(\t, \eta) && \longrightarrow&&
\M_{Y, +/-}(\t, \eta)\\[-2mm]
& \searrow&& \stackrel{\pi_-\quad}{\nearrow}&\\[-2mm]
&&\M_{Y, X_-}(\t, \eta)&&
\end{array}\label{diag:covering}.
\na 

The compositions of $\pi_\pm$ with $\partial_\infty^\pm$ give rise to the 
following asymptotic limit maps:
\[
\pi_\pm\circ \partial_\infty^\pm: \qquad 
\M_{X_\pm}(\s_\pm) \longrightarrow \M_{Y, +/-}(\t, \eta).
\]
Denote by $\M_{X_\pm, X}(\s_\pm)$ the moduli spaces of solutions to
the perturbed Seiberg-Witten equations (\ref{SW:cyl}) on
$(X_\pm(\infty), \s_\pm)$ modulo those gauge transformations
on $(X_\pm(\infty), \s_\pm)$ which can be
extended to $X$.

The following gluing theorem establishes that for a sufficiently large $R$, 
the moduli space $\M_{X(R)}(\s)$ for any $\s\in \spinc (X, \s_\pm)$
can be obtained by gluing certain components of the
moduli spaces $\M_{X_\pm}(\s_\pm)$ on the cylindrical-end 4-manifolds
$(X_\pm (\infty), \s_\pm)$.

\begin{The}\label{gluing:X} For any $\spinc$ structure 
$\s= s_+ \#_{u} s_- \in  \spinc (X, \s_\pm)$ with
$u$ represents a gauge class in $ H^1 (Y, \Z)/(Im(i^*_+) + Im(i^*_-))$, 
there is a sufficiently large $R_0$ such that for any $R>R_0$, there exists the
following orientation preserving diffeomorphism
\[ \M_{X(R)}(\s)\cong
[u]\bigl(\M_{X_+}(\s_+)\bigr) \times_{\M_{Y, +/-}(\t, \eta)} \M_{X_-}(\s_-).
\]
Here $[u]\bigl(\M_{X_+}(\s_+)\bigr) 
\times_{\M_{Y, +/-}(\t, \eta)} \M_{X_-}(\s_-)$
is the fiber product of $\M_{X_+}(\s_+)$ and $\M_{X_-}(\s_-)$
over $\M_{Y, +/-}(\t, \eta)$ with resepct to the 
maps $[u]\circ \pi_+ \circ \partial^+_\infty$ and
$\pi_-\circ \partial^-_\infty$ respectively.
\end{The}
\begin{proof} We only give the proof for $[u]=0$
in $H^1 (Y, \Z)/(Im(i^*_+) + Im(i^*_-))$ here.  Then 
$\s =\s_+ \#\s_-$ is the $\spinc$ structure obtained by gluing
the $\spinc$ structures $\s_\pm$ long $Y$ using a gauge transformation 
which extends to either $X_+$ or $X_-$. The proof 
for general $[u]$ is similar. 

Given Theorem \ref{geo:limit} and
Proposition \ref{compact:cyl}, we only need to prove that for a 
sufficiently large $R$, there exists  an embedding 
\[
\M_{X_+}(\s_+) \times_{\M_{Y, +/-}(\t, \eta)} \M_{X_-}(\s_-)
\longrightarrow \bigcup_{\s\in\spinc (X, \s_\pm)} \M_{X(R)}(\s).
\]
The actual gluing procedure is fairly standard, see Theorem 4.9 in \cite{MW}
(Cf. \cite{Taubes}  and \cite{Fuk} in the instanton case) for the detailed
gluing argument. Here we just set the stage in our case and leave the details
to the reader.

Let $([\AA_+, \Psi_+], [\AA_-, \Psi_-]) \in 
\M_{X_+}(\s_+) \times_{\M_{Y, +/-}(\t, \eta)} \M_{X_-}(\s_-)$, then we can
choose two representatives, still denoted by $(\AA_+, \Psi_+)$
and $(\AA_-, \Psi_-)$ respectively. Then  from Proposition \ref{decay:cyl}
and Proposition \ref{smooth:cyl}, we know that
$(\AA_\pm , \Psi_\pm)$ are solutions to the perturbed
Seiberg-Witten equations (\ref{SW:cyl}) on $(X_\pm (\infty), \s_\pm)$,
decaying exponentially in the $C^2$-topology to
$(A_\pm, \psi_\pm)$ (solutions to the perturbed Seiberg-Witten equations
(\ref{SW:3d}) on $(Y, \t)$) and that the cokernels 
of the operators ${\cal D}_{(\AA_\pm, \Psi_\pm)}$ (\ref{linearization:cyl})
are trivial.  

As 
$$\pi_+\circ \partial_\infty^+ ([\AA_+, \Psi_+]) 
= \pi_-\circ \partial_\infty^-([\AA_-, \Psi_-])$$ there exists a gauge
transformation in $C^\infty(Y, U(1))$ with
$[g]\in Im(i^*_+)+Im(i^*_-)$ such that
$(A_+, \psi_+) = g(A_-, \psi_-)$. 
As $[g]\in Im(i^*_+)+Im(i^*_-)$ then
$g$ can be extended to either $X_+$ or $X_-$.  

Suppose that $[g]\in Im(i^*_+)$. Let
$g$ also represent the extended gauge transformation on $X_+(\infty)$. Apply
the gluing argument to the element
$(g^{-1}(\AA_+, \Psi_+), (\AA_-, \Psi_-))$ as in section 4.2 of
\cite{MW}. Then for each sufficiently large $R$, there is a unique solution
to the perturbed Seiberg-Witten equations (\ref{SW:4dX}) on 
$X(R)$ with some $\spinc$ structure in $\spinc (X, \s_\pm)$.
Write $g^{-1}(\AA_+, \Psi_+)\#_R (\AA_-, \Psi_-) $ for the resulting
monopole on $X(R)$. 
If $[g]\in Im(i^*_-)$, then the same procedure yields a unique
monopole $(\AA_+, \Psi_+)\#_R g(\AA_-, \Psi_-) $. Notice that
\[
g^{-1}(\AA_+, \Psi_+)\#_R (\AA_-, \Psi_-) = (\AA_+, \Psi_+)\#_R g(\AA_-, \Psi_-)
\]
if and only if $[g]\in Im(i^*_+)\cap Im(i^*_-)$. 
This completes the sketch
of the proof.
\end{proof}

Notice that all the  groups $Im(i^*_\pm)$, $Im(i^*_+)\cap Im(i^*_-)$
and $Im(i^*_+)+ Im(i^*_-)$  are subgroups of $Ker (c_1(\t))$, hence, there
is an induced commutative diagram of $\AA(Y)$-equivariant homomorphisms on the
corresponding Seiberg-Witten-Floer homologies:
\ba
\begin{array}{ccccc}
&& HF^{SW}_{*, [Im(i_+^*)]}(Y, \t)&& \\[-2mm]
&\nearrow&&\stackrel{\quad\pi_+}{\searrow}&\\[-2mm]
HF^{SW}_{*, [Im(i_+^*)\cap Im(i_-^*)]}(Y, \t)&&&&
HF^{SW}_{*, [Im(i_+^*)+Im(i_-^*)]}(Y, \t)\\[-2mm]
&\searrow&&\stackrel{\pi_-\quad}{\nearrow} &\\[-2mm]
&&HF^{SW}_{*, [Im(i_-^*)]}(Y, \t).&& \end{array}
\label{diag:SWF}
\na
Recall that the relative invariants $SW_{X_\pm}(\s_\pm, \cdot)$
take values in $HF^{SW}_{*, [Im(i_+^*)]}(Y, \t)$
and $HF^{SW}_{*, [Im(i_-^*)]}(-Y, -\t)$ respectively. Under the 
homomorphism $\pi_\pm$ in (\ref{diag:SWF}), $\pi_\pm(SW_{X_\pm}(\s_\pm, \cdot))$
take values in 
\[
HF^{SW}_{*, [Im(i_+^*)+Im(i_-^*)]}(Y, \t)\qquad
\text{and} \qquad
HF^{SW}_{*, [Im(i_+^*)+Im(i_-^*)]}(-Y, -\t)\]
respectively.

{}From the definition of the relative invariants $SW_{X_\pm}(s_\pm, \cdot)$
and the definition of the  Seiberg-Witten invariant for $(X(R), \s)$
with $\s= \s_+\#_{[u]}\s_-\in \spinc (X, \s_\pm)$, 
Theorem \ref{gluing:X} implies the following gluing formula.

\begin{The}\label{SW:glue}
Suppose that $X$ is a closed 4-manifold with
$b_2^+ \geq 1$,  which splits along a closed 3-manifold $(Y, \t)$
with $b_1(Y) >0$  and $c_1(\t)$ is non-torsion. 
Suppose that there are  $\spinc$ structures
$\s_+$ and $\s_-$ on $X_+$ and $X_-$ respectively such that $\s_+|_Y=
\s_-|_Y =\t$.
Then for a $\spinc$ structure $\s = \s_+\#_{[u]}\s_-$
in $\spinc (X, \s_\pm)$, where $[u]\in H^1(Y, \Z)/(Im (i_+^*)+ Im (i_-^*))$,  
we have the following gluing formula for the Seiberg-Witten invariants:
\ba
SW_X(\s, z_+z_-) = 
 \la [u]\bigl(\pi_+(SW_{X_+}(\s_+, z_+))\bigr),\pi_-(SW_{X_-}(\s_-, z_-))
 \ra 
\label{gluing:formula}
\na
where $[u]$ acts on $HF^{SW}_{_*, [Im(i_+^*)+Im(i_-^*)]}(-Y, -\t)$
as in Remark \ref{u-action}, 
$\pi_\pm$ are the induced $\AA(Y)$-equivariant homomorphisms
from $Im (i_\pm^*) \subseteq Im (i_+^*)+ Im (i_-^*) $ (see
Diagram (\ref{diag:SWF})), and $\la, \ra$ is the pairing 
on
\[
HF^{SW}_{*,[Im(i_+^*)+Im(i_-^*)]}(Y, \t)  \times
HF^{SW}_{_*, [Im(i_+^*)+Im(i_-^*)]}(-Y, -\t),
\] with the degrees in $HF^{SW}_{_*, [Im(i_+^*)+Im(i_-^*)]}(-Y, -\t)$
shifted by 
\[d_X(\s) =  \frac 14 (c_1(\s)^2 - (2 \chi (X) + 3\sigma (X) ))
= deg (z_1) + deg(z_2),
\]
and $z_\pm \in \AA(X_\pm)$. When $b_2^+ =1$, the Seiberg-Witten invariants
$SW_X (\s, \cdot)$ is defined with a fixed
orientation on $H^{2, +}(X, \R)$ such that 
$c_1(\s)\cdot \omega^+ >0$ for an oriented
generator $\omega^+$ of $H^{2, +}(X, \R)$.
\end{The}

Note that $SW_X(\s, z_1z_2) =0$ for any $\s$ with 
$d_X(\s) \neq deg (z_1) + deg(z_2)$, and 
\[
{\cal S} 
= \{\s \in \spinc (X, \s_\pm):
d_X(\s) = deg (z_1) + deg(z_2)\}
\]
 is an affine space over 
$\displaystyle{\frac {Ker (c_1(\t))}{Im(i_+^*)+Im(i_-^*)}}$. 
Summing the gluing formulae (\ref{gluing:formula}) over
all $\spinc$ structures in $\spinc (X, \s_\pm)$, we obtain
that
\ba\label{gluing:formula:sum}
\sum_{\s\in {\cal S}}
SW_X(\s, z_+z_-) =
 \la (\pi (SW_{X_+}(\s_+, z_+)),\pi (SW_{X_-}(\s_-, z_-))
 \ra.
\na
Here $\pi$ is the map (Cf. (\ref{equ:homo}) )
from $HF^{SW}_{*, [Im(i^*_\pm)]}(\pm Y, \pm \eta)$
to \[
HF^{SW}_{*, [Ker (c_1(\eta))]}(\pm Y, \pm \eta) 
\cong HF^{SW}_*(\pm Y, \pm \eta),\]
under the periodicity map (\ref{isomorphism}),  
$\pi (SW_{X_\pm}(\s_\pm, z_\pm))$  take values
in $HF^{SW}_*(\pm Y, \pm\t)$, and the right hand side of 
(\ref{gluing:formula:sum}) the natural pairing on
\[
HF^{SW}_*(Y, \t) \times HF^{SW}_*(-Y, -\t).\]
In particular, in the case of $Im(i_+^*)+Im(i_-^*) = Ker (c_1(\t))$, the
left hand side of (\ref{gluing:formula:sum}) has at most one term.

As an application of these gluing formulae, we briefly mention 
some results in \cite{Munoz}. Applying the gluing theorem, Mu\~noz
and Wang obtained a ring structure on 
$HF_{*}^{SW}(\Sigma_g\times S^1, \s_r)$ 
where $\Sigma_g$ is a closed surface of genus $g\ge 1$ and $\s_r$
is the $\spinc$ structure on $\Sigma_g \times S^1$
with $c_1(\s_r) = 2r PD(S^1)$ ($-(g-1) \le r \le g-1$). Note
that $Ker (c_1(\s_r)) = H^1(\Sigma_g, \Z)$ and 
$H^1(\Sigma_g \times S^1, \Z) /Ker (c_1(\s_r))\cong \Z$.  Set 
$d= g-1 -|r|$, we only discuss the case of  $r\ne 0$ here.
Denote by  $\s_r$ also the $\spinc$ structure on $\Sigma_g \times D^2$
with $c_1(\s_r) = 2r PD ([pt \times D^2])$,
where $D^2$ is a 2-dimensional disc. The relative invariants
for $\Sigma_g \times D^2$ define a map
\[
\begin{array}{c}
\AA(\Sigma_g \times D^2) = \AA (\Sigma_g) \longrightarrow
HF^{SW}_{*}(\Sigma_g \times S^1, \s_r)\\[2mm]
z \mapsto SW_{\Sigma_g\times D^2}(\s_r, z).
\end{array}
\]
Then the gluing theorem (Theorem \ref{SW:glue}) tells us that
\[
\la  SW_{\Sigma_g\times D^2}(\s_r, z_1),  SW_{\Sigma_g\times D^2}(\s_r, z_2)
\ra = \sum_{n \in \Z}  
SW_{\Sigma_g\times S^2}(\s_r + n PD (\Sigma_g), z_1 z_2),
\]
for any $z_1, z_2 \in \AA(\Sigma_g)$. For dimensional reasons 
there is at most one $n \in \Z$ such that 
$SW_{\Sigma_g \times S^2}(\s_r + 2n PD ([\Sigma_g]), z_1z_2)$ is non-zero. 
If $deg (z_1) + deg(z_2) = 2d$,  then
\[
\la  SW_{\Sigma_g\times D^2}(\s_r, z_1),  SW_{\Sigma_g\times D^2}(\s_r, z_2)
\ra = \la z_1 z_2, [Sym^{d}(\Sigma_g)]\ra.
\]
These pairing, along with others in \cite{Munoz}, can be
used to study the structure of $HF_{*}^{SW}(\Sigma_g \times S^1,
\s_r)$, see \cite{Munoz} for explicit calculations.

\vskip .2in

\noindent {\bf A.L. Carey},

\noindent School of Mathematical Sciences,
Australian National University, Canberra ACT, Australia\par
\noindent acarey@maths.anu.edu.au\par

\vskip .2in
\noindent {\bf B.L. Wang},

\noindent Department of Pure
Mathematics, University of Adelaide, Adelaide SA 5005 \par
\noindent bwang@maths.adelaide.edu.au \par
\noindent Max-Planck-Institut f\"ur Mathematik, D-53111 Bonn, Germany
\par
\noindent bwang@mpim-bonn.mpg.de

\end{document}